\documentclass[preprint,12pt]{elsarticle}
\usepackage{xcolor}
\usepackage{amsmath}
\usepackage{amsthm}
\usepackage{amssymb}
\usepackage{algorithm}
\usepackage{algpseudocode}
\usepackage{setspace}
\usepackage{float}
\usepackage{booktabs}
\usepackage{makecell}
\usepackage{caption}
\usepackage{subcaption}
\usepackage{adjustbox} 
\usepackage{rotating}
\usepackage{enumitem}

\usepackage{listings}
\usepackage[framed,numbered,autolinebreaks,useliterate]{mcode}
\journal{Computer Methods in Applied Mechanics and Engineering}
\theoremstyle{remark}
\newtheorem{rem}{Remark}
\theoremstyle{scheme}
\newtheorem{scheme}{Scheme}
\usepackage{comment} %
\begin{document}

\begin{frontmatter}



\title{A GPU-Accelerated Matrix-Free FAS Multigrid Solver for Navier–Stokes Equations with Memory-Efficient Implementations}

\author[inst1]{Jiale Meng}
%\ead{mengjl@csrc.ac.cn}

\author[inst1]{Shuqi Tang}
%\ead{tangsq@csrc.ac.cn}

\author[inst2]{Steven M. Wise}
%\ead{swise1@utk.edu}

\author[inst1]{Zhenlin Guo\corref{cor1}}
\ead{zguo@csrc.ac.cn}
\cortext[cor1]{Corresponding author}

\affiliation[inst1]{organization={Beijing Computational Science Research Center},
                    city={Beijing},
                    postcode={100193},
                    country={China}}

\affiliation[inst2]{organization={Department of Mathematics, The University of Tennessee},
                    city={Knoxville},
                    postcode={TN 37996},
                    country={USA}}

\begin{abstract}
We develop a matrix-free Full Approximation Storage (FAS) multigrid solver based on staggered finite differences and implemented on GPU in MATLAB. To enhance performance, intermediate variables are reused, and an X-shape Multi-Color Gauss–Seidel (X-MCGS) smoother is introduced. This scheme eliminates the conditional branching required to distinguish red and black nodes in the standard two-color Red–Black Gauss–Seidel (RBGS) method by partitioning the grid into four submatrices according to row–column parity. In addition, restriction and prolongation operators are implemented with GPU acceleration. Algebraic and asymptotic convergence tests verify the solver’s robustness and accuracy, while benchmarks on large-scale problems demonstrate substantial performance gains: for the 2D heat equation on an $8192^2$ grid, the NVIDIA RTX~4090 GPU (24~GB) achieves a 61× speedup over a single core of the Intel Core i9-13900KS CPU, and in 3D at $512^3$, it reaches 46×. To overcome GPU memory limitations, we design a memory-efficient implementation of the first- and second-order projection schemes for the Navier–Stokes equations using a dynamic reuse strategy, which reduces GPU-resident variables from 12 (1st-order scheme) and 15 (2nd-order scheme) to only 8. lowering memory footprint and improving performance by 20–30\%. This enables $512^3$ computations of the Navier–Stokes equations on a single RTX~4090, where classical implementations run out of memory. The solver’s applicability is further demonstrated through large-scale simulations. Grain growth simulations on a $512^2$ grid accommodate up to $q=1189$ orientations in 2D and $q=123$ in 3D, with fitted growth exponents reproducing the expected scaling laws. Moreover, the memory-efficient Navier–Stokes implementations, coupled with the Cahn–Hilliard equations, enable air–water two-bubble coalescence simulations on a $256 \times 256 \times 1024$ grid using a single RTX~4090 GPU, yielding results in close agreement with experimental observations.

\end{abstract}

%Research highlights
\begin{highlights}
\item Memory-efficient GPU-parallel matrix-free FAS multigrid solver on staggered grids in MATLAB.
\item X-shape Multi-Color Gauss–Seidel smoother eliminating branching, enhancing GPU parallelism.
\item Heat equation on cell-centered grids: up to 61× (2D) and 46× (3D) speedup vs Intel i9-13900KS on RTX 4090.
\item First- and second-order Navier–Stokes schemes using 8 GPU variables, 20–30\% faster on $512^3$ grids.
\item Phase-field simulations of grain growth, phase separation, and two-bubble coalescence confirm solver accuracy.
\end{highlights}

%% Keywords

\end{frontmatter}

\section{Introduction}
\label{sec1}
%\iffalse
The geometric multigrid method is among the most efficient solving techniques for discrete algebraic systems arising from elliptic partial differential equations. It originated from the pioneering work of Fedorenko in the 1960s \cite{bakhvalov1966,fedorenko1962,fedorenko1964} and was further developed by Brandt in the 1970s \cite{brandt1977,hackbusch1977}. Today, it is widely recognized in computational mathematics (see \cite{hackbusch1984,fulton1986} for reviews). By solving discrete systems—typically obtained from finite difference or finite element discretizations—on fine meshes while exploiting a hierarchy of coarser grids, the geometric multigrid method accelerates convergence by efficiently transferring error components across different levels of resolution. To extend multigrid methods to nonlinear problems, Brandt (1977) proposed the Full Approximation Storage (FAS) scheme \cite{brandt1977,trottenberg2001}, which has since become a standard approach for inter-grid communication. Unlike correction schemes designed primarily for linear systems, FAS transfers approximate solutions directly across grid levels, thereby making multigrid methods naturally applicable to nonlinear problems.

In recent years, Graphics Processing Units (GPUs) have emerged as a powerful platform for accelerating large-scale scientific computing, owing to their massive parallelism and high memory bandwidth. Geometric multigrid solvers, with their algorithmic scalability, structured data access, and localized operations, are particularly amenable to GPU parallelization, making the combination especially effective for large-scale simulations\cite{Bolz2003,ljungkvist2017,feng2014numerical,antepara2024,cui2025,chen2024}. Among multigrid variants, the Full Approximation Storage (FAS) framework is especially well-suited, as it extends multigrid methods to nonlinear problems while preserving the recursive structure and locality of operations, thereby offering both algorithmic flexibility and GPU-friendly parallelism \cite{shi2020,BraedstrupEgholm2014,Gorobets2024}. In 2020, Shi \cite{shi2020} presented an FAS multigrid solver integrated with multi-GPU parallelization, demonstrating its effectiveness in simulating unsteady incompressible flows. Most recently, Gorobets (2024)~\cite{Gorobets2024} developed a heterogeneous parallel implementation of FAS multigrid in the Noisette code, combining MPI, OpenMP, and GPU acceleration to handle large-scale compressible flow simulations. Their approach efficiently exploits hybrid cluster architectures, significantly improving parallel efficiency and highlighting the scalability of FAS multigrid solvers on modern GPU systems.

Although GPUs provide substantial computational power, their limited onboard memory often becomes a major bottleneck for large-scale matrix computations \cite{antepara2024,Ansari2025,Recasens2025,griebel2010}. A widely adopted solution is to employ multiple GPUs to expand the available memory pool \cite{feng2014numerical,shi2020,Kashi2025}. However, to fully exploit the performance of a single GPU, it remains essential to develop memory-efficient strategies. One promising direction is the use of matrix-free methods within geometric multigrid frameworks \cite{ljungkvist2017,cui2025,chen2024}, which avoid the explicit storage of large sparse matrices while leveraging the structured data access of geometric multigrid frameworks. This approach not only reduces memory consumption but also improves data locality—both well aligned with GPU parallel execution. For instance, Chen et al.~\cite{chen2024} proposed a matrix-free geometric multigrid preconditioned conjugate gradient (MGCG) solver based on the summation-by-parts simultaneous-approximation-term (SBP–SAT) method, incorporating specialized GPU kernels. Compared with conventional sparse matrix–vector multiplication (SpMV), their implementation not only reduced the memory footprint but also delivered substantial computational acceleration. In large-scale problems involving approximately 67 million degrees of freedom, the matrix-free GPU MGCG achieved an overall 5× speedup, underscoring the effectiveness of matrix-free geometric multigrid strategies for stable and efficient large-scale PDE simulations on modern GPU architectures.

In matrix-free frameworks, iterative methods are generally preferred for solving large-scale linear systems. The classical Jacobi iteration is straightforward to parallelize but suffers from slow convergence \cite{saad2003}. The Gauss–Seidel iteration, as an improved variant, achieves faster convergence with fewer iterations due to its higher efficiency \cite{feng2014numerical,mazumder2016}, yet its inherently sequential nature makes parallel implementation challenging. To address this issue, the Red–Black Gauss–Seidel (RBGS) method introduces a two-color partitioning of the computational grid, enabling parallel updates while largely preserving the convergence properties of the original scheme \cite{wesseling1992,briggs2000}. Owing to its simplicity, RBGS has been widely adopted in CPU- and MPI-based parallel solvers. On GPU architectures, however, RBGS typically requires conditional branching to distinguish between red and black nodes, which introduces thread divergence, hinders vectorization, and thus significantly undermines computational efficiency. To overcome the drawbacks of RBGS on GPUs, multi-color Gauss–Seidel (MCGS) methods~\cite{feng2014numerical,ZHAO2022,li2020} have been developed. Instead of a two-color partition, the computational grid is divided into four disjoint sets in two dimensions or eight sets in three dimensions, based on the parity of the grid indices. In \cite{feng2014numerical}, geometric multigrid methods on CPU–GPU heterogeneous computers are studied, and comparisons of weighted Jacobi, two-color Gauss–Seidel, and four-color Gauss–Seidel smoothers in 2D show that the four-color scheme yields superior convergence properties on GPUs. In \cite{ZHAO2022}, a Constrained Pressure Residual (CPR)-preconditioned algebraic multigrid method with a four-color Gauss–Seidel smoother using a Z-shape ordering sequence is proposed to efficiently solve the fully implicit black oil model and improve parallel performance on both CPU and GPU architectures.     
%???cite other multi color papers???

Building upon the matrix-free paradigm and modern GPU architectures,this study proposes a matrix-free Full Approximation Storage (FAS) multigrid solver based on staggered finite differences and implemented in MATLAB with GPU acceleration. The solver employs an X-shape multi-color Gauss–Seidel (X-MCGS) smoother which, unlike the traditional two-color Red–Black GS that suffers from thread divergence on GPUs, partitions the grid into four (2D) or eight (3D) disjoint sets following an X-shaped ordering.This design removes conditional branching for color checks, promotes regular memory access, and enables fully vectorized updates. Although each X-MCGS sweep involves more substeps than RBGS, the improved memory regularity and parallel efficiency typically offset the additional work, making X-MCGS an effective smoother for GPU-accelerated multigrid solvers. Numerical experiments further indicate that the X-shape ordering requires only about half as many iterations as the U-shape variant and roughly one-third as many as the Z-shape variant to achieve comparable convergence tolerances.

To further enhance GPU performance, we introduce memory-efficient strategies that reuse intermediate variables and merge computational steps, thereby reducing the number of temporary arrays during multigrid operations. In time-dependent simulations, such as those with Crank–Nicolson schemes, solution variables at the current step overwrite those from the previous step, further lowering memory usage. For incompressible Navier–Stokes equations, we design memory-efficient implementations of both first- and second-order projection schemes using a dynamic reuse strategy, which reduces GPU-resident variables from 12 (first-order) and 15 (second-order) to only eight. This optimization lowers the memory footprint and improves execution efficiency by 20–30\%, enabling large-scale simulations up to $512^3$ on a single NVIDIA RTX~4090 GPU (24~GB), whereas classical implementations would exceed device memory.

Motivated by these developments, we apply our matrix-free FAS multigrid solver to a variety of large-scale multiphase flow and material problems, including grain growth \cite{fan1997} on a $512^2$ grid with a maximum of 1189 orientation variables in 2D and a $128^3$ grid with a maximum of 123 orientation variables in 3D, phase separation~\cite{wang2022lipid,ratz_pdes_2006,teigen2009diffuse, li2009solving} on a $256^3$ grid, and two-bubble coalescence in water–air systems \cite{guo2022, denefle2014,brereton1991} on a $256 \times 256 \times 1024$ grid, all computed using a single RTX~4090 GPU. The results show excellent agreement with experimental data, demonstrating strong scalability across different physical models while maintaining a compact and accessible MATLAB implementation.

The paper is organized as follows. In \S~\ref{sec2}, we present the algorithmic framework of the GPU-accelerated, matrix-free FAS multigrid solver on cell-centered grids, including the design of smoothing, restriction, and prolongation operators. To enhance GPU parallelism, we introduce the X-MCGS smoothing operator, and compare its convergence behaviour with alternative orderings such as Z-shape and U-shape, demonstrating that the X-MCGS achieves convergence in the fewest iterations. In \S~\ref{sec3}, numerical validations are presented through convergence studies and execution performance tests, followed by two large-scale physical applications: grain growth and vesicle phase separation. \S~\ref{sec4} extends the matrix-free FAS multigrid solver to staggered grids and focuses on GPU memory-efficient implementations for solving the Navier–Stokes equations with both first-order and second-order projection schemes. Finally, \S~\ref{sec5} concludes the paper and discusses potential directions for future improvements.

\section{GPU-accelerated FAS multigrid solver}
\label{sec2}

In this section, we present a GPU-accelerated, matrix-free Full Approximation Storage (FAS) multigrid solver based on finite difference discretization on two- and three-dimensional staggered grids. To maximize GPU performance, we introduce an X-shape Multi-Color Gauss–Seidel (X-MCGS) smoothing operator—employing four-color in 2D and eight-color in 3D—which avoids the conditional branching required to distinguish red and black nodes in the standard two-color Red–Black Gauss–Seidel (RBGS) smoother. For clarity, we first illustrate the solver with two-grid in 2D as a representative example, before extending it to full multilevel implementations.

\subsection{Definitions and notations for two-grid solver}\label{sec2.1}
We first introduce the fundamental definitions and notations used in the finite difference discretization on two-dimensional uniform staggered grids. Let $\Omega =\prod_{l=1}^{2} [a_l, b_l] \subset \mathbb{R}^2$ denote the continuous computational domain, where $\mathbb{R}$ represents the real number space. The two-grid solver involves two levels of discretization: a fine grid and a coarse grid. We denote the corresponding discrete domains by $\Omega_1$ (fine grid) and $\Omega_0$ (coarse grid), respectively. Let \( M_1 = 2^{m}\) and \( N_1 = 2^{n} \) represent the numbers of grid intervals in the horizontal and vertical directions on the fine grid, where \( m \) and \( n \) are positive integers. To ensure uniform discretization, the fine grid spacings in both directions must satisfy
\begin{align}
h_1 = \frac{b_1 - a_1}{M_1}= \frac{b_2 - a_2}{N_1}.
\end{align}
For the coarse grid, let \( M_0 \) and \( N_0 \) denote the numbers of grid intervals in the horizontal and vertical directions, respectively. These are related to the corresponding fine-grid intervals by $M_0 = M_1/2$ and $N_0 = N_1/2$. It then follows directly that the coarse-grid spacing is $h_0 = 2h_1$.  
 
\begin{figure}[H]
\centering
\includegraphics[scale=0.7]{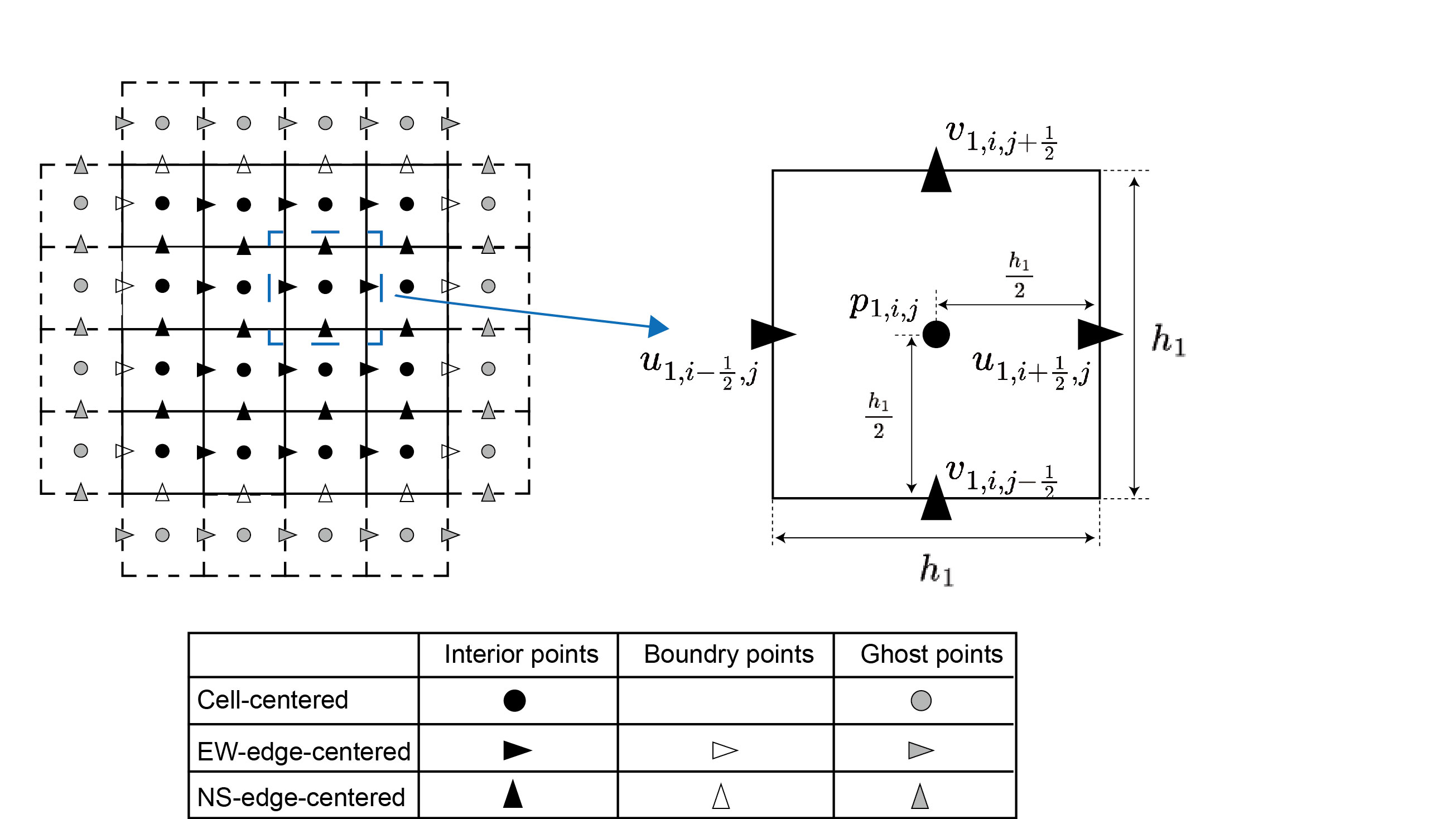}
\caption{Illustration of the fine grid $\Omega_1$ in the two-grid framework, which can be further generalized to multilevel grid hierarchies. The figure presents a 2D uniform staggered discretization with cell-centered and edge-centered variables. The left panel depicts the global grid configuration, whereas the right panel provides a magnified view of a representative control volume with mesh size $h_1$. The symbols denote the types of variables as specified in the legend. See \S~\ref{sec2.1} for details.}
\label{grid}
\end{figure}

We next introduce several finite difference index sets associated with the fine grid $\Omega_1$:
\begin{align}
E_{M_1}      &= \left\{ x_{i+\frac{1}{2}} \,\middle|\, i=1,\dots,M_1-1 \right\}, 
&\quad \bar{E}_{M_1} &= \left\{ x_{i+\frac{1}{2}} \,\middle|\, i=0,\dots,M_1 \right\}, \\
E_{N_1}      &= \left\{ y_{j+\frac{1}{2}} \,\middle|\, j=1,\dots,N_1-1 \right\},  
&\quad \bar{E}_{N_1} &= \left\{ y_{j+\frac{1}{2}} \,\middle|\, j=0,\dots,N_1 \right\}, \\
C_{M_1}      &= \left\{ x_i \,\middle|\, i=1,\dots,M_1 \right\}, 
&\quad \bar{C}_{M_1} &= \left\{ x_i \,\middle|\, i=0,\dots,M_1+1 \right\}, \\
C_{N_1}      &= \left\{ y_j \,\middle|\, j=1,\dots,N_1 \right\}, 
&\quad \bar{C}_{N_1} &= \left\{ y_j \,\middle|\, j=0,\dots,N_1+1 \right\},
\end{align}
where the grid coordinates are defined by $x_{i+1/2} = a_1 + i \cdot h_1$, $y_{j+1/2} = a_2 + j \cdot h_1 $, and $x_i = a_1 + \left(i - 1/2\right) \cdot h_1$, $y_j = a_2 + \left(j - 1/2\right) \cdot h_1$. The staggered grid consists of interlaced cell-centered and edge-centered points, as illustrated in Figure~\ref{grid}. In the $x$-direction, $E_{M_1}$ and $\bar{E}_{M_1}$ represent uniform partitions of $[a_1, b_1]$ with $M_1$ divisions. The points in $E_{M_1}$ and $\bar{E}_{M_1}$ are edge-centered points, while the additional two points in $\bar{E}_{M_1} \setminus E_{M_1}$ are boundary points of the edge-centered grid. Similarly, $C_{M_1}$ and $\bar{C}_{M_1}$ denote cell-centered points, with the extra points in $\bar{C}_{M_1} \setminus C_{M_1}$ serving as ghost points for the cell-centered grid. Analogous definitions apply in the $y$-direction. The above construction introduces a single layer of ghost points. Note that additional layers of ghost points may be introduced at both cell-centered and edge-centered locations as needed.

We further define the following function spaces with respect to the fine grid $\Omega_1$: 
\begin{align}
\mathcal{C}_{\Omega_1} &= \left\{ p_{1,i,j} : C_{M_1} \times C_{N_1} \rightarrow \mathbb{R} \right\},&
\bar{\mathcal{C}}_{\Omega_1} &= \left\{ p_{1,i,j} : \bar{C}_{M_1} \times\bar{C}_{N_1} \rightarrow \mathbb{R} \right\}, \label{cell center} \\
\mathcal{E}^{EW}_{\Omega_1} &= \left\{ u_{1,i+\frac{1}{2},j} : E_{M_1} \times C_{N_1} \rightarrow \mathbb{R} \right\}, 
&\bar{\mathcal{E}}^{EW}_{\Omega_1} &= \left\{ u_{1,i+\frac{1}{2},j} : \bar{E}_{M_1} \times \bar{C}_{N_1} \rightarrow \mathbb{R} \right\}, \label{ew edge center}\\
\mathcal{E}^{NS}_{\Omega_1} &= \left\{ v_{1,i,j+\frac{1}{2}} : C_{M_1} \times E_{N_1} \rightarrow \mathbb{R} \right\},
&\bar{\mathcal{E}}^{NS}_{\Omega_1} &= \left\{ v_{1,i,j+\frac{1}{2}}:\bar{C}_{M_1} \times \bar{E}_{N_1} \rightarrow \mathbb{R} \right\},\label{ns edge center}
\end{align}
where $\mathcal{C}_{\Omega_1}$ and $\bar{\mathcal{C}}_{\Omega_1}$ denote the spaces of cell-centered functions, $\mathcal{E}^{EW}_{\Omega_1}$ and $\bar{\mathcal{E}}^{EW}_{\Omega_1}$ denote the spaces of east-west (EW) edge-centered functions, and $\mathcal{E}^{NS}_{\Omega_1}$ and $\bar{\mathcal{E}}^{NS}_{\Omega_1}$ denote the spaces of north-south (NS) edge-centered functions. The corresponding index sets $E_{M_0}$, $\bar{E}_{M_0}$, $E_{N_0}$, $\bar{E}_{N_0}$, $C_{M_0}$, $\bar{C}_{M_1}$, $C_{N_0}$, $\bar{C}_{N_0}$, together with the function spaces $\mathcal{C}_{\Omega_0}$, $\bar{\mathcal{C}}_{\Omega_0}$, $\mathcal{E}^{EW}_{\Omega_0}$, $\bar{\mathcal{E}}^{EW}_{\Omega_0}$, $\mathcal{E}^{NS}_{\Omega_0}$, and $\bar{\mathcal{E}}^{NS}_{\Omega_0}$, are defined analogously on the coarse grid.

\subsection{Two-grid solver on 2D cell-centered grids}\label{sec2.2}

In this section, we use the cell-centered grids as a representative example to illustrate our GPU-accelerated matrix-free FAS solver with two-grid. Extensions to staggered grids are discussed in detail in \S~\ref{sec4}. Specifically, we consider the application of the two-grid solver to the following time-independent problem defined on $\Omega = [0, 1]^2$:
\begin{equation}\label{PDE}\left\{\begin{aligned}
    p-\Delta p &= f \quad\quad \text{in} \,\,\quad\Omega,\\
    p &= 0  \quad \quad\text{on} \quad\partial\Omega.
\end{aligned}\right.\end{equation}

We begin by discretizing Equation~\eqref{PDE} on the fine grid $\Omega_1$. 
Specifically, we prescribe the source term \( f_{1,i,j} \in \mathcal{C}_{\Omega_1} \) and seek a solution 
\( p_{1,i,j} \in \bar{\mathcal{C}}_{\Omega_1} \) such that
\begin{align}
&L_{h_1}(p_{1,i,j}) = p_{1,i,j} - \Delta_{h_1} p_{1,i,j} = f_{1,i,j},  \label{PDE_discret}\\
&p_{1,0,j} = -p_{1,1,j},\quad p_{1,N+1,j} = p_{1,N,j},  \label{boundry1}\\
&p_{1,i,0} = -p_{1,i,1},\quad p_{1,i,N+1} = p_{1,i,N},  \label{boundry2}
\end{align}
where $L_{h_1}(p_{1,i,j}) \in \mathcal{C}_{\Omega_1}$ denotes the discrete operator applied 
to $p_{1,i,j}$ on the fine grid, and $\Delta_{h_1}$ is the 5-point Laplacian mapping from 
the extended function space $\bar{\mathcal{C}}_{\Omega_1}$ to $\mathcal{C}_{\Omega_1}$. 
The corresponding point-wise residual is defined by
\begin{align}
r_{1,i,j} = f_{1,i,j} - L_{h_1}(p_{1,i,j}),
\end{align}
and its discrete $L^2$ norm, which is used as the residual tolerance, is given by
\begin{align}
\text{res} = \| r_{1,i,j} \|_{L^2(\Omega_1)} 
= h_1 \left( \sum_{i,j} |r_{1,i,j}|^2 \right)^{1/2}.
\label{residual}
\end{align}

The corresponding finite-difference discretization on the coarse grid $\Omega_0$ takes the same form as Equation~\eqref{PDE_discret}, except that the source term is replaced by \( f_{0,i,j} \), which is obtained by restricting the fine-grid residual \( r_{1,i,j} \) and substituting the restricted solution \( p_{1,i,j} \) into the coarse-grid operator $L_{h_0}(p_{0,i,j}) \in \mathcal{C}_{\Omega_0}$. All other components remain analogous to those on the fine grid and are omitted here for brevity.

We now present our matrix-free two-grid solver in Algorithm~\ref{two-grid algorithm}. For notational simplicity, we omit the discrete indices \(i, j\), and denote the solution variable simply as \(p_{\text{level}}\) rather than \(p_{\text{level}, i, j}\) on the corresponding grid level. The same convention applies to all other variables.

\begin{algorithm}[H]
\begin{spacing}{1.5} 
\caption{$p_1=\mathrm{TMGSolver}(p_1^{1,0},f_1,L_{h_1},L_{h_0},k_{\text{Max}},s,tol)$}
\label{two-grid algorithm}
\begin{algorithmic}[1]
\State \textbf{input:} initial guess $p_1^{1,0}$, source term $f_1$ 
%\Comment{Define source term and initial solution}
%\Comment{Initialize iteration counter and residual}
%\While{$res>tol$ \text{or} $k\leq k_{\text{Max}}$ }     
\For {$k=1:k_{\text{Max}}$}
    \For{$m=1:s$}
        \State  $p_1^{k,m} := \mathrm{S}(f_1, p_1^{k,m-1})$ \Comment{Fine grid pre-smoothing }
    \EndFor
        \State $r_1 := f_1 - L_{h_1}(p_1^{k,m})$\Comment{Fine grid point-wise residual}
        \State $p_0^{k,0} := \mathrm{R}(p_1^{k,m}), \quad r_0 := \mathrm{R}(r_1)$ \Comment{Coarse grid approximation}
        \State $f_0 := r_0 + L_{h_0}(p_0^{k,0})$ \Comment{Coarse source term}
     \For{$m=1:s$}
         \State $p_0^{k,m} := \mathrm{S}(f_0, p_0^{k,m-1})$\Comment{Coarse grid smoothing }
    \EndFor 
        \State $c_0 := p_0^{k,m} - p_0^{k,0}$\Comment{Coarse grid correction}
        \State $p_1^{k,0} := p_1^{k,m} + \mathrm{P}(c_0)$\Comment{Fine grid correction}
    \For{$m=1:s$}
      \State $p_1^{k,m} := \mathrm{S}(f_1, p_1^{k,m-1})$\Comment{Fine grid post-smoothing }
    \EndFor  
        \State$\text{res}:= \|f_1 - L_{h_1}(p_1^{k,m})\|_{L^2(\Omega_1)}$ \Comment{Fine grid residual}
    \If{$\text{res}> tol$}
         \State $p_1^{k+1,0} := p_1^{k,m}$ \Comment{Update for $(k+1)$-th iteration}
    \Else
         \State \textbf{output:} final approximation $p_1 = p_1^{k,m}$
         \State \textbf{break} \Comment{Exit the for-loop}
    \EndIf
%\EndWhile
\EndFor

\end{algorithmic}
\end{spacing}
\end{algorithm}

For clarity, the symbols used in our matrix-free FAS two-grid solver are summarized in detail in Table~\ref{tab:FASsymbols}.
\begin{table}[H]
\centering
\rule{\linewidth}{0.4pt} % caption 下方横线
\begin{tabular}{c @{\hskip 3.5em} l}
\\
\textbf{Built-in functions and parameters} & \textbf{Description} \\
$\Omega_0, \Omega_1$ & Coarse and fine grids. \\
$p_0^{k,m}, p_1^{k,m}$ & Approximate solutions on coarse and fine grids. \\
$f_0, f_1$ & Source terms on coarse and fine grids. \\
$L_{h_0}, L_{h_1}$ & Discretized operators on coarse and fine grids. \\
$r_0, r_1$ & Point-wise residuals on coarse and fine grids. \\
$\text{res}$ & Residual after completion of the multigrid V-cycle. \\
$c_0$ & Coarse grid correction. \\
$\mathrm{S}$ & Smoothing operator. \\
$\mathrm{R}$ & Restriction operator. \\
$\mathrm{P}$ & Prolongation operator. \\
$k$ & Iteration counter. \\
\\
\textbf{User-defined parameters} &  \\
$s$ & Number of smoothing steps. \\
$k_{\text{Max}}$ & Maximum number of iterations. \\
$tol$ & Convergence tolerance. \\
$meshLevel$ & Depth of recursion in a full multigrid solver. \\
\end{tabular}
\caption{Definition of symbols employed in Algorithm~\ref{two-grid algorithm}}
\label{tab:FASsymbols}
\rule{\linewidth}{0.4pt} % caption 下方横线
\end{table}

To reduce memory overhead on the GPU, the solver reuses intermediate variables across computational steps. For example, in Step \textbf{9}, the operation $f_0 = L_{h_0}(p^0_0) + r_0$ is reformulated as $r_0 = L_{h_0}(p^0_0) + r_0$, thereby avoiding additional storage requirements. The smoothing operator $\mathrm{S}$ will be described in detail in the following subsection, while the restriction operator $\mathrm{R}$ and prolongation operator $\mathrm{P}$—which link the cell-centered coarse and fine grids in Algorithm~\ref{two-grid algorithm}—are presented in \S\ref{app3} together with their explicit formulations. It is worth emphasizing that both $\mathrm{R}$ and $\mathrm{P}$ can be expressed in matrix form, making them highly amenable to GPU-based parallelization and enabling substantial improvements in computational efficiency and execution time performance.

\begin{rem}
In conventional two-grid solvers, the coarse-grid level (Steps \textbf{9-11}) is usually solved by a direct method \cite{FENG2018}. In contrast, to make the matrix-free framework better suited for GPU execution, we replace the direct solve with X-shaped multi-color smoothing.
\end{rem}

\begin{rem}
A full matrix-free multigrid solver is obtained by recursively applying the two-grid solver, where Steps \textbf{9-11} are replaced with a recursive call on the next coarser grid. The recursion depth is controlled by the variable $meshLevel$ (with $meshLevel = 1$ for a two-grid scheme). In the multigrid case, additional recursion levels may be introduced to more effectively eliminate low-frequency error components and improve convergence, continuing until a user-defined coarsest level is reached. 
 \end{rem}

\subsection{GPU-optimized smoothing operator $\mathrm{S}$: X-shape Multi-Color Gauss–Seidel}\label{sec2.3}

We now describe our GPU-optimized X-shape Multi-Color Gauss–Seidel smoothing operator $\mathrm{S}$, designed to fully exploit parallelism on GPU and applied as a central component of the two-grid solver in Steps~\textbf{4}, \textbf{10}, and \textbf{15} of Algorithm~\ref{two-grid algorithm}. 

In the finite-difference, matrix-free framework, iterative methods are generally preferred for solving large-scale linear systems. The classical Jacobi iteration is straightforward to parallelize but suffers from slow convergence \cite{saad2003}. The Gauss–Seidel iteration, as an improved variant, achieves faster convergence with fewer iterations due to its higher efficiency \cite{feng2014numerical,mazumder2016}. However, its inherently sequential nature makes parallel implementation challenging. To overcome the sequential nature of Gauss–Seidel iteration, the Red–Black Gauss–Seidel (RBGS) method partitions the grid into two colors, enabling parallel updates while preserving convergence properties \cite{wesseling1992,briggs2000}. Its simplicity has led to widespread use in CPU- and MPI-based solvers; however, on GPUs the need for conditional branching to distinguish red and black nodes introduces thread divergence and hinders vectorization, resulting in significant performance degradation. A matrix-free formulation of RBGS is outlined in \S\ref{app1}, with a MATLAB implementation provided in Listing~\ref{RBGS smoothing}.

To further improve GPU efficiency, we propose an X-shape Multi-Color Gauss–Seidel (X-MCGS) scheme—using four-color in 2D and eight-color in 3D—which partitions the grid into four/eight disjoint sets following an X-shaped ordering. This design eliminates conditional branching during parallel execution (see Listing~\ref{FCGS smoothing} for a 2D example), promotes more regular memory access patterns, and reduces thread divergence, thereby improving GPU utilization. Although each iteration involves four sub-steps instead of two, the enhanced parallel efficiency compensates for the extra cost, making the X-MCGS scheme particularly well-suited for large-scale GPU-based multigrid solvers.

\begin{figure}[H] 
    \centering
    \includegraphics[scale=0.5]{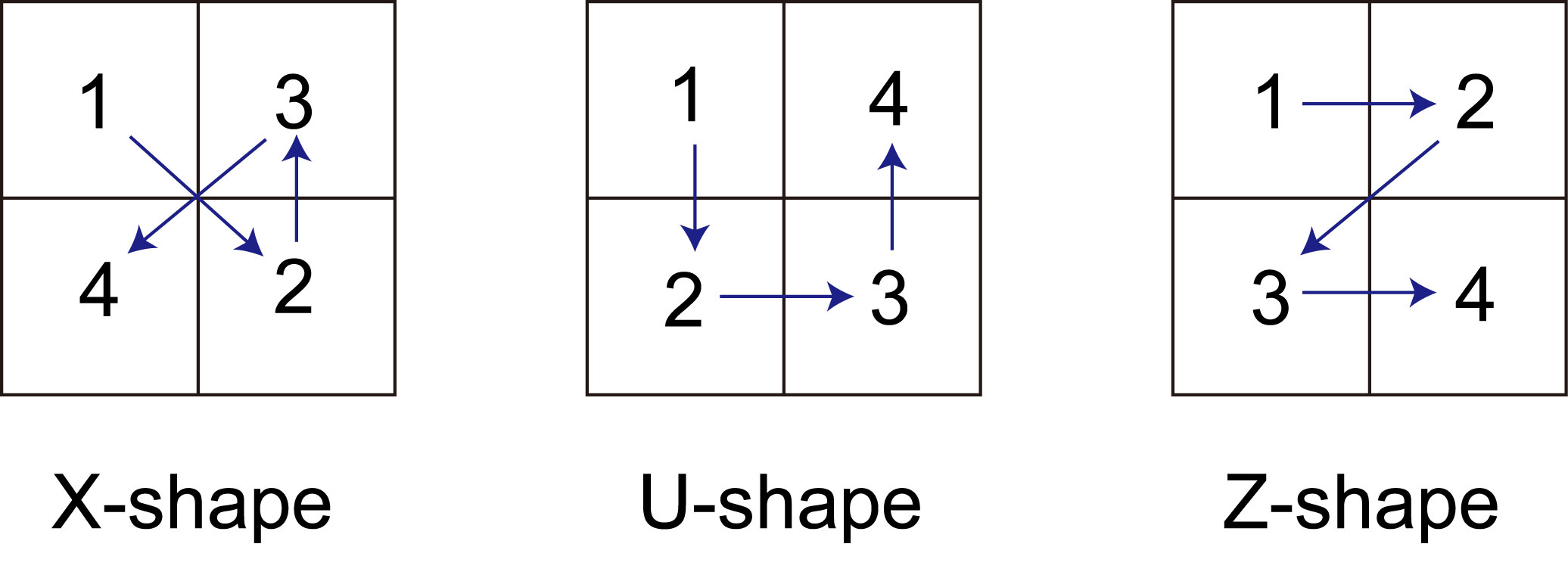}
    \caption{Illustration of three MCGS partitioning patterns used for smoothing on 2D uniform grids: (left) X-shape, (middle) U-shape, and (right) Z-shape. The numbers (1–4) indicate the sub-block update order within each scheme, and the arrows show the sweep directions.}
    \label{four-color-smoothing}
\end{figure}

In what follows, we focus on the X-MCGS smoother in 2D; the extension to the 3D eight-color case is straightforward. The performance of the proposed X-MCGS operator is evaluated by solving Equation \eqref{PDE_discret} within a multigrid V-cycle and comparing it with the U-shape and the existing Z-shape variants \cite{ZHAO2022,li2020} of the 2D MCGS smoother. The experiment uses a $2048^2$ grid, with other details given in Section~\ref{sec3.1.2}. On each grid level, smoothing is performed using one of the three partitioning patterns—X-shape, U-shape, or Z-shape—as illustrated in Figure~\ref{four-color-smoothing}, with each smoothing step consisting of two MCGS iterations. Two ordering sequences are tested: 1234–1234, corresponding to two consecutive forward sweeps, and 1234–4321, corresponding to a forward sweep followed by a backward sweep. The latter may be regarded as a symmetric variant of MCGS, analogous to the classical Symmetric Gauss–Seidel iteration \cite{alvarez2022,bhatti2023}, and is generally expected to enhance both stability and smoothing performance in multigrid cycles.

\begin{figure}[H] 
    \centering
    \includegraphics[scale=0.45]{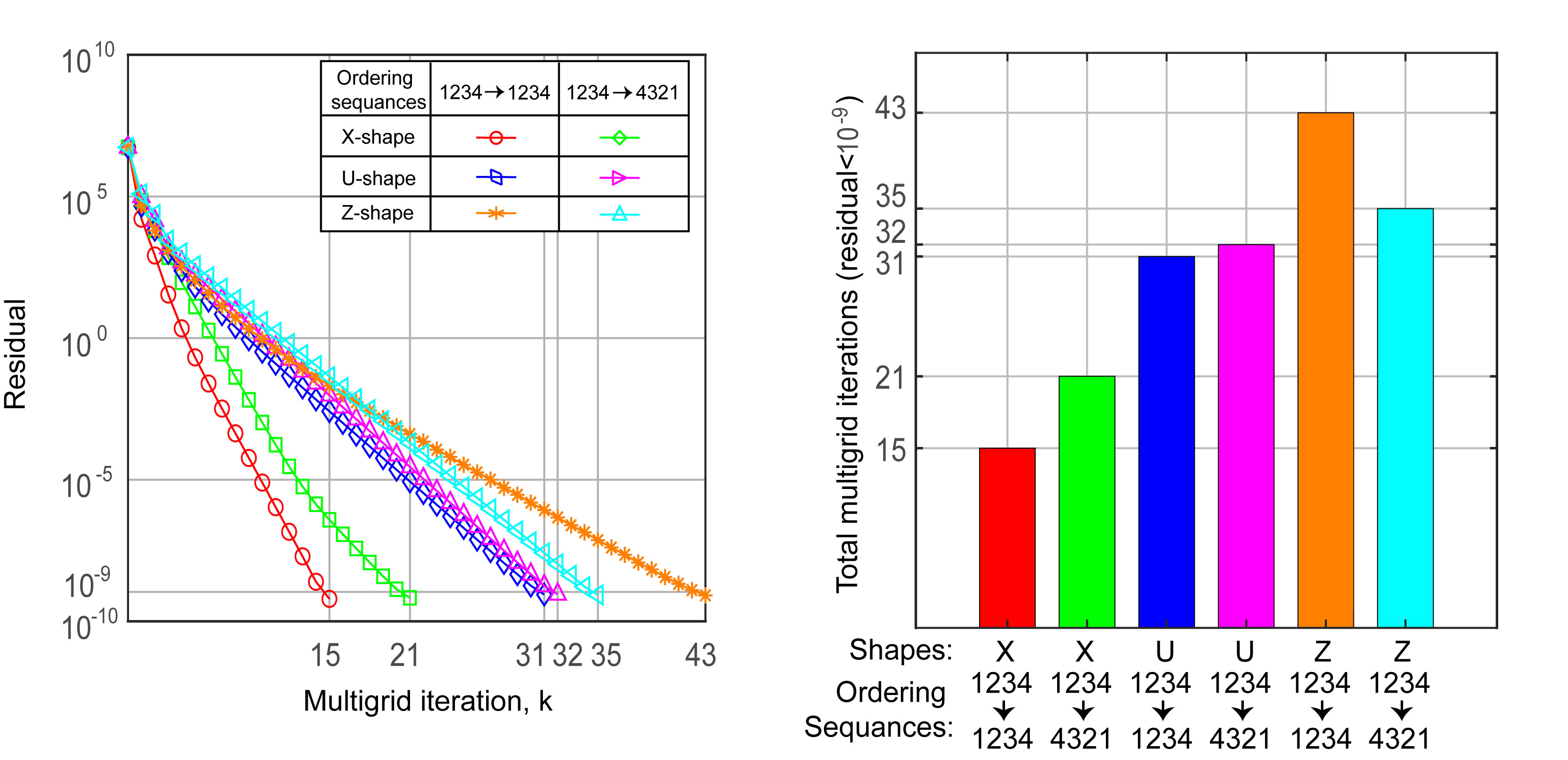}
    \caption{Convergence performance of six MCGS smoothing operators within the multigrid V-cycle in 2D. The experiment uses a $2048^2$ grid, with other details given in Section~\ref{sec3.1.2}. Left: Residual histories for the X-shape, U-shape, and Z-shape partitioning patterns under two ordering sequences (1234–1234 and 1234–4321). The iterations proceed until $\text{res}<10^{-9}$. Right: Total number of multigrid iterations required by each configuration to reach the stopping criterion. The results show that the X-shape scheme with the 1234–1234 ordering consistently achieves the fastest convergence, requiring only about half as many iterations as the slowest U-shape schemes and roughly one-third as many as the slowest Z-shape schemes. See \S~\ref{sec2.3} for details.}
    \label{smoothing-test}
\end{figure}

The convergence behavior of all six configurations is shown in Figure~\ref{smoothing-test}. In each case, we apply Algorithm~\ref{two-grid algorithm} to compute Equation~\eqref{PDE_discret} with the user-defined parameters in Table~\ref{tab:FASsymbols} set to $tol = 10^{-9}$, $k_{\text{Max}} = 100$, $s = 2$, and $meshLevel = 10$. The number of iterations $k$ required by each smoothing strategy is then recorded for comparison. The results show that the X-shape smoother with the 1234–1234 ordering consistently converges the fastest, requiring the fewest iterations among all tested configurations. Specifically, it needs only about half as many iterations as the slowest U-shape schemes and roughly one-third as many as the slowest Z-shape schemes. Based on these observations, we adopt the X-shape partition with the 1234–1234 ordering as the default smoothing operator for our GPU-accelerated multigrid solver. We further confirm that the same convergence behavior holds in the 3D case, although the results are omitted here for brevity.

\begin{figure}[H] 
    \centering
    \includegraphics[scale=0.55]{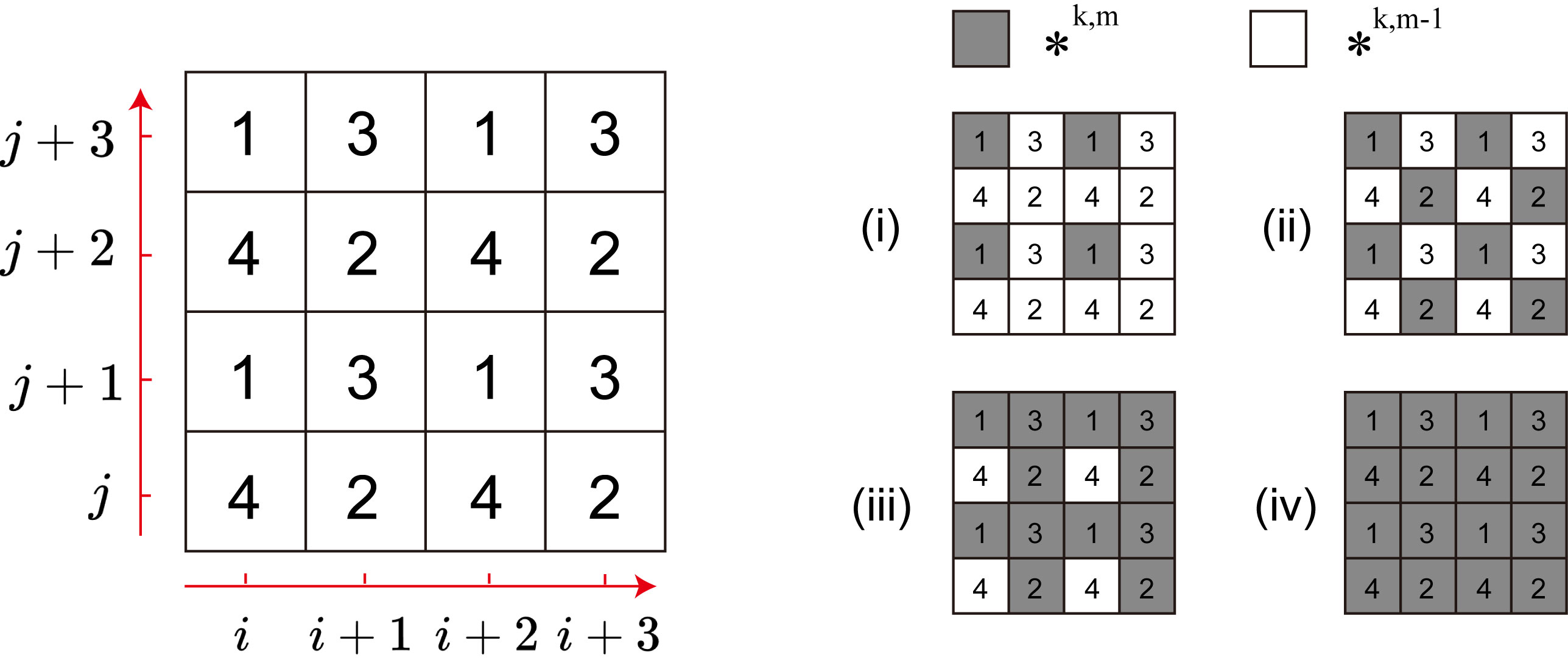}
    \caption{Illustration of the X-MCGS partitioning in 2D. \textbf{Left:} The computational grid is decomposed into four disjoint sets (1--4) based on the parity of the row and column indices. 
    \textbf{Right:} Sequential update order: (i) set~1, (ii) set~2, (iii) set~3, and (iv) set~4. White cells correspond to the values $*^{k,m-1}$ from the preceding smoothing step, whereas gray cells denote the new values $*^{k,m}$ updated from them during the current $m$-th smoothing of the $k$-th V-cycle. This decomposition enables vectorized sub-block updates, maximizes GPU parallelism, and eliminates conditional branching during parallel execution. See \S~\ref{sec2.3} for details.}
    \label{G-S}
\end{figure}

We now provide a more detailed description of the X-MCGS smoothing operator in 2D  in Steps~\textbf{4}, \textbf{10}, and \textbf{15} of Algorithm~\ref{two-grid algorithm}, in which the grid is partitioned into four submatrices according to the parity of the row and column indices (see Figure~\ref{G-S}, Left). This design eliminates conditional branching during parallel execution (see Listing~\ref{FCGS smoothing} for a MATLAB implementation of the 2D X-MCGS; a 3D MATLAB implementation is provided in \S\ref{app2}), allowing each update to be carried out in a fully vectorized form. Such a structure enables efficient batch operations on GPUs, thereby fully exploiting hardware parallelism. Specifically, grid points with $(i=\text{odd},\, j=\text{even})$ and $(i=\text{even},\, j=\text{odd})$ are assigned to sets~1 and~2, whereas points with $(i,j)$ both odd and both even are assigned to sets~3 and~4. 
The updates are then performed sequentially in the order: (i)~set~1, (ii)~set~2, (iii)~set~3, and (iv)~set~4 (see also Figure~\ref{G-S}, Right).

According to Algorithm 1, in the $k$-th multigrid V-cycle, once $p_{1,i,j}^{k,m-1}$ has been obtained after the $(m-1)$-th smoothing step, the update to $p_{1,i,j}^{k,m}$ is given by:

\begin{enumerate}[label=(\roman*)]
    \item Parallel update of set $1$:
    \begin{equation}\label{FCGS-1}
    \begin{minipage}{0.15\linewidth}
    \includegraphics[width=\linewidth]{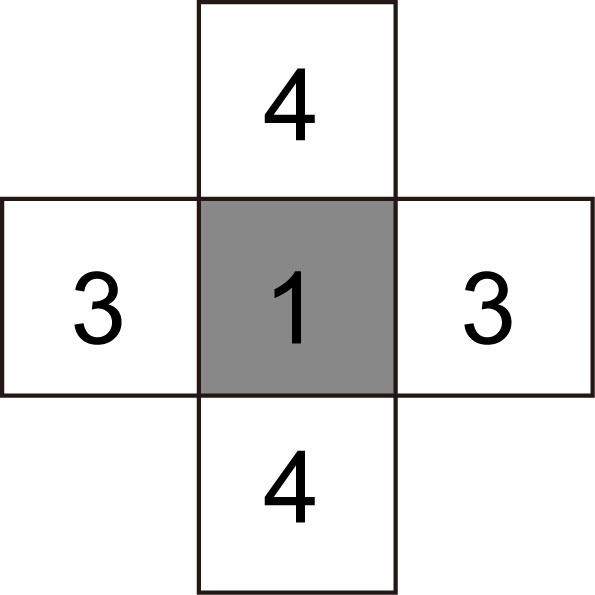} 
    \end{minipage}%
    \begin{minipage}{0.75\linewidth}
    \vspace{-1.4em}
    \[
    \left\{
    \begin{aligned}
    p_{1,i,j}^{k,m} &= \frac{p_{1,i+1,j}^{k,m-1} + p_{1,i-1,j}^{k,m-1} + p_{1,i,j+1}^{k,m-1} + p_{1,i,j-1}^{k,m-1} + h_1^2 f_{1,i,j}}{4 + h_1^2}, \\
    i &= \text{odd}, \ j = \text{even}.
    \end{aligned}
    \right.
    \]
    \end{minipage}
    \end{equation}

    \item Parallel update of set $2$:
    \begin{equation}\label{FCGS-2}
    \begin{minipage}{0.15\linewidth}
    \includegraphics[width=\linewidth]{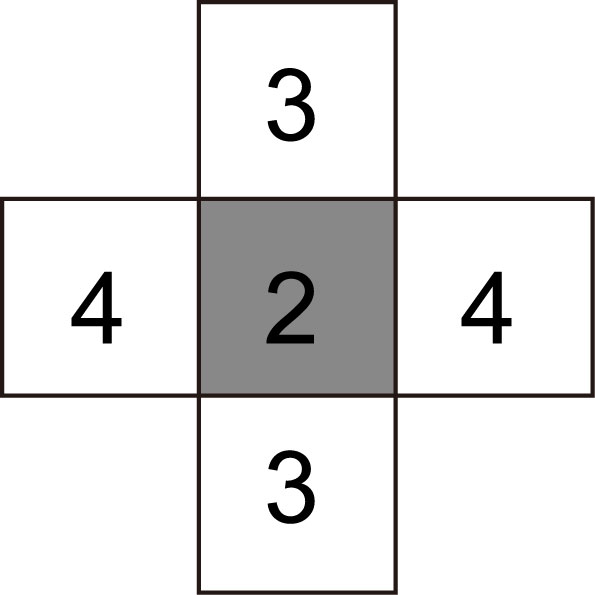} 
    \end{minipage}%
    \begin{minipage}{0.75\linewidth}
    \vspace{-1.4em}
    \[
    \left\{
    \begin{aligned}
    p_{1,i,j}^{k,m} &= \frac{p_{1,i+1,j}^{k,m-1} + p_{1,i-1,j}^{k,m-1} + p_{1,i,j+1}^{k,m-1} + p_{1,i,j-1}^{k,m-1} + h_1^2 f_{1,i,j}}{4 + h_1^2}, \\
    i &= \text{even}, \ j = \text{odd}.
    \end{aligned}
    \right.
    \]
    \end{minipage}
    \end{equation}

    \item Parallel update of set $3$:
    \begin{equation}\label{FCGS-3}
    \begin{minipage}{0.15\linewidth}
    \includegraphics[width=\linewidth]{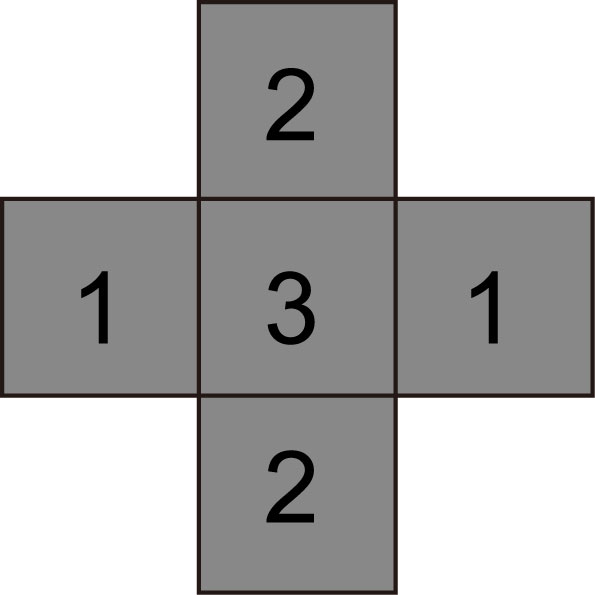} 
    \end{minipage}%
    \begin{minipage}{0.75\linewidth}
    \vspace{-1.4em}
    \[
    \left\{
    \begin{aligned}
    p_{1,i,j}^{k,m} &= \frac{p_{1,i+1,j}^{k,m} + p_{1,i-1,j}^{k,m} + p_{1,i,j+1}^{k,m} + p_{1,i,j-1}^{k,m} + h_1^2 f_{1,i,j}}{4 + h_1^2}, \\
    i &= \text{even}, \ j = \text{even}.
    \end{aligned}
    \right.
    \]
    \end{minipage}
    \end{equation}

    \item Parallel update of set $4$:
    \begin{equation}\label{FCGS-4}
    \begin{minipage}{0.15\linewidth}
    \includegraphics[width=\linewidth]{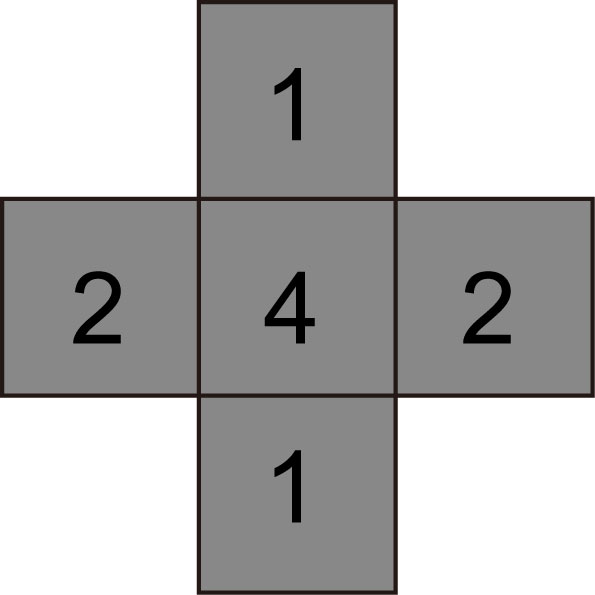} 
    \end{minipage}%
    \begin{minipage}{0.75\linewidth}
    \vspace{-1.4em}
    \[
    \left\{
    \begin{aligned}
    p_{1,i,j}^{k,m} &= \frac{p_{1,i+1,j}^{k,m} + p_{1,i-1,j}^{k,m} + p_{1,i,j+1}^{k,m} + p_{1,i,j-1}^{k,m} + h_1^2 f_{1,i,j}}{4 + h_1^2}, \\
    i &= \text{odd}, \ j = \text{odd}.
    \end{aligned}
    \right.
    \]
    \end{minipage}
    \end{equation}
\end{enumerate}

\begin{lstlisting}[caption={MATLAB implementation of the smoothing operator $\mathrm{S}$ using the X-shape multi-color Gauss–Seidel scheme in 2D}, label={FCGS smoothing}]
function p = Smoothing(f, p, M_1, N_1)
% X-shape Multi-Color Gauss-Seidel smoother in 2D

% Parallel update of set 1
vector_i=1:2:M_1-1; vector_j=2:2:N_1;
p(vector_i,vector_j) = GaussSeidel(f, p, vector_i, vector_j); 

% Parallel update of set 2
vector_i=2:2:M_1; vector_j=1:2:N_1-1;
p(vector_i,vector_j) = GaussSeidel(f, p, vector_i, vector_j); 

% Parallel update of set 3
vector_i=2:2:M_1; vector_j=2:2:N_1;
p(vector_i,vector_j) = GaussSeidel(f, p, vector_i, vector_j); 

% Parallel update of set 4
vector_i=1:2:M_1-1; vector_j=1:2:N_1-1;
p(vector_i,vector_j) = GaussSeidel(f, p, vector_i, vector_j); 

end
\end{lstlisting}

\section{Convergence and GPU performance of the matrix-free FAS multigrid solver on cell-centered grids}\label{sec3}
In this section, we evaluate the performance of the matrix-free FAS multigrid solver on cell-centered grids. Algebraic and asymptotic convergence tests are conducted to verify its robustness and accuracy. Execution efficiency is assessed by measuring the wall-clock time per iteration on different hardware platforms: NVIDIA RTX~4090 GPU (desktop), NVIDIA RTX~4060 GPU (desktop), NVIDIA RTX~3070 GPU (laptop), and a single Intel Core i9-13900KS CPU core. Finally, we showcase two large-scale GPU simulations: a grain-growth model with large crystal populations and a phase-separation model on spherical vesicles that captures the evolution of multicomponent giant unilamellar vesicles. These results establish the foundation for extending the solver to staggered grids and for addressing the incompressible Navier–Stokes equations in the subsequent sections.

\subsection{Convergence tests}\label{sec3.1}
We first perform algebraic and asymptotic convergence tests for the time-independent Equations~\eqref{PDE_discret}-\eqref{boundry2} in both 2D and 3D. In 2D, the domain is $\Omega = [0,1]^2$, with the exact solution $p_{1,i,j}^{\text{exact}} = \sin(\pi\sin(\pi x_i))\sin(\pi\sin(\pi y_j))$. In 3D, the domain is $\Omega = [0,1]^3$ with the exact solution 
    \[
p^{\text{exact}} = \sin(\pi\sin(\pi x_i))\sin(\pi\sin(\pi y_j))\sin(\pi\sin(\pi z_l)).
    \]
The finest-level grid is defined as $N_1^2$ in 2D or $N_1^3$ in 3D. We apply Algorithm~\ref{two-grid algorithm} with the user-defined parameters in Table~\ref{tab:FASsymbols} set to $tol = 10^{-9}$, $k_{\text{Max}} = 20$, $s = 2$, and $meshLevel = \log_2 N_1 - 1$.

\subsubsection{Algebraic convergence test}\label{sec3.1.1}
\begin{figure}[H] 
    \centering
    \begin{subfigure}{0.45\textwidth} 
        \centering
        \includegraphics[width=\linewidth]{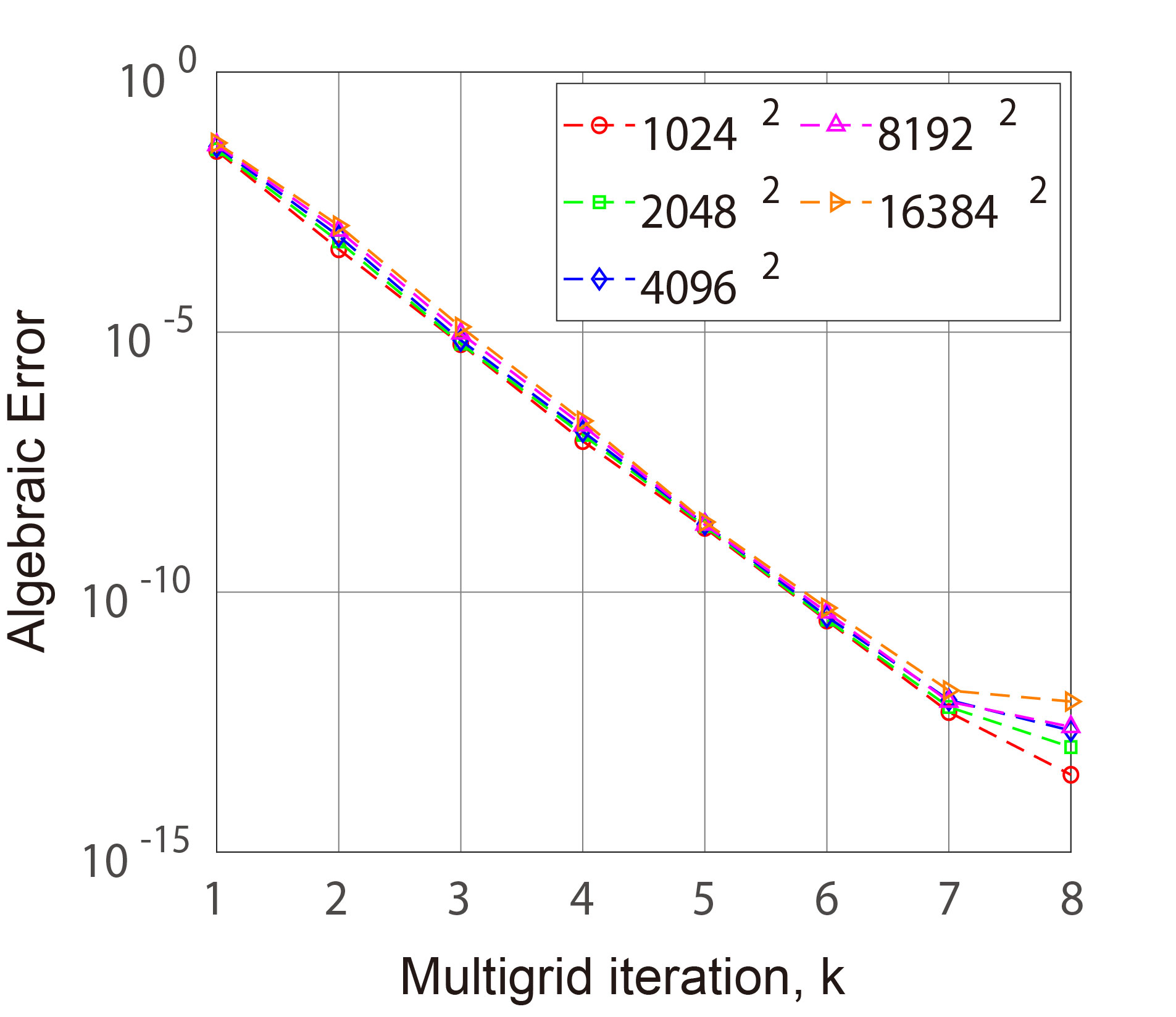}
    \end{subfigure}
    \hfill
    \begin{subfigure} {0.45\textwidth}
        \centering
        \includegraphics[width=\linewidth]{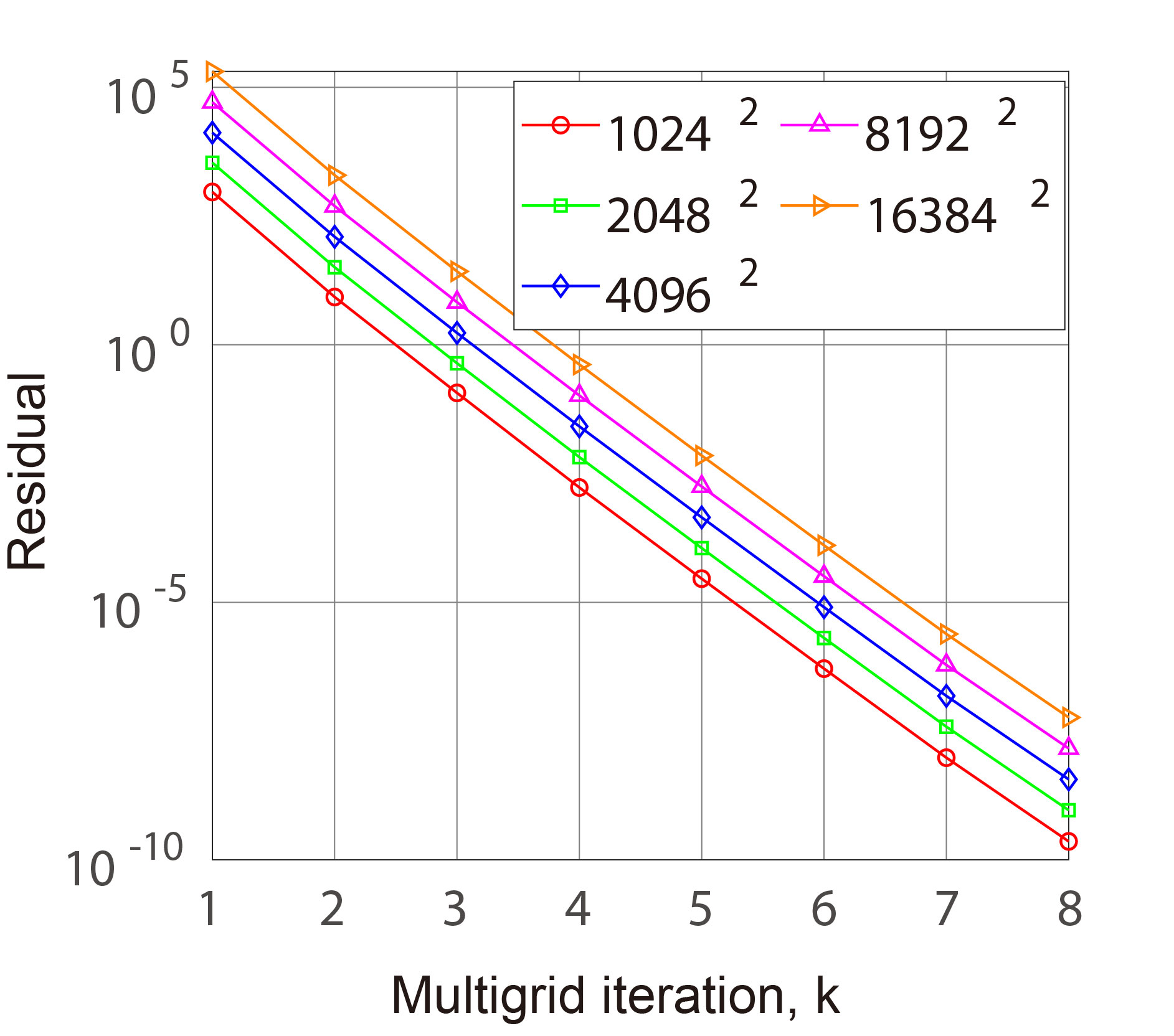}
    \end{subfigure}
    
    \caption{Algebraic convergence of the matrix-free FAS multigrid solver in 2D, where the Algebraic error, defined as $ \| p_{1,i,j}^{\text{exact}} - p_{1,i,j}^{k,m} \|_{L^2(\Omega_1)}$, and the residual, defined as Equation~\eqref{residual}  both show nearly $h$-independent reduction as the iteration count $k$ increases. The legend indicates different grid sizes $N_1^2$ by colored lines. See \S~\ref{sec3.1.1} for details.}
    \label{algebraic2D}
\end{figure}

\begin{figure}[H] 
    \centering
    \begin{subfigure}{0.45\textwidth} 
        \centering
        \includegraphics[width=\linewidth]{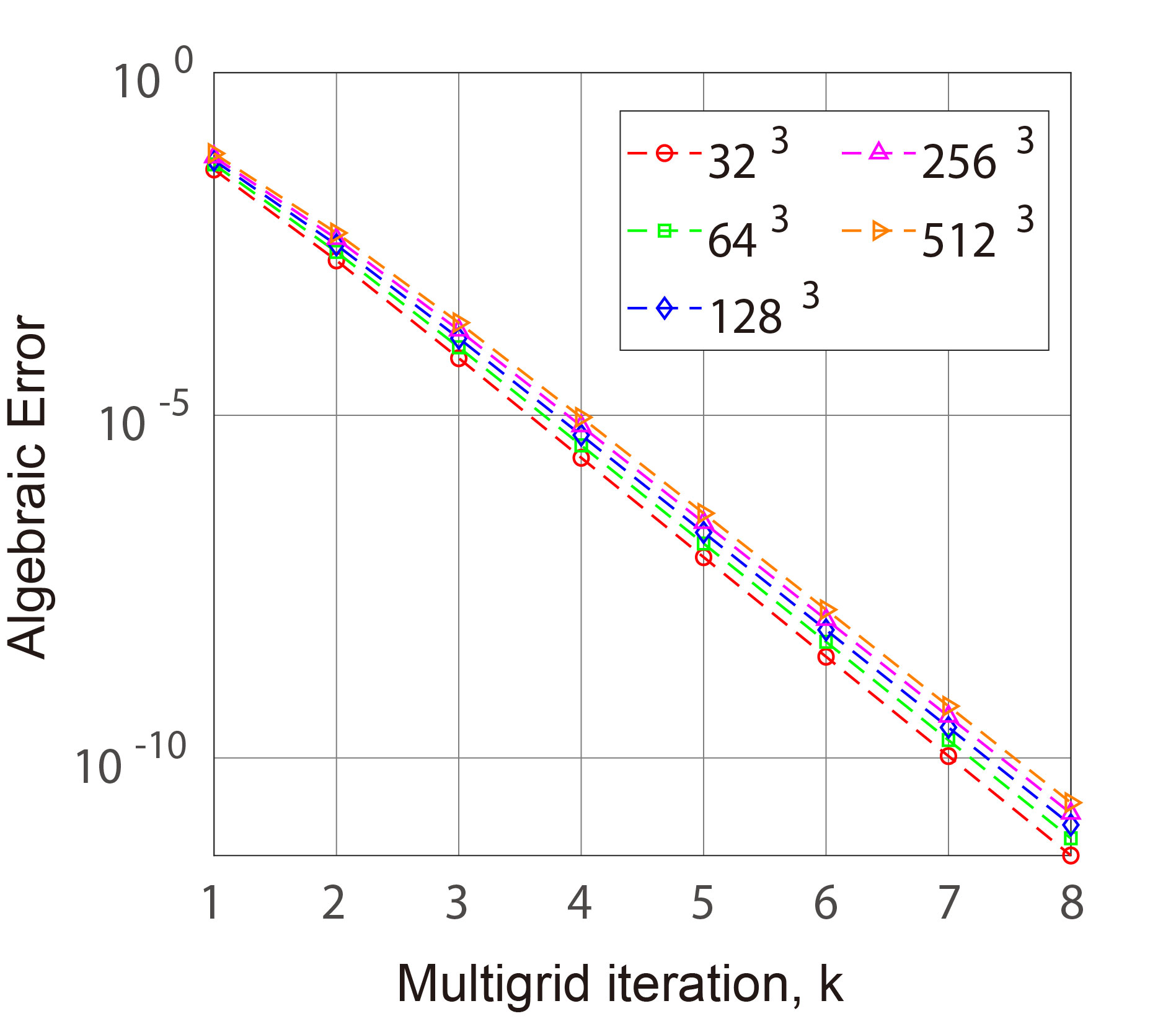}
    \end{subfigure}
    \hfill
    \begin{subfigure} {0.45\textwidth}
        \centering
        \includegraphics[width=\linewidth]{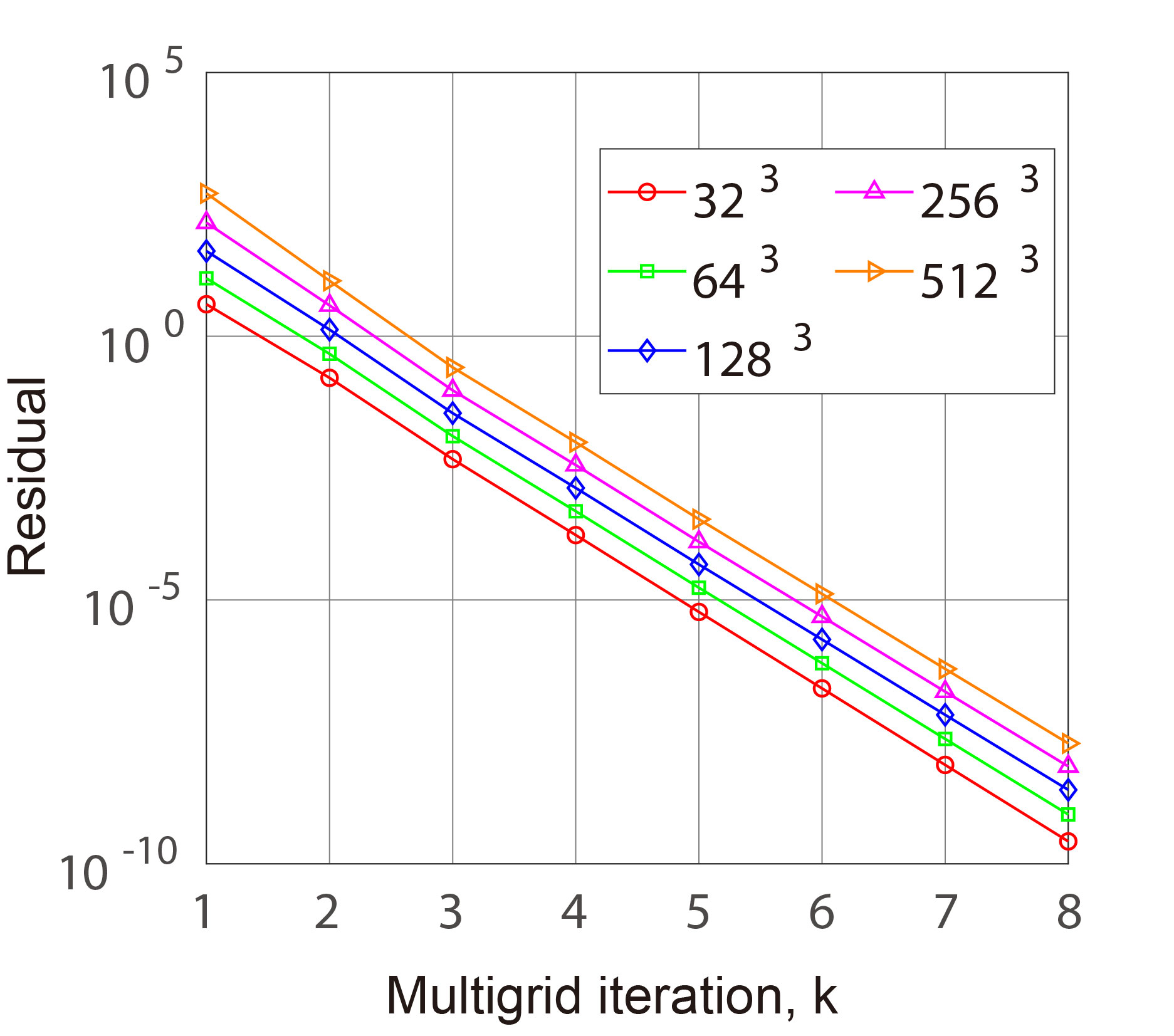}
    \end{subfigure}
    
    \caption{Algebraic convergence of the matrix-free FAS multigrid solver in 3D, where the Algebraic error, defined as $\| p_{1,i,j,l}^{\text{exact}} - p_{1,i,j,l}^{k,m} \|_{L^2(\Omega_1)}$, and the residual defined as $\| f_{1,i,j,l} - L_{h_1}(p_{1,i,j,l}^{(m)}) \|_{L^2(\Omega_1)}$ both show nearly $h$-independent reduction as the iteration count $k$ increases. The legend indicates different grid sizes $N_1^3$ by colored lines. See \S~\ref{sec3.1.1} for details.}
    \label{algebraic3D}
\end{figure}

%{\bf Algebraic convergence test:} 
This test examines whether the multigrid algorithm achieves grid-independent convergence. In 2D, the source term is obtained by substituting the exact solution $p^{\text{exact}}$ into the discrete Equation~\eqref{PDE_discret}, defining $f_{1,i,j} = L_{h_1}(p_{1,i,j}^{\text{exact}})$, so that errors arise only from machine precision. The initial condition $p_{1,i,j}^0$ is set to random values in $[0,1]$. The 3D case follows the same procedure. As shown in Figures~\ref{algebraic2D} and \ref{algebraic3D}, both the error and residual decay at similar rates across different grids as the iteration count $k$ increases, demonstrating grid-independent convergence in both 2D and 3D.

\subsubsection{Asymptotic convergence tests}\label{sec3.1.2}
%{\bf Asymptotic convergence tests:}
In 2D, the source term $f_{1,i,j}$ is obtained by substituting the exact solution $p^{\text{exact}}$ into the continuous Equation~\eqref{PDE}, with the initial condition $p_{1,i,j}^{1,0} = 0$ at all grid points. The 3D case follows the same procedure. This test evaluates the asymptotic error reduction rate. As shown in Tables~\ref{err_2D} and \ref{err_3D}, both 2D and 3D results demonstrate second-order convergence, confirming the accuracy of the FAS multigrid solver.

\begin{table}[H]
\tabcolsep=1.2cm
\centering%把表居中
\begin{tabular}{ccc}%四个c代表该表一共四列，内容全部居中
\toprule%第一道横线
Fine grid size $N_1^2$&Error&Order\\
\midrule%第二道横线 
$1024^2$&$2.58\times10^{-6}$&-\\
$2048^2$&$6.46\times10^{-7}$&$2.00$\\
$4096^2$&$1.61\times10^{-7}$&$2.00$\\
$8192^2$&$4.03\times10^{-8}$&$2.00$\\
$16384^2$&$1.01\times10^{-8}$&$2.00$\\
\bottomrule%第三道横线
\end{tabular}
\caption{Asymptotic Convergence test results for the 2D matrix-free FAS multigrid solver on cell-centered grids ranging from $1024^2$ to $16384^2$. The error, measured as $\| p_{1,i,j}^{\text{exact}} - p_{1,i,j}^{k,m} \|_{L^2(\Omega_1)}$, confirms the expected second-order convergence. See \S~\ref{sec3.1.2} for details.}
\label{err_2D}
\end{table}

\begin{table}[H]

\tabcolsep=1.2cm
\centering%把表居中
\begin{tabular}{ccc}%四个c代表该表一共四列，内容全部居中
\toprule%第一道横线
Fine grid size $N_1^3$&$\mathrm{Error}$&Order\\
\midrule%第二道横线 
$32^3$&$1.74\times10^{-3}$&-\\
$64^3$&$4.33\times10^{-4}$&$2.01$\\
$128^3$&$1.08\times10^{-4}$&$2.00$\\
$256^3$&$2.70\times10^{-5}$&$2.00$\\
$512^3$&$6.75\times10^{-6}$&$2.00$\\
\bottomrule%第三道横线
\end{tabular}
\caption{Asymptotic convergence test results for the 3D matrix-free FAS multigrid solver on cell-centered grids ranging from $32^3$ to $512^3$. The error, measured as  $ \| p_{1,i,j,l}^{\text{exact}} - p_{1,i,j,l}^{k,m} \|_{L^2(\Omega_1)}$, confirms the expected second-order convergence. See \S~\ref{sec3.1.2} for details.}
\label{err_3D}
\end{table}

\subsection{Execution performance}\label{sec3.2}
We next evaluate the performance of the proposed matrix-free FAS multigrid solver on different GPU devices and a CPU for the same benchmark problem. GPU performance is largely determined by core count, frequency, bandwidth, and memory capacity, the latter often limiting the maximum grid size and number of variables handled in parallel. To compare architectures, we conducted experiments on three GPUs: RTX~4090, RTX~4060 and RTX~3070. The main hardware specifications are summarized in Table~\ref{GPU Specifications}.

\begin{table}[H]
\centering
\tabcolsep=0.30cm
\begin{tabular}{cccc}
\toprule
\textbf{GPU platform} & \textbf{Desktop} & \textbf{Desktop} & \textbf{Laptop} \\
\midrule
Model   & RTX~4090 & RTX~4060 & RTX~3070\\
Memory (GB) & 24 & 8  & 8      \\
CUDA cores  &  16384  & 3072 & 5120   \\
Bandwidth (GB/s)& 1008 & 272  &448 \\
\bottomrule%第三道横线
\end{tabular}
\caption{Comparison of GPU architectures across different platforms. See \S~\ref{sec3.2} for details.}
\label{GPU Specifications}
\end{table}
%From the above table, we understand that the RTX 3070 and RTX 4060 are equipped with 8 GB of memory each, whereas the RTX 4090 provides a significantly larger capacity of 24 GB. Taking the single-variable Equation ~\eqref{PDE} as an example, the maximum solvable grid size for 8 GB of memory reaches $8192^2$ in 2D and $256^3$ in 3D. In addition, with 24 GB of memory, the maximum grid size increases to $16384^2$ in 2D and $512^3$ in 3D. 

%We then compare the performance of our code on CPU and GPU architectures. Specifically, we use a single CPU core (Intel\textsuperscript{\textregistered} Core\texttrademark~i9-13900KS) for the CPU-side evaluation. It is worth noting that such comparisons are inherently challenging. First, there is no universally accepted standard for evaluating and comparing the performance of CPU and GPU systems. Second, code performance can vary significantly across different hardware platforms, making it difficult to conduct unbiased tests. Consequently, the results presented below should be considered as an initial approximation.

We apply Algorithm~\ref{two-grid algorithm} with the user-defined parameters in Table~\ref{tab:FASsymbols} set to $tol = 10^{-9}$, $k_{\text{Max}} = 1$, $s = 2$, and $meshLevel = \log_2 N_1 - 1$, and then compare solver performance on GPU and CPU architectures. For consistency, ten multigrid iterations are executed in each test, and the average execution time per iteration is recorded. The results for 2D and 3D cases across different mesh sizes are summarized in Tables~\ref{runningTime_2D} and \ref{runningTime_3D}, respectively.

\begin{table}[H]
\tabcolsep=0.12cm
\centering%把表居中
\begin{tabular}{ccccc}%四个c代表该表一共四列，内容全部居中
\toprule%第一道横线
\textbf{2D} & \multicolumn{4}{c}{\textbf{Execution time for 1 iteration step (s)}} \\
\midrule
Fine grid size $N_1^2$&i9-13900KS &\thead{{RTX~3070 (Laptop)}\\{(8~GB)}}&\thead{{RTX~4060}\\{(8~GB)}}&\thead{{RTX~4090}\\{(24~GB)}}\\
\midrule%第二道横线 
$1024^2$&$0.0582$&$0.1235$&$0.0370$&$0.0351$\\
$2048^2$&$0.3932$&$0.1376$&$0.0503$&$0.0395$\\
$4096^2$&$2.1319$&$0.1651$&$0.1386$&$0.0530$\\
$8192^2$&$9.1859$&$0.6721$&$0.6394$&$0.1484$\\
$16384^2$&$36.6231$&Out of memory&Out of memory&$0.9396$\\
\bottomrule%第三道横线
\end{tabular}
\caption{Per-iteration execution time of the 2D matrix-free FAS method on cell-centered grids, comparing a single CPU core (i9-13900KS) with different GPU devices. Each value represents the average of ten runs. See \S~\ref{sec3.2} for details.}
\label{runningTime_2D}
\end{table}

\begin{table}[H]
\tabcolsep=0.12cm
\centering%把表居中
\begin{tabular}{ccccc}%四个c代表该表一共四列，内容全部居中
\toprule%第一道横线
\textbf{3D} & \multicolumn{4}{c}{\textbf{Execution time for 1 iteration step (s)}} \\
\midrule
Fine grid size $N_1^3$&i9-13900KS &\thead{{RTX~3070 (Laptop)}\\{(8~GB)}}&\thead{{RTX~4060}\\{(8~GB)}}&\thead{{RTX~4090}\\{(24~GB)}}\\
\midrule%第二道横线 
$128^3$&$0.1531$&$0.1669$&$0.0548$&$0.0511$\\
$256^3$&$2.7970$&$0.2203$&$0.1626$&$0.0710$\\
$512^3$&$21.5990$&Out of memory&Out of memory&$0.4633$\\
\bottomrule%第三道横线
\end{tabular}
\caption{Per-iteration execution time of the 3D matrix-free FAS method on cell-centered grids, comparing a single CPU core (i9-13900KS) with different GPU devices. Each value represents the average of ten runs. See \S~\ref{sec3.2} for details.}%标题
\label{runningTime_3D}
\end{table}
It is evident that GPU parallelization substantially accelerates the computations. In 2D with a $1024^2$ grid, the RTX~4060 and RTX~4090 are only slightly faster than the CPU, whereas the RTX~3070 is about twice as slow. For $2048^2$ and larger grids, the advantage of GPUs becomes evident: on the $8192^2$ grid, the RTX~4090 achieves a speedup of approximately 61×, while the RTX~4060 and RTX~3070 deliver comparable performance, each about 14× faster. On the $16384^2$ grid, limited GPU memory allows only the RTX~4090 to handle the problem, reaching about 39× the CPU speed.

In 3D, similar trends are observed. With a $128^3$ grid, the RTX~4060 and RTX~4090 are about 2.9× faster than the CPU, while the RTX~3070 performs at roughly the same speed. At $256^3$, all three GPUs show clear advantages: the RTX~4090, RTX~4060, and RTX~3070 are about 39×, 17×, and 12× faster, respectively. At $512^3$, GPU memory constraints leave only the RTX~4090 capable of completing the computation, with a speedup of about 46× over the CPU.

\subsection{Numerical examples}\label{sec3.3}
We next present two numerical examples performed on a single RTX~4090 GPU to demonstrate the solver's capability in handling multicomponent systems and large-scale grids.

\subsubsection {Grain growth} \label{sec3.3.1}
Grain growth strongly affects the properties of polycrystalline materials, including conductivity, thermal transport, and corrosion resistance. Since grain boundaries and junctions are difficult to probe experimentally, computational modeling—particularly the phase-field method—is widely employed to capture interface motion and topological transitions. The Fan–Chen model \cite{fan1997} represents grain orientations with multiple order parameters; although real materials exhibit a continuous spectrum of orientations ($q=\infty$), simulations use a finite $q$, and a sufficiently large value is required to reproduce growth kinetics. The governing equations follow a simplified phase-field formulation, where orientation variables ${\eta_i(\mathbf{r},t)}$ evolve via the Ginzburg–Landau equation:
\begin{equation}\label{grain-growth-equations}
\left\{
\begin{aligned}
\frac{\partial \eta_i}{\partial t} &= -L_i \left( 
    -\alpha \eta_i + \beta \eta_i^3 + 2 \gamma \eta_i \sum_{\substack{j=1 \\ j \ne i}}^{p} \eta_j^2 
    - \kappa_i \nabla^2 \eta_i 
\right), \\
\varphi(\mathbf{r}) &= \sum_{i=1}^{p} \eta_i^2(\mathbf{r}), \qquad i = 1, 2, \ldots, q, 
\end{aligned}
\right.
\end{equation}
where parameters are chosen as $L_i=1.0$, $\alpha=\beta=\gamma=1.0$, $\kappa_i=2.0$.

In our simulations, periodic boundary conditions are imposed in all directions. The initial condition is specified by small random values in the range $-0.001 < \eta_i < 0.001$, representing a disordered, liquid-like configuration. We simulate the system in both 2D and 3D on a single RTX~4090 GPU (24~GB) and evaluate the maximum number of orientation variables $q$ that can be accommodated on fixed grids.
 
\begin{figure}[H] 
    \centering
    \includegraphics[scale=0.22]{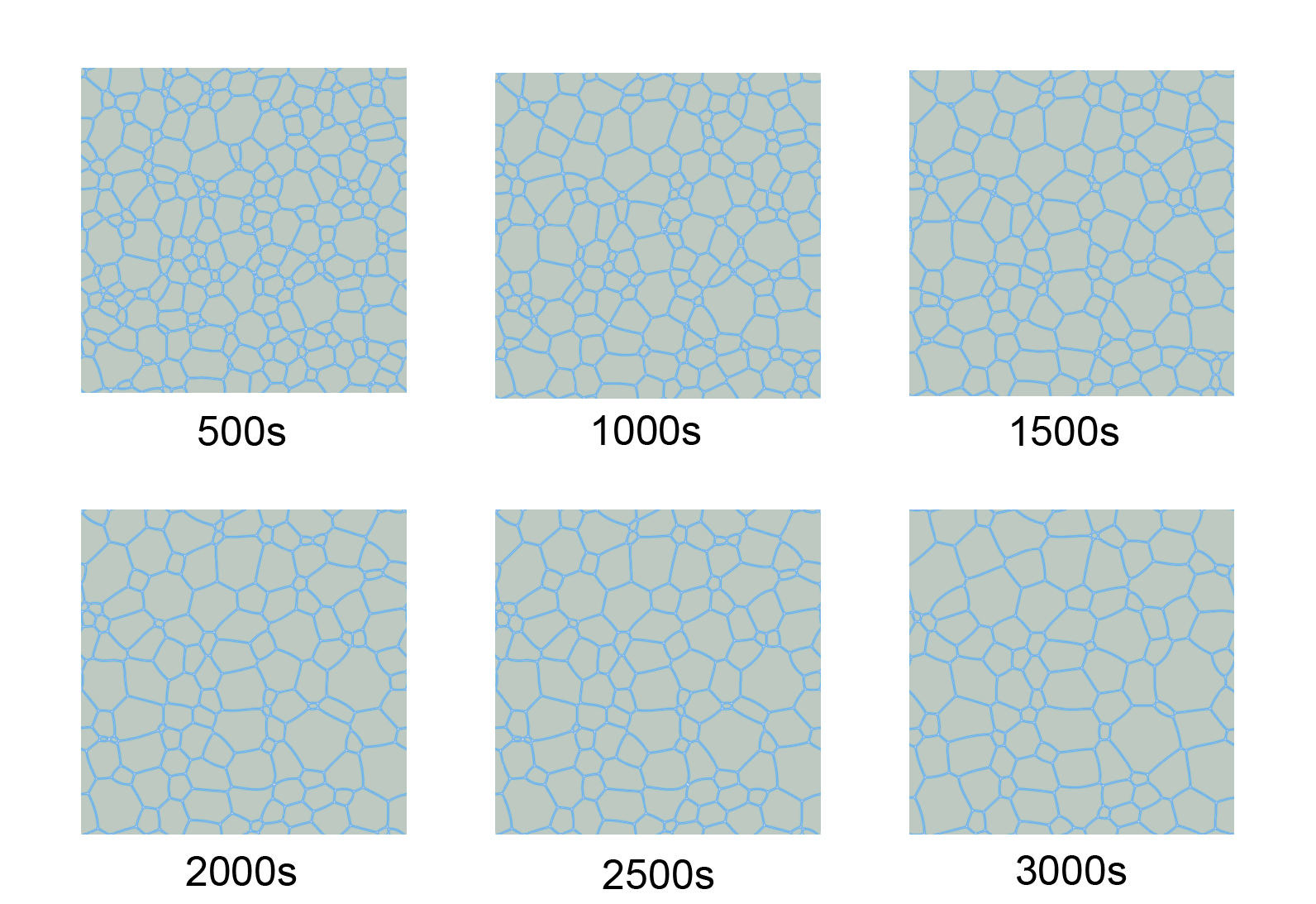}
    \caption{Microstructural evolution of grains on a $512^2$ grid with $q = 1189$ orientation variables, simulated on a single RTX 4090 GPU (24~GB). Snapshots at $500$, $1000$, $1500$, $2000$, $2500$ and $3000~s$, illustrate the progressive growth and coarsening of grains over simulation time. See \S~\ref{sec3.3.1} for details.}
    \label{grain_growth_2D}
\end{figure}

In 2D with a $512^2$ grid, up to $q=1189$ orientations can be represented. The computational domain is $[0,1024]^2$, with the finest grid spacing along Cartesian axes set to $h_1 = 2.0$, and a time step of $\Delta t = 0.25$ is used. Figure~\ref{grain_growth_2D} shows the microstructural evolution for $q=1189$: the early stage corresponds to crystallization, followed by growth and coarsening of grains over time. Additional simulations with $q=4,36,50,100,200$ exhibit similar qualitative behavior.

\begin{figure}[H] 
    \centering
    \includegraphics[width=\linewidth]{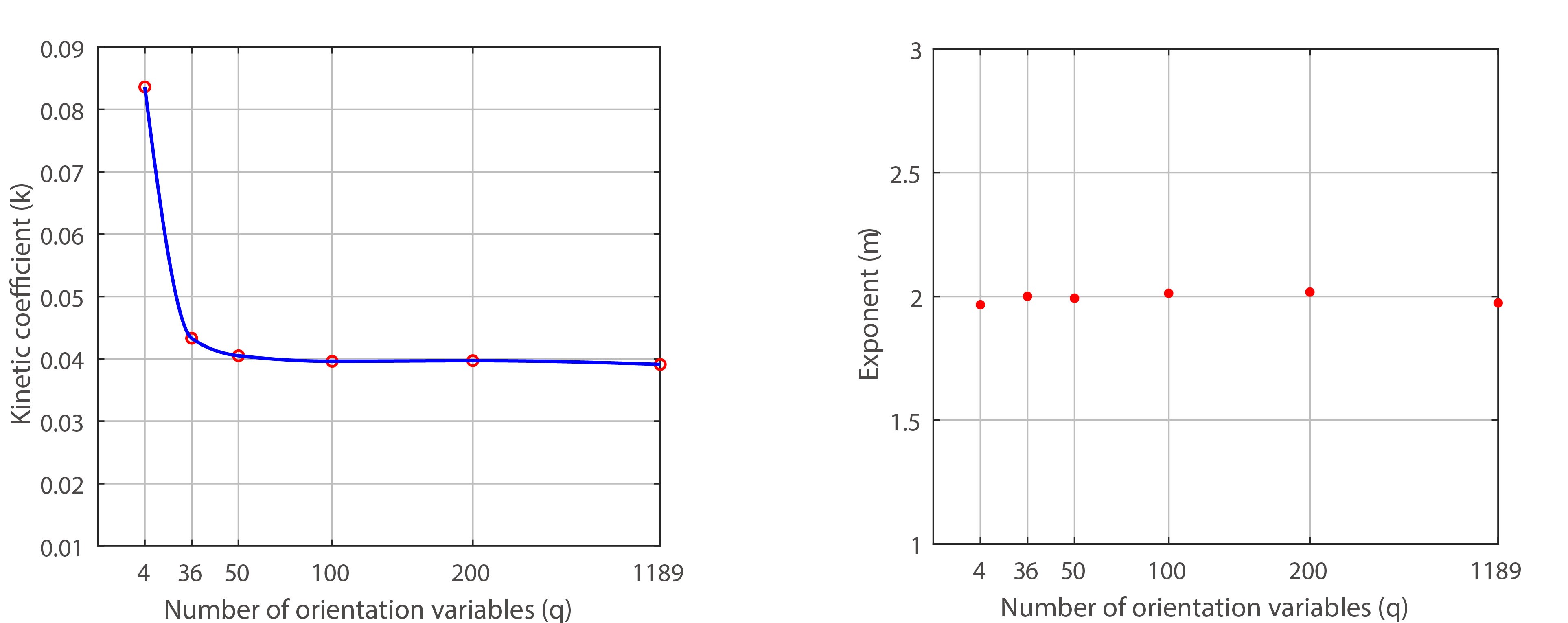}
    \caption{Fitted grain growth exponent $m$ and coefficient $k$ for different numbers of orientation variables $q$ in 2D. See \S~\ref{sec3.3.1} for details. }
\label{2Dgrain-size}
\end{figure}

Using the method of \cite{fan1997}, we assess correctness by fitting the average grain size $R$ to the form $R_t^m-R_0^m=kt$. As shown in Figure~\ref{2Dgrain-size}, the kinetic coefficient $k$ converges as $q$ increases, while the growth exponent $m$ remains nearly constant at $2.0$, consistent with previous results.

\begin{figure}[H] 
    \centering
    \includegraphics[scale=0.85]{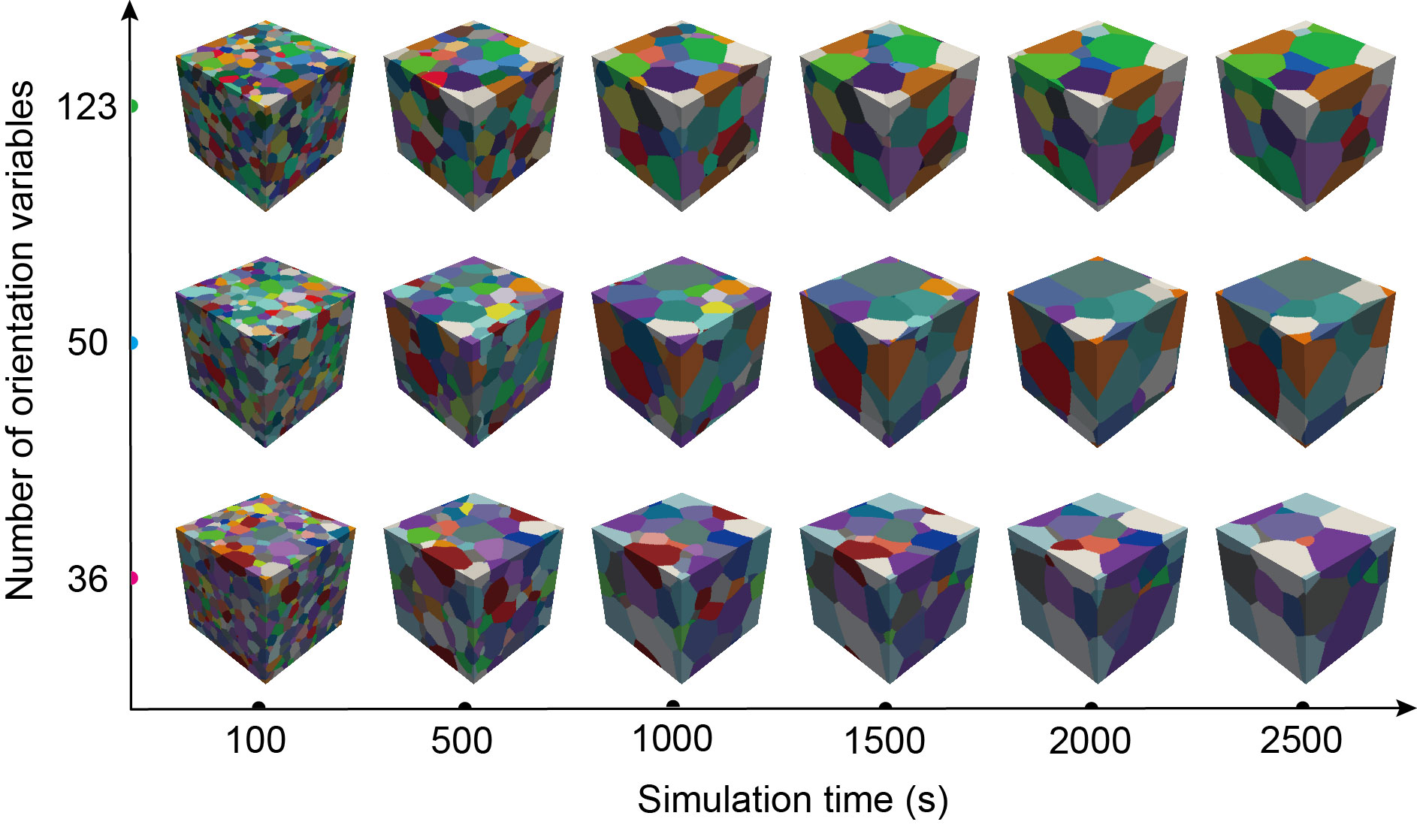}
    \caption{Microstructural evolution of grains on a $128^3$ grid with $q = 36,\,50,\,123$ orientation variables, simulated on a single RTX 4090 GPU (24~GB). Snapshots at $100$, $500$, $1000$, $1500$, $2000$ and $2500~s$, illustrate the growth and coarsening of grains over simulation time. See \S~\ref{sec3.3.1} for details. }
    \label{grain_growth_3D}
\end{figure}

\begin{figure}[H] 
    \centering
    \includegraphics[width=\linewidth]{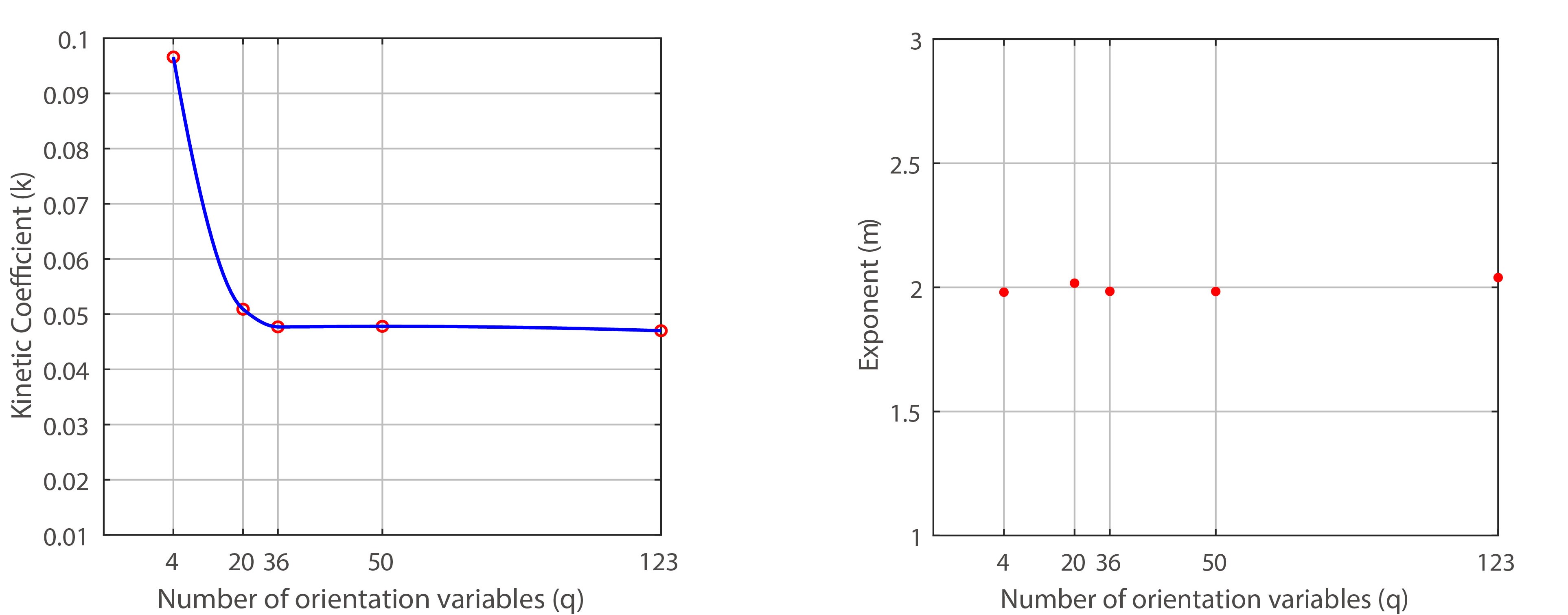}
    \caption{Fitted grain growth exponent $m$ and coefficient $k$ for different numbers of orientation variables $q$ in 3D. See \S~\ref{sec3.3.1} for details. }
\label{3Dgrain-size}
\end{figure}

In 3D, simulations are performed on a $128^3$ grid, and the system can reliably accommodate up to $q=123$ orientations. The computational domain is $[0,256]^3$, with the finest grid spacing along Cartesian axes set to $h_1 = 2.0$, and a time step of $\Delta t = 0.25$. Figure~\ref{grain_growth_3D} illustrates grain evolution for $q=36,50,123$, again showing the formation of well-defined grains followed by coarsening. The fitted values of $m$ and $k$, shown in Figure~\ref{3Dgrain-size}, agree well with \cite{fan1997}, further validating the solver's accuracy.

\subsubsection{Phase Separation on Multicomponent Vesicles}\label{sec3.3.2}

We further demonstrate the effectiveness of our solver by simulating phase separation on a fixed spherical multicomponent vesicle, motivated by the experiments of Wang et al.~\cite{wang2022lipid}, which reported coarsening of lipid domains on giant unilamellar vesicles (GUVs). While their study employed a surface Navier–Stokes–Cahn–Hilliard model that couples membrane hydrodynamics with phase separation, we adopt a simplified bulk phase-field formulation by neglecting membrane flow. The resulting nondimensional model is
\begin{align}
    \delta(\phi)\partial_t c &= \frac{1}{\text{Pe}} \boldsymbol{\nabla} \cdot \left( \delta(\phi) m(c) \boldsymbol{\nabla} \omega_c \right), \\ \label{omegac}
    \delta(\phi) \omega_c &=  \delta(\phi)\eta_c \frac{f^\prime(c)}{\varepsilon_c} -  \eta_c \varepsilon_c \boldsymbol{\nabla}\cdot(\delta(\phi) \boldsymbol{\nabla} c),
\end{align}
where the diffuse domain method \cite{ratz_pdes_2006,teigen2009diffuse, li2009solving} is employed through $\delta(\phi) = \eta_\phi {\varepsilon_\phi} ^{-1}{f(\phi)} + \eta_\phi \varepsilon_\phi {|\nabla \phi|^2}/{2}$, 
a regularized delta function characterizing the vesicle membrane. The parameter $\eta_\phi$ scales the diffuse-interface formulation to its sharp-interface counterpart. The phase-field variable $\phi(\boldsymbol{x},t)$ localizes the vesicle membrane, with $\phi=1$ inside, $\phi=-1$ outside, and a smooth transition across an interfacial layer of thickness $\varepsilon_\phi$. The second phase-field variable $c(\boldsymbol{x},t)$ represents the membrane composition, with $c=1$ and $c=0$ corresponding to the liquid-disordered ($l_d$) and liquid-ordered ($l_o$) phases, respectively, separated by a diffuse interface of thickness $\varepsilon_c$. Its chemical potential is given by Equation~\eqref{omegac}, where the double-well potential is $f(c)=c^2(1-c)^2/4$. The mobility is taken as $m(c)=c^2(1-c)^2+\varepsilon_c$. Finally, the nondimensional Peclet number $\text{Pe}$ characterizes the relative rate of lipid diffusion compared with the characteristic time scale.

In our simulations, the computational domain is ${\Omega} = [-1.5, 1.5]^3$, discretized on a $256^3$ grid. The fixed vesicle radius is $r_0 = 1$ with membrane thickness $\varepsilon_{\phi} = 0.01$. The initial phase-field variable is prescribed as
\begin{align}\label{initial shape}
\phi(\boldsymbol{x}, 0) = \tanh\left(\frac{r_0 - \sqrt{x^2+y^2+z^2}}{2\sqrt{2}\,\varepsilon_\phi}\right), \quad \boldsymbol{x} \in \bar{\Omega},
\end{align}
representing a smooth diffuse interface centered on the vesicle surface.
To model an initially homogeneous liposome consistent with \cite{wang2022lipid}, we generate random numbers $ c_{\text{rand}}(\boldsymbol{x})$ uniformly distributed in $[0,1]$. The initial surface fraction $ c_\Gamma(\boldsymbol{x}, 0) $ is then specified by
\begin{align} \label{ic-phase-sep-c_gamma}
    c_{\Gamma}(\boldsymbol{x}, 0) = 
    \begin{cases}
        ~~1 \quad (l_{d}), & \text{if } c_{\text{rand}}(\boldsymbol{x}) \in [0, a_{ld}], \\
        ~~0 \quad (l_{o}), & \text{if } c_{\text{rand}}(\boldsymbol{x}) \in (a_{ld}, 1],
    \end{cases}
\end{align}
where $a_{ld}$ and $a_{lo} = 1 - a_{ld}$ represent the prescribed area fractions of the $l_d$ (red) and $l_o$ (blue) phases. In addition, no-flux boundary conditions are imposed on $\phi$, $c$, and $\omega_c$ at the domain boundary $\partial \bar{\Omega}$.

Following the experimental parameters in \cite{wang2022lipid}, we set the interfacial thickness for membrane components to $\varepsilon_c=0.01$, the membrane thickness to $\varepsilon_\phi=0.01$, and the Peclet number to $\text{Pe}=1$. The time step is chosen as $\Delta t=10^{-2}$.

We consider two area fractions, $a_{ld} = 0.3$ and $a_{ld} = 0.7$, corresponding to the experimental conditions in \cite{wang2022lipid}. To account for stochastic variability, 10 independent simulations are performed with different random seeds for initializing $c$ at each area fraction. The realization most closely resembling the experimental observations is selected for comparison.

\begin{figure}[H]
   \centering
   \includegraphics[width=1\linewidth]{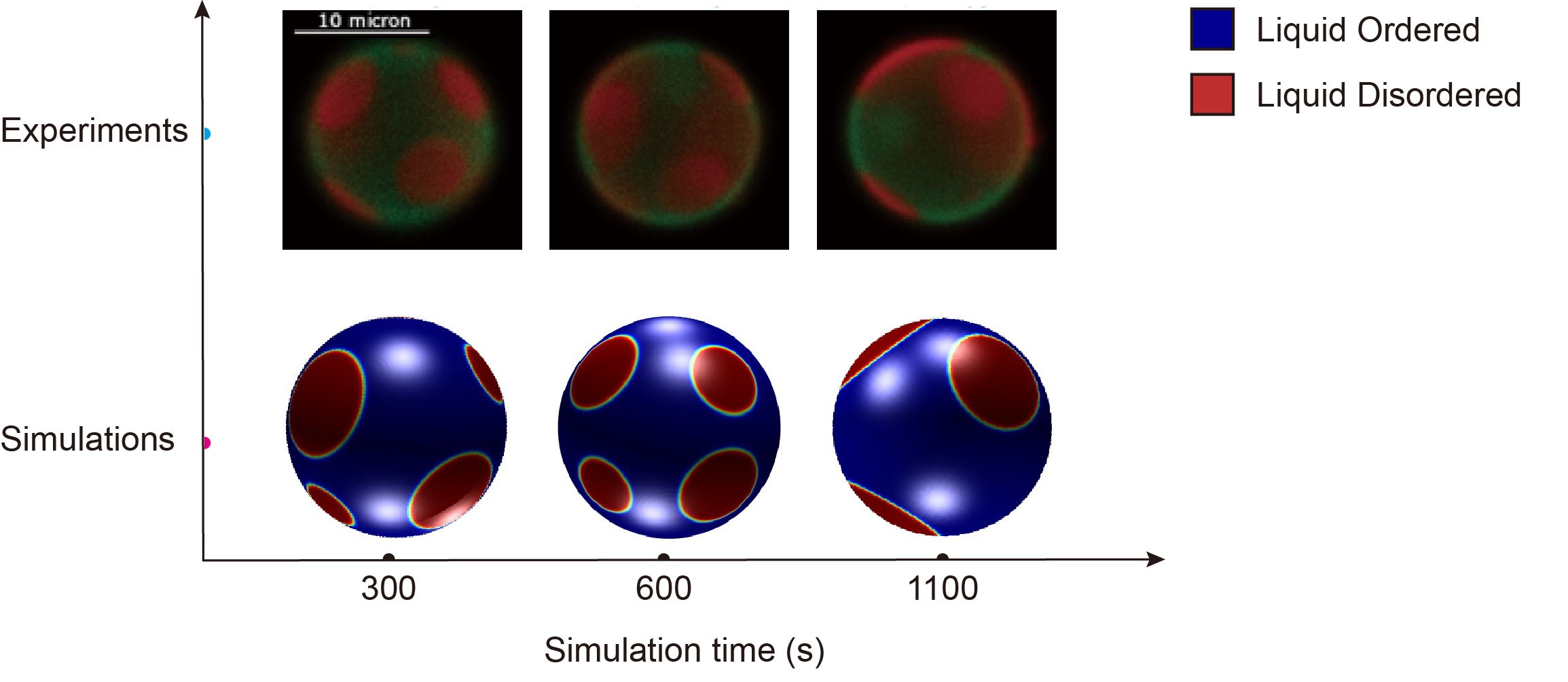}
   \caption{Comparison of phase separation dynamics for area fraction $a_{l_d} = 0.3$: experimental images (top) from \cite{wang2022lipid}; numerical simulation results (bottom) illustrating the formation of ${l_d}$ domains (red) within the ${l_o}$ phase (blue). See \S~\ref{sec3.3.2} for details.}
   \label{fig.37od}
\end{figure}

\begin{figure}[H]
\centering
\includegraphics[width=1\linewidth]{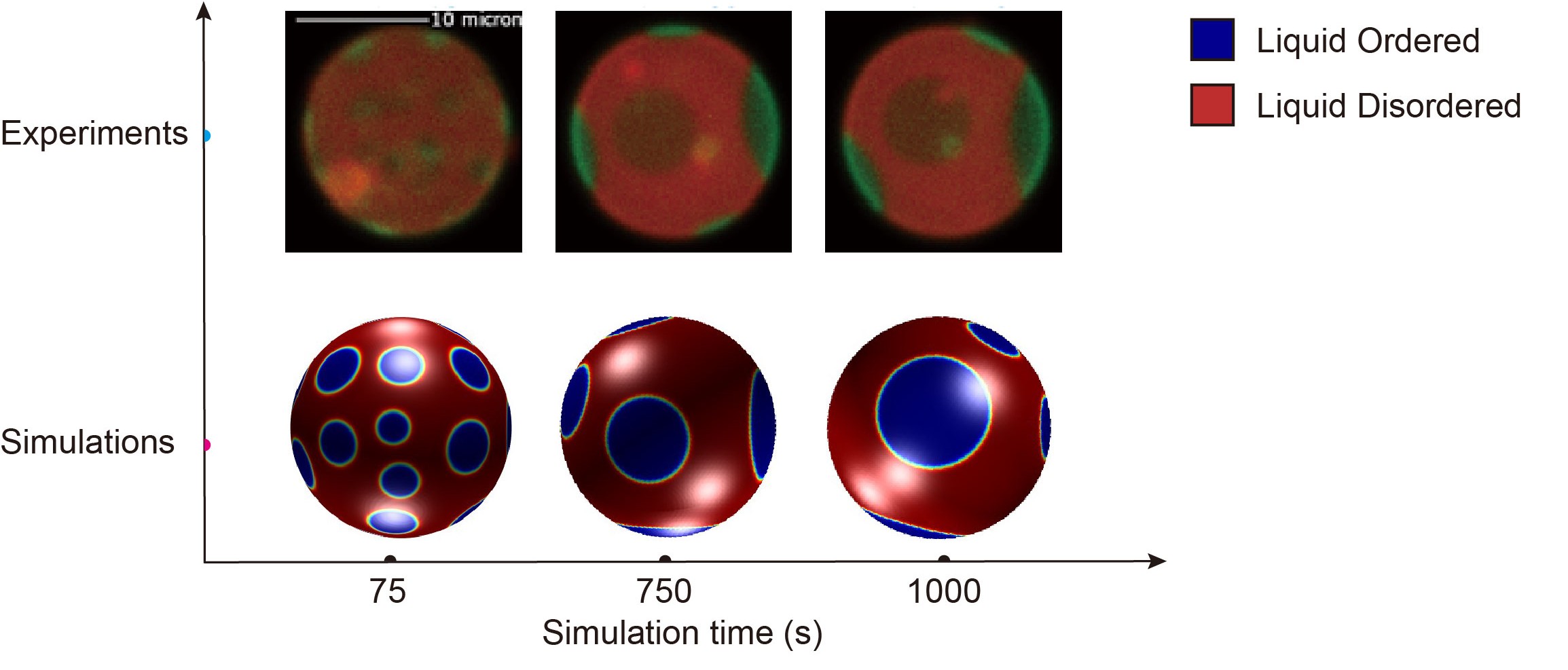}
\caption{Comparison of phase separation dynamics for area fraction $a_{l_d} = 0.7$: experimental images (top) from \cite{wang2022lipid}; numerical simulation results (bottom) illustrating the formation of $l_o$ domains (blue) within the $l_d$ phase (red). See \S~\ref{sec3.3.2} for details.}
\label{fig.73od}
\end{figure}

Figure~\ref{fig.73od} and~\ref{fig.37od} present the phase separation dynamics on the vesicle membrane for $a_{ld} = 0.3$ and $a_{ld} = 0.7$, respectively, alongside the corresponding experimental images from \cite{wang2022lipid}. Domains of the minority lipid component with the same color form, merge, and coarsen into fewer but larger raft-like domains, as shown in the rightmost panels of each figure. Overall, the numerical results exhibit good qualitative agreement with the experiments.

\section{Memory-efficient implementations for solving Navier–Stokes equations}
\label{sec4}
In this section, we extend the cell-centered matrix-free FAS multigrid solver to staggered grids and develop memory-efficient implementations of the projection schemes for the Navier–Stokes equations. The proposed methods are systematically compared with classical first- and second-order schemes in terms of GPU memory footprint, execution efficiency, and convergence behavior. By introducing dynamic reuse of GPU-resident variables, the memory footprint is significantly reduced without sacrificing accuracy. As a result, simulations on grids up to $512^3$ in three dimensions can be carried out on a single RTX~4090 GPU (24~GB) for both first- and second-order schemes. Finally, the effectiveness of the approach is demonstrated through a benchmark simulation of air–water two-bubble coalescence.
\subsection{Multigrid solver on 2D staggered grids}\label{sec4.1}

We first extend the cell-centered matrix-free FAS multigrid solver to staggered grids and apply it to solve the two-dimensional incompressible Navier–Stokes equations:
\begin{align}
\frac{\partial {\bf u}}{\partial t} +{\bf u}\cdot\boldsymbol{\nabla}{\bf u}-\frac{1}{Re}\boldsymbol{\triangle} {\bf u}+\boldsymbol{\nabla}  p= 0 \quad &\text{in } \Omega \times J, \label{Ns_1}\\
\boldsymbol{\nabla} \cdot {\bf u} = 0  \quad &\text{in } \Omega \times J, \label{Ns_2}\\
{\bf u} = {\bf g} \quad &\text{on } \partial\Omega \times J, \label{Ns_3}
\end{align}

Here, $\Omega = \prod_{l=1}^{2} [a_l, b_l]$ denotes the spatial domain, $J = [0, T]$ the time interval, and $Re$ the Reynolds number. In two dimensions, the velocity field is $\mathbf{u} = (u, v)$, while $p$ denotes the pressure. We likewise use the two-grid notation in \S~\ref{sec2.1} to introduce the variables. As illustrated in Figure~\ref{grid}, the $u$-component is located at the east–west (EW) edge-centered grid points of the fine grid $\Omega_1$, i.e., $u_{1,i+1/2,j} \in \bar{\mathcal{E}}^{EW}{\Omega_1}$, while the $v$-component is located at the north–south (NS) edge-centered points, i.e., $v_{1,i,j+1/2} \in \bar{\mathcal{E}}^{NS}{\Omega_1}$. The pressure is defined at the cell-centered grid points, i.e., $p_{1,i,j} \in \bar{\mathcal{C}}_{\Omega_1}$. This arrangement naturally leads to a staggered grid formulation for the numerical solution of the Navier–Stokes equations. Details of the GPU matrix-free, memory-efficient implementations of first- and second-order projection schemes for these equations are given in \S~\ref{sec4.2}-\ref{sec4.3}.

The multigrid solver for edge-centered variables is similar to that used for cell-centered variables. For clarity, we focus here on the operators for EW edge-centered variables, including operator $\mathrm{S}$,  $\mathrm{R}$, and  $\mathrm{P}$, noting that the treatment in the NS edge-centered variables is entirely analogous.

 Operator $\mathrm{S}$ on edge-centered grids is performed in the same X-MCGS manner as for the cell-centered points in 2D. In particular, points are grouped into four sets according to the parity of $(i,j)$ for $u_{1,i+1/2,j}$, and each set is updated separately according to Equations~\eqref{FCGS-1}–\eqref{FCGS-4}. For the boundary conditions in Equation~\eqref{Ns_3}, the smoothing updates interior edge-centered points and ghost points (the latter via prescribed conditions), while boundary points remain unchanged.

For the operators $\mathrm{R}$ and $\mathrm{P}$, the edge-centered case is inherently more complex than the cell-centered case due to the different variable locations. Their detailed formulations are given in \S\ref{app4}. It is worth emphasizing that restriction and prolongation are essentially matrix operations, which makes them particularly well-suited for GPU-based parallelization, thereby significantly improving computational efficiency and reducing execution time.

In addition, alternative boundary conditions for the velocity components of the Navier–Stokes equations, such as periodic boundaries, can be naturally accommodated within our Multigrid solver by introducing additional ghost points at the edge-centered grid. 
\subsection{Memory-efficient implementation of first-order projection scheme} \label{sec4.2}
We begin with the classical first-order projection scheme~\cite{li2022} for the Navier–Stokes equations~\eqref{Ns_1}–\eqref{Ns_3}:
\begin{scheme}\label{first scheme}
Let $\Delta t$ be the time step size, $\mathbf{u}^n=(u^n,v^n,w^n)$ the velocity field, and $p^n$ the pressure at time step $n$. In the projection method, auxiliary corrections $\tilde{\mathbf{u}}^n=(\tilde{u}^n,\tilde{v}^n,\tilde{w}^n)$ and $\tilde{p}^n$ are introduced. Given $(\cdot)^n$, the next step $({\tilde{\mathbf{u}}}^{,n+1},\mathbf{u}^{,n+1},p^{,n+1})$ is obtained by solving
\begin{align}
\frac{\tilde{\mathbf{u}}^{n+1} - \mathbf{u}^n}{\Delta t}+ \mathbf{u}^n \cdot \boldsymbol{\nabla} \mathbf{u}^n- \frac{1}{Re} \boldsymbol{\Delta} \tilde{\mathbf{u}}^{n+1}+ \boldsymbol{\nabla} p^n &= 0, \label{first-1}\\
\Delta t \boldsymbol{\nabla} \cdot (\boldsymbol{\nabla} \tilde{p}^{n+1})&=\boldsymbol{\nabla}\cdot\tilde{\mathbf{u}}^{n+1},\label{first-2}\\
\mathbf{u}^{n+1}&= \tilde{\mathbf{u}}^{n+1} - \Delta t \boldsymbol{\nabla} \tilde{p}^{n+1},\label{first-3}\\
p^{n+1} &= p^n + \tilde{p}^{n+1}.\label{first-4}
\end{align}

In total, the first-order scheme involves 12 variables: $ u^{n+1}, u^n,\tilde{u}^{n+1},\\ v^{n+1}, v^n,\tilde{v}^{n+1},w^{n+1}, w^n, \tilde{w}^{n+1}, p^{n+1}, p^n, \tilde{p}^{n+1}.$
\end{scheme}
We next describe the classical implementation of Scheme~\ref{first scheme}, with its computational procedure summarized in Table~\ref{first-order-classical}. This implementation follows Scheme~\ref{first scheme} directly and requires 12 GPU-resident variables—$U_\text{new}$, $U_\text{old}$, $\tilde{U}_\text{new}$, $V_\text{new}$, $V_\text{old}$, $\tilde{V}_\text{new}$, $W_\text{new}$, $W_\text{old}$, $\tilde{W}_\text{new}$, $P_\text{new}$, $P_\text{old}$, and $\tilde{P}_\text{new}$—which are updated throughout the computation. Temporary GPU variables created within Algorithm~\ref{two-grid algorithm} are not included in this count.

In table \ref{first-order-classical}, the \textbf{Step} column specifies the sequence of computational operations, while the \textbf{Projection method} column details the corresponding numerical procedure. Specifically, in Steps 1, 3, and 5, $f_u$, $f_v$, and $f_w$ compute the source terms (see Algorithm~\ref{two-grid algorithm}) associated with the velocity components $u$, $v$, and $w$, respectively, while Steps 2, 4, and 6 apply the relaxation updates for these components. The \textbf{Update} column identifies the variables modified at each step. Finally, the \textbf{12 GPU-resident variables} store the current values of the corresponding physical quantities throughout the implementation.

\begin{sidewaystable}
{\renewcommand{\arraystretch}{1.8}
\begin{adjustbox}{max width=\textheight, center}
\begin{tabular}{|c|c|c|c|c|c|c|c|c|c|c|c|c|c|c|}
\hline
\multicolumn{3}{|c|}{\textbf{Classical implementation of first-order scheme}} & \multicolumn{12}{|c|}{\textbf{12 GPU-resident variables }}\\
\hline
\textbf{Step} & \makecell[c]{\textbf{Projection method}}  &\textbf{Update}
&$U_\text{new}$&$U_\text{old}$&$\tilde{U}_\text{new}$ &$V_\text{new}$&$V_\text{old}$&$\tilde{V}_\text{new}$ 
&$W_\text{new}$&$W_\text{old}$&$\tilde{W}_\text{new}$ 
&$P_\text{new}$&$P_\text{old}$&$\tilde{P}_\text{new}$\\ 
\hline
1 &\parbox[c][1.8cm][c]{7cm}{$f_u=-\Delta t\left(u^{n}\frac{\partial u^{n}}{\partial x} + v^{n}\frac{\partial u^{n}}{\partial y}+ w^{n}\frac{\partial u^{n}}{\partial z}\right) $\\$\quad- \Delta t(\nabla p^n)_x+u^n$}& $f_u$
&  & $u^{n}$ & & &$v^{n}$ & & & $w^{n}$ & & & $p^{n}$& \\
\hline
2&$\tilde{u}^{n+1}- \frac{\Delta t}{Re} \Delta \tilde{u}^{n+1}=f_u$&$\tilde{u}^{n+1}$ & & &$\tilde{u}^{n+1}$ & & & & & & & & & \\
\hline
3 &\parbox[c][1.8cm][c]{7cm}{$f_v=-\Delta t\left(u^{n}\frac{\partial v^{n}}{\partial x} + v^{n}\frac{\partial v^{n}}{\partial y}+ w^{n}\frac{\partial v^{n}}{\partial z}\right)$\\$\quad- \Delta t(\nabla p^n)_y+v^n$}& $f_v$
&  & $u^{n}$ & & &$v^{n}$ & & & $w^{n}$ & & & $p^{n}$& \\
\hline
4&$\tilde{v}^{n+1}- \frac{\Delta t}{Re} \Delta \tilde{v}^{n+1}=f_v$&$\tilde{v}^{n+1}$ & & & & & &$\tilde{v}^{n+1}$& & & & & & \\
\hline
5  &\parbox[c][1.8cm][c]{7cm}{$f_w=-\Delta t\left(u^{n}\frac{\partial w^{n}}{\partial x} + v^{n}\frac{\partial w^{n}}{\partial y}+ w^{n}\frac{\partial w^{n}}{\partial z}\right)$\\$\quad- \Delta t(\nabla p^n)_z+w^n$}& $f_w$
&  & $u^{n}$ & & &$v^{n}$ & & & $w^{n}$ & & & $p^{n}$& \\
\hline
6&$\tilde{w}^{n+1}- \frac{\Delta t}{Re} \Delta \tilde{w}^{n+1}=f_w$&$\tilde{w}^{n+1}$ & & & & & & & & &$\tilde{w}^{n+1}$ & & & \\
\hline
7 & $\Delta t\nabla\cdot(\nabla\tilde{p}^{n+1}) =\frac{\partial\tilde{u}^{n+1}}{\partial x}+\frac{\partial\tilde{v}^{n+1}}{\partial y}+\frac{\partial\tilde{w}^{n+1}}{\partial z}$
& $\tilde{p}^{n+1}$& & &$\tilde{u}^{n+1}$ & & & $\tilde{v}^{n+1}$ & & & $\tilde{w}^{n+1}$ & & & $\tilde{p}^{n+1}$ \\
\hline
8 & 
$u^{n+1} =\tilde{u}^{n+1}-\Delta t(\nabla\tilde{p}^{n+1})_x$
&$u^{n+1}$&$u^{n+1}$& &$\tilde{u}^{n+1}$ & & & & & & & & &$\tilde{p}^{n+1}$\\
\hline
9 & 
$v^{n+1} =\tilde{v}^{n+1}-\Delta t(\nabla\tilde{p}^{n+1})_y$
&$v^{n+1}$& & & & $v^{n+1}$& & $\tilde{v}^{n+1}$& & & & & &$\tilde{p}^{n+1}$\\
\hline
10 & 
$w^{n+1} =\tilde{w}^{n+1}-\Delta t(\nabla\tilde{p}^{n+1})_z$
&$w^{n+1}$& & & && & &$w^{n+1}$ & &$\tilde{w}^{n+1}$ & & &$\tilde{p}^{n+1}$\\
\hline
11 & 
$p^{n+1}=p^n+\tilde{p}^{n+1}$& $p^{n+1}$
& & & & & & & & & &$p^{n+1}$ &$p^{n}$&$\tilde{p}^{n+1}$\\
\hline
12 &$u^{n}=u^{n+1}$&$u^{n}$&$u^{n+1}$&$u^{n}$& & & & & & & & & &\\
\hline
13 &$v^{n}=v^{n+1}$&$v^{n}$& & & &$v^{n+1}$&$v^{n}$& & & & & & &\\
\hline
14 &$w^{n}=w^{n+1}$&$w^{n}$& & & & & & &$w^{n+1}$&$w^{n}$& & & &\\
\hline
15 &$p^{n}=p^{n+1}$&$u^{n}$& & & & & & & & & &$p^{n+1}$&$p^{n}$&\\
\hline
\end{tabular}
\end{adjustbox}}
\caption{Classical implementation of first-order projection scheme for the 3D Navier–Stokes equations. See \S~\ref{sec4.2} for details.} 
\label{first-order-classical}
\end{sidewaystable}

Next, we introduce the memory-efficient implementation of the first-order scheme, summarized in Table~\ref{first-order-memory-efficient}. This optimized version reduces the number of GPU-resident variables from 12 to 8 by introducing a dynamic reuse strategy, while fully preserving the structure of the original scheme. A closer inspection of the computational steps~\eqref{first-1}–\eqref{first-4} shows that they are executed independently, thereby allowing a single GPU-resident variable to represent multiple quantities at different stages. The key idea is to eliminate the conventional time-stepping update $(\cdot)_\text{old} = (\cdot)_\text{new}$ (see Steps 12–15 in Table~\ref{first-order-classical}). Instead, by separating the source term $f$ and the operator term from Algorithm~\ref{two-grid algorithm}, $u^{n+1}$ can be stored directly in $U_\text{old}$ and $\tilde{u}^{n+1}$ in $U_\text{new}$ without affecting the numerical results (see Step 8 in Table~\ref{first-order-memory-efficient}). The same reuse strategy applies to $v^{n+1}, \tilde{v}^{n+1}, w^{n+1}, \tilde{w}^{n+1}, p^{n+1}$, and $\tilde{p}^{n+1}$ (see Steps 9–11 in Table~\ref{first-order-memory-efficient}).

\begin{sidewaystable}
{\renewcommand{\arraystretch}{1.8}
\begin{adjustbox}{max width=\textheight, center}
\begin{tabular}{|c|c|c|c|c|c|c|c|c|c|c|}
\hline
\multicolumn{3}{|c|}{\textbf{ Memory-efficient implementation of first-order scheme}} & \multicolumn{8}{|c|}{\textbf{8 GPU-resident variables}}\\
\hline
\textbf{Step} & \makecell[c]{\textbf{Projection method}}  &\textbf{Update}
& $U_\text{new}$ & $U_\text{old} $& $V_\text{new}$ & $V_\text{old}$ & $W_\text{new}$ & $W_\text{old}$ & $P_\text{new}$ & $P_\text{old} $ \\ 
\hline
1 &\parbox[c][1.8cm][c]{7cm}{$f_u=-\Delta t\left(u^{n}\frac{\partial u^{n}}{\partial x} + v^{n}\frac{\partial u^{n}}{\partial y}+ w^{n}\frac{\partial u^{n}}{\partial z}\right) $\\$\quad- \Delta t(\nabla p^n)_x+u^n$}
& $f_u$&  & $u^{n}$ &  &$v^{n}$ & & $w^{n}$ & & $p^{n}$ \\
\hline
2&$\tilde{u}^{n+1}- \frac{\Delta t}{Re} \Delta \tilde{u}^{n+1}=f_u$&$\tilde{u}^{n+1}$ &$\tilde{u}^{n+1}$ &&&&&&&\\
\hline
3 &\parbox[c][1.8cm][c]{7cm}{$f_v=-\Delta t\left(u^{n}\frac{\partial v^{n}}{\partial x} + v^{n}\frac{\partial v^{n}}{\partial y}+ w^{n}\frac{\partial v^{n}}{\partial z}\right)$\\$\quad- \Delta t(\nabla p^n)_y+v^n$}& $f_v$
&  & $u^{n}$ &  &$v^{n}$ & & $w^{n}$ & & $p^{n}$ \\
\hline
4&$\tilde{v}^{n+1}- \frac{\Delta t}{Re} \Delta \tilde{v}^{n+1}=f_v$&$\tilde{v}^{n+1}$ &$\tilde{v}^{n+1}$ &&&&&&&\\
\hline
5  &\parbox[c][1.8cm][c]{7cm}{$f_w=-\Delta t\left(u^{n}\frac{\partial w^{n}}{\partial x} + v^{n}\frac{\partial w^{n}}{\partial y}+ w^{n}\frac{\partial w^{n}}{\partial z}\right)$\\$\quad- \Delta t(\nabla p^n)_z+w^n$}& $f_w$
&  & $u^{n}$ &  &$v^{n}$ & & $w^{n}$ & & $p^{n}$ \\
\hline
6&$\tilde{w}^{n+1}- \frac{\Delta t}{Re} \Delta \tilde{w}^{n+1}=f_w$&$\tilde{w}^{n+1}$ &$\tilde{w}^{n+1}$ &&&&&&&\\
\hline
7 & 
$\Delta t\nabla\cdot(\nabla\tilde{p}^{n+1}) =\frac{\partial\tilde{u}^{n+1}}{\partial x}+\frac{\partial\tilde{v}^{n+1}}{\partial y}+\frac{\partial\tilde{w}^{n+1}}{\partial z}$
& $\tilde{p}^{n+1}$& $\tilde{u}^{n+1}$ & & $\tilde{v}^{n+1}$ && $\tilde{w}^{n+1}$ & & $\tilde{p}^{n+1}$& \\
\hline
8 & 
$u^{n+1} =\tilde{u}^{n+1}-\Delta t(\nabla\tilde{p}^{n+1})_x$
& $u^{n+1}$& $\tilde{u}^{n+1}$ & $u^{n+1}$& & & & & $\tilde{p}^{n+1}$ &\\
\hline
9 & 
$v^{n+1} =\tilde{v}^{n+1}-\Delta t(\nabla\tilde{p}^{n+1})_y$&$v^{n+1}$ &
& & $\tilde{v}^{n+1}$ &$v^{n+1}$ & & & $\tilde{p}^{n+1}$ &\\
\hline
10 & 
$w^{n+1} =\tilde{w}^{n+1}-\Delta t(\nabla\tilde{p}^{n+1})_z$& $w^{n+1}$& & &
& & $\tilde{w}^{n+1}$ &$w^{n+1}$ & $\tilde{p}^{n+1}$ &\\
\hline
11 & 
$p^{n+1}=p^n+\tilde{p}^{n+1}$& $p^{n+1}$
& & & & & & &$\tilde{p}^{n+1}$ & $p^{n+1}, p^{n}$\\
\hline
\end{tabular}
\end{adjustbox}}
\caption{Memory-efficient implementation of first-order projection scheme for the 3D Navier–Stokes equations suitable for GPU implementation. See \S~\ref{sec4.2} for details.}
\label{first-order-memory-efficient}
\end{sidewaystable}

It is worth noting that, in the first-order scheme, a third-order WENO discretization~\cite{liu1994,jiang1996,shu1999} is employed to stabilize the nonlinear convection term $\mathbf{u}^n \cdot \nabla \mathbf{u}^n$. The standard WENO procedure involves numerous conditional statements and logical branching to select smooth stencils, compute nonlinear weights, and handle boundary stencils, particularly near domain boundaries or in regions with steep gradients. To better leverage GPU acceleration, the conditional logic in the WENO computation is replaced by matrix operations, and the use of intermediate variables is minimized. The same strategy is adopted in the second-order scheme, as discussed in the next subsection.

The GPU memory footprint, execution performance, and convergence tests for both the classical and memory-efficient implementations of the first-order scheme are presented in \S~\ref{sec4.4}.

\subsection{Memory-efficient implementation of second-order projection
scheme}\label{sec4.3}
We next introduce the second-order projection schemes~\cite{li2022} for the Navier–Stokes Equations~\eqref{Ns_1}–\eqref{Ns_3}:
\begin{scheme}\label{second scheme}
Assume $(\cdot)^n$ is known and $(\cdot)^{-1} = (\cdot)^0$. The solution at the next time step $({\bf \tilde{u}}^{n+1},{\bf u}^{n+1},p^{n+1})$ is obtained by solving
\begin{align}
\frac{\tilde{\mathbf{u}}^{n+1} - \mathbf{u}^n}{\Delta t} + \mathbf{u}^{n+\frac{1}{2}} \cdot \boldsymbol{\nabla} \mathbf{u}^{n+\frac{1}{2}} - \frac{1}{Re} \boldsymbol{\Delta} \tilde{\mathbf{u}}^{n+\frac{1}{2}} + \boldsymbol{\nabla} p^n &= 0, \label{second-1}\\
\Delta t \boldsymbol{\nabla} \cdot (\boldsymbol{\nabla} \tilde{p}^{n+1}) &= \boldsymbol{\nabla} \cdot \tilde{\mathbf{u}}^{n+1}, \label{second-2}\\
\mathbf{u}^{n+1} &= \tilde{\mathbf{u}}^{n+1} - \Delta t \boldsymbol{\nabla} \tilde{p}^{n+1}, \label{second-3}\\
p^{n+1} &= p^n + \tilde{p}^{n+1}, \label{second-4}
\end{align}

where the notations are consistent with those used in Scheme~\ref{first scheme}.  The interpolated values $\mathbf{u}^{n+1/2}$ and $\tilde{\mathbf{u}}^{n+1/2}$ are defined as
$\mathbf{u}^{n+1/2} = (3\mathbf{u}^n - \mathbf{u}^{n-1})/2$, $\tilde{\mathbf{u}}^{n+1/2} = (\mathbf{u}^n + \tilde{\mathbf{u}}^{n+1})/2$. In total, the second-order projection scheme involves 15 variables: $u^{n+1}, u^n, u^{n-1}, \tilde{u}^{n+1}, v^{n+1}, v^n, v^{n-1}, \tilde{v}^{n+1}, w^{n+1}, w^n, w^{n-1}, \tilde{w}^{n+1}, p^{n+1}, p^n, \tilde{p}^{n+1}$.
\end{scheme}
The classical implementation (Table~\ref{second-order-classical}) requires 15 GPU-resident variables, including three additional ones ($U_{\text{2old}}, V_{\text{2old}}, W_{\text{2old}}$) to store $u^{n-1}, v^{n-1}, w^{n-1}$.

\begin{sidewaystable}
{\renewcommand{\arraystretch}{1.8}
\begin{adjustbox}{max width=\textheight, center}
\begin{tabular}{|c|c|c|c|c|c|c|c|c|c|c|c|c|c|c|c|c|c|}
\hline
\multicolumn{3}{|c|}{\textbf{Classical implementation of second-order  scheme}} & \multicolumn{15}{|c|}{\textbf{15 GPU-resident variables}}\\
\hline
\textbf{Step} & \makecell[c]{\textbf{Projection method}}  &\textbf{Update}
&$U_\text{new}$&$U_\text{old}$&$U_\text{2old}$&$\tilde{U}_\text{new}$ &$V_\text{new}$&$V_\text{old}$&$V_\text{2old}$&$\tilde{V}_\text{new}$ 
&$W_\text{new}$&$W_\text{old}$&$W_\text{2old}$&$\tilde{W}_\text{new}$ 
&$P_\text{new}$&$P_\text{old}$&$\tilde{P}_\text{new}$\\ 
\hline
1 &  \parbox[c][2.2cm][c]{9.5cm}{$f_u=-\Delta t\left(\frac{3u^{n}-u^{n-1}}{2}\frac{\partial(3u^{n}-u^{n-1})}{2\partial x} + \frac{3v^{n}-v^{n-1}}{2}\frac{\partial(3u^{n}-u^{n-1})}{2\partial y}\right.$\\
$\left.+ \frac{3w^{n}-w^{n-1}}{2}\frac{\partial(3u^{n}-u^{n-1})}{2\partial z}\right)+ \frac{\Delta t}{2Re} \Delta u^n - \Delta t(\nabla p^n)_x+u^n$}
& $f_u$&  & $u^{n}$ &$u^{n-1}$& & &$v^{n}$ &$v^{n-1}$& & & $w^{n}$ &$w^{n-1}$ & & & $p^{n}$& \\
\hline
2&$\tilde{u}^{n+1}- \frac{\Delta t}{2Re} \Delta \tilde{u}^{n+1}=f_u$&$\tilde{u}^{n+1}$ & & & & $\tilde{u}^{n+1}$& & & & & & & & & & &\\
\hline
3 &  \parbox[c][2.2cm][c]{9.5cm}{$f_v=-\Delta t\left(\frac{3u^{n}-u^{n-1}}{2}\frac{\partial(3u^{v}-v^{n-1})}{2\partial x} + \frac{3v^{n}-v^{n-1}}{2}\frac{\partial(3v^{n}-v^{n-1})}{2\partial y}\right.$\\
$\left.+ \frac{3w^{n}-w^{n-1}}{2}\frac{\partial(3v^{n}-v^{n-1})}{2\partial z}\right)+ \frac{\Delta t}{2Re} \Delta v^n - \Delta t(\nabla p^n)_y+v^n$}
& $f_v$&  & $u^{n}$ &$u^{n-1}$& & &$v^{n}$ &$v^{n-1}$& & & $w^{n}$ &$w^{n-1}$ & & & $p^{n}$& \\
\hline
4&$\tilde{v}^{n+1}- \frac{\Delta t}{2Re} \Delta \tilde{v}^{n+1}=f_v$&$\tilde{v}^{n+1}$ & & & & 
& & & & $\tilde{v}^{n+1}$& & & & & & &\\
\hline
5 & \parbox[c][2.2cm][c]{9.5cm}{$f_w=-\Delta t\left(\frac{3u^{n}-u^{n-1}}{2}\frac{\partial(3w^{n}-w^{n-1})}{2\partial x} + \frac{3v^{n}-v^{n-1}}{2}\frac{\partial(3w^{n}-w^{n-1})}{2\partial y}\right.$\\
$\left.+ \frac{3w^{n}-w^{n-1}}{2}\frac{\partial(3w^{n}-w^{n-1})}{2\partial z}\right)+ \frac{\Delta t}{2Re} \Delta w^n - \Delta t(\nabla p^n)_z+w^n$}
& $f_w$&  & $u^{n}$ &$u^{n-1}$& & &$v^{n}$ &$v^{n-1}$& & & $w^{n}$ &$w^{n-1}$ & & & $p^{n}$& \\
\hline
6&$\tilde{w}^{n+1}- \frac{\Delta t}{2Re} \Delta \tilde{w}^{n+1}=f_w$&$\tilde{w}^{n+1}$ & & & & 
& & & & & & & &$\tilde{w}^{n+1}$ & & &\\
\hline
7 & $\Delta t\nabla\cdot(\nabla\tilde{p}^{n+1}) =\frac{\partial\tilde{u}^{n+1}}{\partial x}+\frac{\partial\tilde{v}^{n+1}}{\partial y}+\frac{\partial\tilde{w}^{n+1}}{\partial z}$
& $\tilde{p}^{n+1}$& & & &$\tilde{u}^{n+1}$ & & & & $\tilde{v}^{n+1}$ & & & & $\tilde{w}^{n+1}$ & & & $\tilde{p}^{n+1}$ \\
\hline
8 & 
$u^{n+1} =\tilde{u}^{n+1}-\Delta t(\nabla\tilde{p}^{n+1})_x$
&$u^{n+1}$&$u^{n+1}$& & &$\tilde{u}^{n+1}$ & & & & & & & & & & &$\tilde{p}^{n+1}$\\
\hline
9 & 
$v^{n+1} =\tilde{v}^{n+1}-\Delta t(\nabla\tilde{p}^{n+1})_y$
&$v^{n+1}$& & & & & $v^{n+1}$& & & $\tilde{v}^{n+1}$& & & & & & &$\tilde{p}^{n+1}$\\
\hline
10 & 
$w^{n+1} =\tilde{w}^{n+1}-\Delta t(\nabla\tilde{p}^{n+1})_z$
&$w^{n+1}$& & & & & & & & &$w^{n+1}$ & & &$\tilde{w}^{n+1}$ & & &$\tilde{p}^{n+1}$\\
\hline
11 & 
$p^{n+1}=p^n+\tilde{p}^{n+1}$& $p^{n+1}$
& & & & & & & & & & & & &$p^{n+1}$ &$p^{n}$&$\tilde{p}^{n+1}$\\
\hline
12 &$u^{n-1}=u^{n}$, $u^{n}=u^{n+1}$&$u^{n-1},u^{n}$&$u^{n+1}$&$u^{n}$&$u^{n-1}$& & & & & & & & & & & &\\
\hline
13 &$v^{n-1}=v^{n}$, $v^{n}=v^{n+1}$&$v^{n-1},v^{n}$& & & & &$v^{n+1}$&$v^{n}$&$v^{n-1}$& & & & & & & &\\
\hline
14 &$w^{n-1}=w^{n}$, $w^{n}=w^{n+1}$&$w^{n-1},w^{n}$& & & & & & & & &$w^{n+1}$&$w^{n}$&$w^{n-1}$& & & &\\
\hline
15 &$p^{n}=p^{n+1}$&$p^{n}$& & & & & & & & & & & & &$p^{n+1}$&$p^{n}$&\\
\hline
\end{tabular}
\end{adjustbox}}
\caption{Classical implementation of second-order projection scheme for the 3D Navier–Stokes equations. See \S~\ref{sec4.3} for details.} 
\label{second-order-classical}
\end{sidewaystable}

We then propose a memory-efficient implementation (Table~\ref{second-order-memory-efficient}), which reduces the number of GPU-resident variables from 15 to 8 while preserving second-order accuracy. Unlike the first-order case, the time-stepping update $(\cdot)\text{old} = (\cdot)\text{new}$ is retained but deferred until the source terms $f$ for $(u,v,w)$ at the next step have been defined (see Steps 2, 5, and 8). This lagged update allows $u^{n-1}$, $u^{n}$, and $\tilde{u}^{n+1}$ to be cyclically assigned to only two GPU-resident variables $U_\text{new}$ and $U_\text{old}$, ensuring that $u^{n+1/2}$ is always avaiable when constructing WENO terms. For example, $u^{n+1/2} = (3u^{n} - u^{n-1})/2$ is used in Step 1, while $u^{n+1/2} = (u^{n} + \tilde{u}^{n+1})/2$ is employed in Steps 4 and 7. The same procedure applies to $v^{n+1/2}$ and $w^{n+1/2}$, while for pressure, the update $P_\text{old} = P_\text{new}$ is again omitted (see Step 14).

\begin{sidewaystable}
{\renewcommand{\arraystretch}{1.8}
\begin{adjustbox}{max width=\textheight, center}
\begin{tabular}{|c|c|c|c|c|c|c|c|c|c|c|}
\hline
\multicolumn{3}{|c|}{\textbf{Memory-efficient implementation of second-order scheme}} & \multicolumn{8}{|c|}{\textbf{8 GPU-resident variables}}\\
\hline
\textbf{Step} & \makecell[c]{\textbf{Projection method}}  &\textbf{Update}
& $U_\text{new}$ & $U_\text{old} $& $V_\text{new}$ & $V_\text{old}$ & $W_\text{new}$ & $W_\text{old}$ & $P_\text{new}$ & $P_\text{old} $ \\ 
\hline
1 &\parbox[c][2.2cm][c]{9.5cm}{$f_u=-\Delta t\left(\frac{3u^{n}-u^{n-1}}{2}\frac{\partial(3u^{n}-u^{n-1})}{2\partial x} + \frac{3v^{n}-v^{n-1}}{2}\frac{\partial(3u^{n}-u^{n-1})}{2\partial y}\right.$\\$\left.+ \frac{3w^{n}-w^{n-1}}{2}\frac{\partial(3u^{n}-u^{n-1})}{2\partial z}\right)+ \frac{\Delta t}{2Re} \Delta u^n - \Delta t(\nabla p^n)_x+u^n$}& $f_u$
& $u^{n}$ & $u^{n-1}$ &$v^{n}$ & $v^{n-1}$ &$w^{n}$ & $w^{n-1}$ & & $p^{n}$ \\
\hline
2&$u^{n-1}=u^{n}$&$u^{n-1}$&$u^{n}$& $u^{n-1}$&&&&&&\\
\hline
3&$\tilde{u}^{n+1}- \frac{\Delta t}{2Re} \Delta \tilde{u}^{n+1}=f_u$&$\tilde{u}^{n+1}$ &$\tilde{u}^{n+1}$ &&&&&&&\\
\hline
4 & \parbox[c][2.2cm][c]{9.5cm}{$f_v=-\Delta t\left(\frac{u^{n}+\tilde{u}^{n+1}}{2}\frac{\partial(3v^{n}-v^{n-1})}{2\partial x} + \frac{3v^{n}-v^{n-1}}{2}\frac{\partial(3v^{n}-v^{n-1})}{2\partial y}\right.$\\$\left. +\frac{3w^{n}-w^{n-1}}{2}\frac{\partial(3v^{n}-v^{n-1})}{2\partial z}\right)+ \frac{\Delta t}{2Re} \Delta v^n - \Delta t(\nabla p^n)_y+v^n$}& $f_v$
& $\tilde{u}^{n+1}$ & $u^{n}$ &$v^{n}$ & $v^{n-1}$ &$w^{n}$ & $w^{n-1}$&& $p^{n}$\\
\hline
5&$v^{n-1}=v^{n}$&$v^{n-1}$&$v^{n}$& $v^{n-1}$&&&&&&\\
\hline
6&$\tilde{v}^{n+1}- \frac{\Delta t}{2Re} \Delta \tilde{v}^{n+1}=f_v$&$\tilde{v}^{n+1}$ &$\tilde{v}^{n+1}$ &&&&&&&\\
\hline
7 & \parbox[c][2.2cm][c]{9.5cm}{$f_w=-\Delta t\left(\frac{u^{n}+\tilde{u}^{n+1}}{2}\frac{\partial(3w^{n}-w^{n-1})}{2\partial x} + \frac{v^{n}+\tilde{v}^{n+1}}{2}\frac{\partial(3w^{n}-w^{n-1})}{2\partial y}\right.$\\$\left.+ \frac{3w^{n}-w^{n-1}}{2}\frac{\partial(3w^{n}-w^{n-1})}{2\partial z}\right)+ \frac{\Delta t}{2Re} \Delta w^n - \Delta t(\nabla p^n)_z+w^n$}& $f_w$
& $\tilde{u}^{n+1}$ & $u^{n}$ &$\tilde{v}^{n+1}$ & $v^{n}$ &$w^{n}$ & $w^{n-1}$&& $p^{n}$\\
\hline
8&$w^{n-1}=w^{n}$&$w^{n-1}$&$w^{n}$& $w^{n-1}$&&&&&&\\
\hline
9&$\tilde{w}^{n+1}- \frac{\Delta t}{2Re} \Delta \tilde{w}^{n+1}=f_w$&$\tilde{w}^{n+1}$ &$\tilde{w}^{n+1}$ &&&&&&&\\
\hline
10 & 
$\Delta t\nabla\cdot(\nabla\tilde{p}^{n+1}) =\frac{\partial\tilde{u}^{n+1}}{\partial x}+\frac{\partial\tilde{v}^{n+1}}{\partial y}+\frac{\partial\tilde{w}^{n+1}}{\partial z}$
& $\tilde{p}^{n+1}$& $\tilde{u}^{n+1}$ & & $\tilde{v}^{n+1}$ && $\tilde{w}^{n+1}$ & & $\tilde{p}^{n+1}$& \\
\hline
11 & 
$u^{n+1} =\tilde{u}^{n+1}-\Delta t(\nabla\tilde{p}^{n+1})_x$
& $u^{n+1}$& $\tilde{u}^{n+1},u^{n+1}$ & & & & & & $\tilde{p}^{n+1}$ &\\
\hline
12 & 
$v^{n+1} =\tilde{v}^{n+1}-\Delta t(\nabla\tilde{p}^{n+1})_y$& $v^{n+1}$&
& & $\tilde{v}^{n+1},v^{n+1}$ & & & & $\tilde{p}^{n+1}$ &\\
\hline
13 & 
$w^{n+1} =\tilde{w}^{n+1}-\Delta t(\nabla\tilde{p}^{n+1})_z$&$w^{n+1}$ & & &
& & $\tilde{w}^{n+1},w^{n+1}$ & & $\tilde{p}^{n+1}$ &\\
\hline
14 & 
$p^{n+1}=p^n+\tilde{p}^{n+1}$& $p^{n+1}$
& & & & & & &$\tilde{p}^{n+1}$ & $p^{n+1}, p^{n}$\\
\hline
\end{tabular}
\end{adjustbox}}
\caption{Memory-efficient implementation of second-order projection scheme for the 3D Navier–Stokes equations suitable for GPU implementation. See \S~\ref{sec4.3} for details.}
\label{second-order-memory-efficient}
\end{sidewaystable}

\subsection{Memory-efficient vs. Classical implementation: GPU memory, performance, and convergence}
\label{sec4.4}
 We compare the memory-efficient and classical implementations of the first- and second-order schemes in terms of GPU memory footprint and execution time. Across devices, the memory-efficient variants reduce both memory usage and wall-clock time by approximately 20–30\%. To further assess accuracy and stability, we perform: (a) an asymptotic-in-time convergence test following Li~\cite{li2022}, confirming first- and second-order temporal accuracy; and (b) a benchmark comparison with the lid-driven cavity results of Ghia et al.~\cite{ghia1982} over multiple Reynolds numbers, demonstrating reliability and robustness.
 
\subsubsection{GPU memory footprint}\label{sec4.4.1}
\begin{figure}[H] 
\centering
\includegraphics[width=1\linewidth]{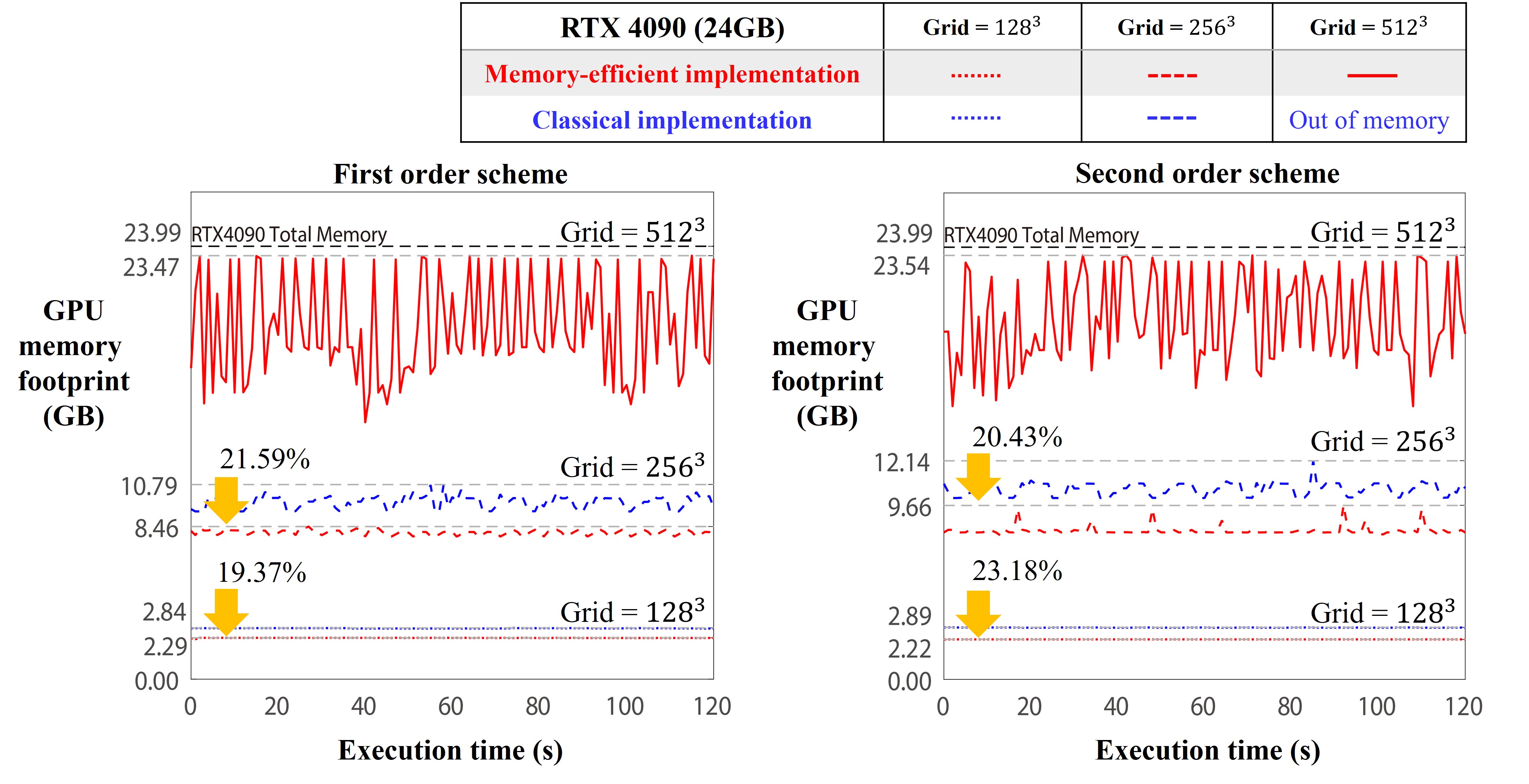}
\caption{GPU memory footprint versus execution time on RTX~4090 (24~GB). Left: first-order; right: second-order. The memory-efficient implementations successfully perform the $512^3$ case within 24~GB memory limit. For smaller grids, peak memory is reduced by $\approx$19-23\% compared with the classical implementations. See \S~\ref{sec4.4.1} for details.}
\label{memory-efficient4090}
\end{figure}

GPU memory footprint was measured on a single RTX~4090 (24~GB) for both orders using Ghia’s cavity setup (details below) and compared with the classical implementations. Figure~\ref{memory-efficient4090} shows the memory profile during execution (red: memory-efficient; blue: classical). With the memory-efficient approach, the 3D $512^3$ case fits within 24~GB, whereas the classical versions fail due to excessive resident variables. For $256^3$ and $128^3$, the memory-efficient first- and second-order schemes reduce peak footprint by about 20–25\% relative to their classical counterparts. Notably, for a fixed grid size, the first- and second-order memory-efficient variants exhibit nearly identical memory usage.

\begin{figure}[H] 
\centering
\includegraphics[width=1\linewidth]{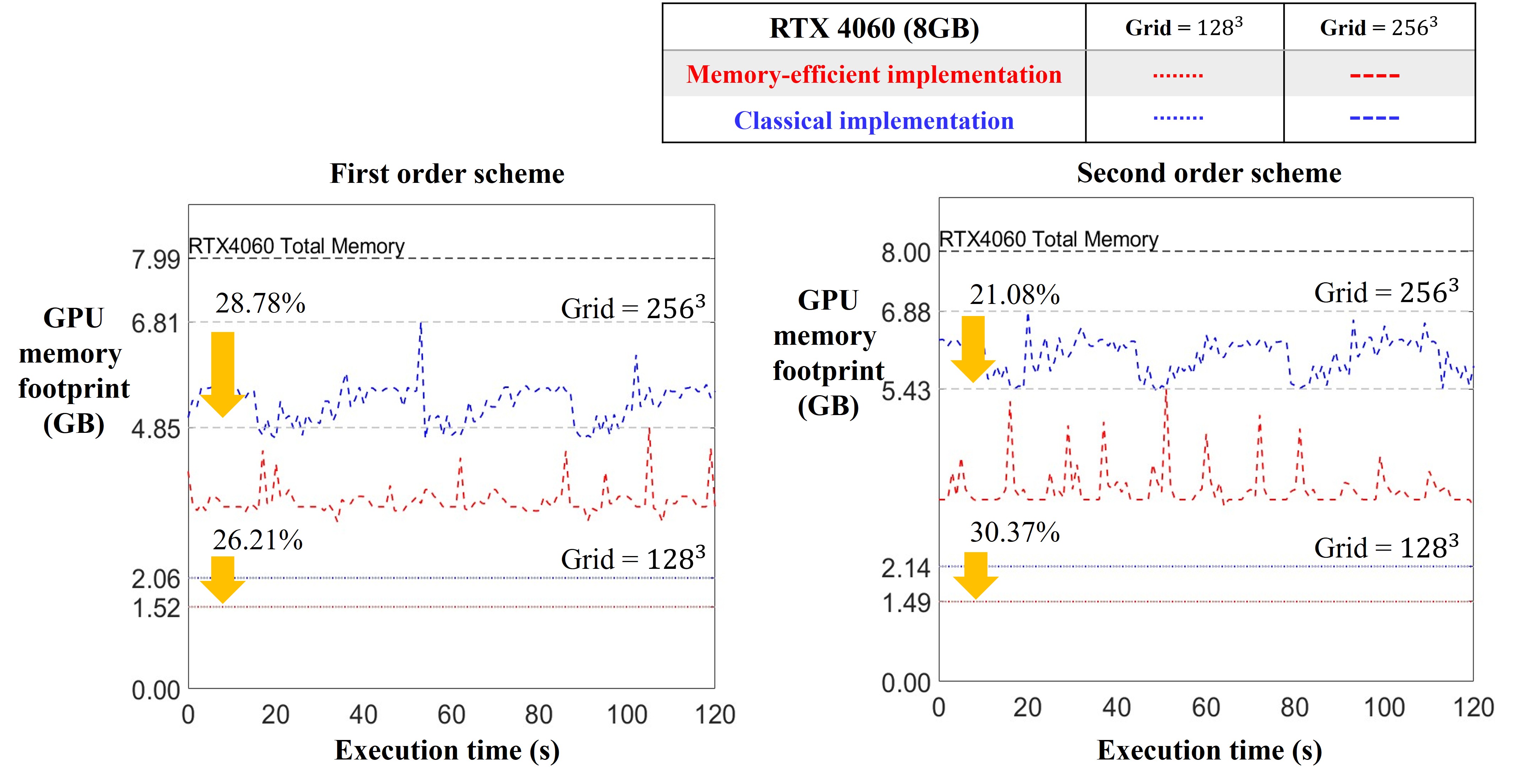}
\caption{GPU memory footprint versus execution time on RTX 4060 (8 GB). Left: first-order; right: second-order. All implementations run up to the $256^3$ grid within 8~GB memory limit, while the memory-efficient versions reduce peak memory usage by about 21–30$\%$ compared with the classical ones. See \S~\ref{sec4.4.1} for details.}
\label{memory-efficient4060}
\end{figure}

Because memory allocation strategies differ across GPUs, we also evaluate the memory–time trade-off on an RTX~4060 (8~GB), see Figure~\ref{memory-efficient4060}. Relative to the classical implementations, the proposed method achieves even larger savings, with the second-order scheme reducing memory by up to $\approx$30\% on the $128^3$ grid.

\begin{rem}
Compared with Figure~\ref{memory-efficient4090}, Figure~\ref{memory-efficient4060} indicates that, for the same grid size, MATLAB reduces memory footprint automatically on lower-memory GPUs, at the expense of slower execution. This underscores the intrinsic trade-off between memory consumption and speed.
\end{rem}
\subsubsection{Execution performance}\label{sec4.4.2}
\begin{figure}[H] 
\centering
\includegraphics[width=1\linewidth]{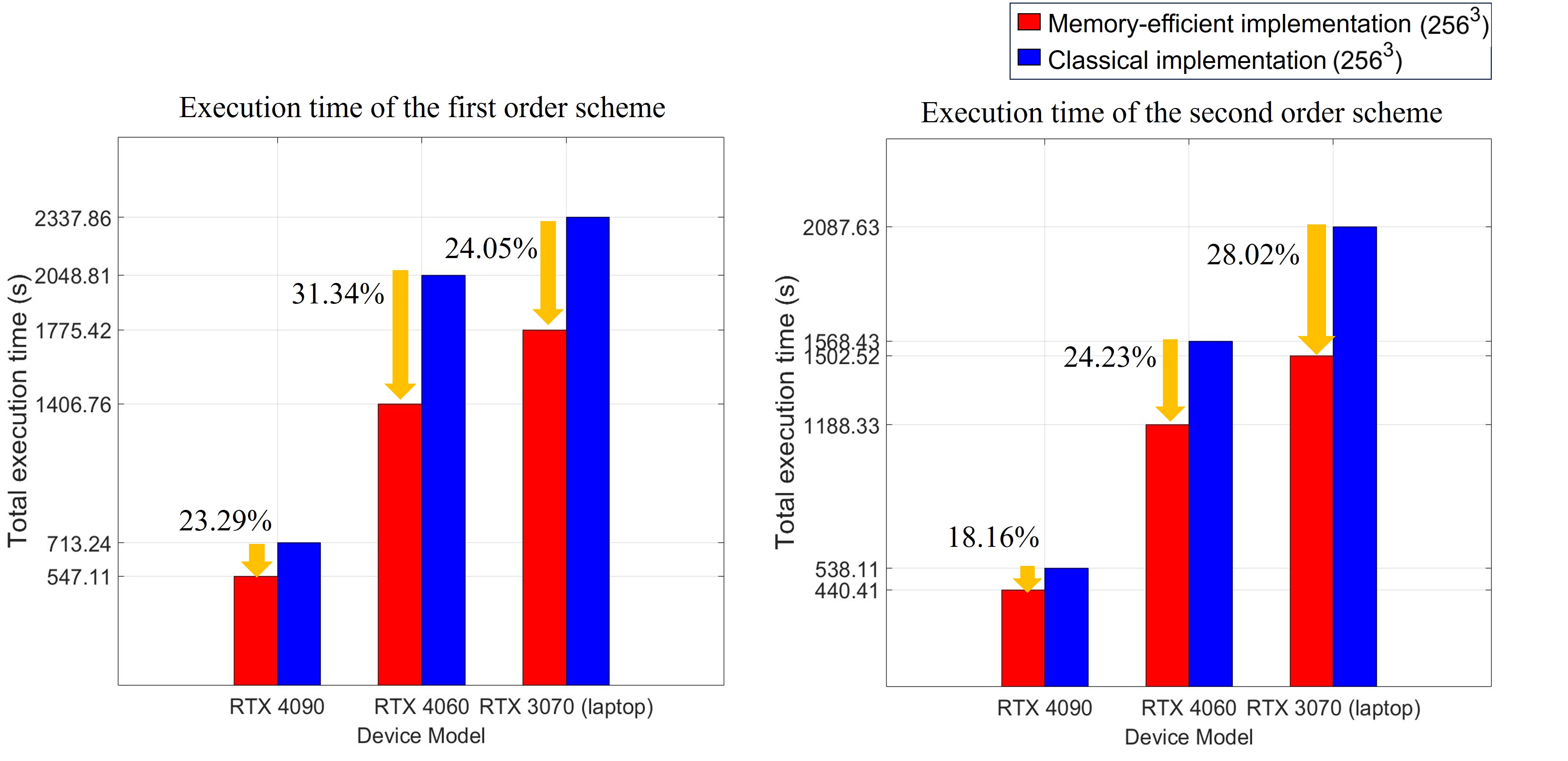}
\caption{Total execution time (from $t=0$ to $t=1$ with $\Delta t=0.01$) for first- and second-order schemes on RTX~4090 (24~GB), RTX~4060 (8~GB), and RTX~3070 (laptop~8~GB) at a $256^3$ grid. The memory-efficient implementations reduce execution time by $\approx$18–31\% relative to the classical ones. See \S~\ref{sec4.4.2} for details.}
\label{memory-efficient-sudu}
\end{figure}

Since execution time on GPUs is often correlated with memory footprint, we benchmark execution time for a fixed grid $256^3$ with $\Delta t=0.01$, advancing from $t=0$ to $t=1$. Figure~\ref{memory-efficient-sudu} reports total execution time on RTX~4090 (24~GB), RTX~4060 (8~GB), and RTX~3070 (laptop~8~GB) for both orders and both implementations. The memory-efficient variants consistently outperform the classical ones, yielding $\approx$ 18–31\% speedups depending on the device and scheme.
% Next, we examine the convergence behavior to ensure that our memory-efficient schemes maintain the expected level of accuracy.

\subsubsection{Asymptotic convergence test}\label{sec4.4.3}

We next carry out an asymptotic convergence test following Li~\cite{li2022}, designed to verify whether the schemes achieve first- and second-order accuracy, respectively. To construct an exact solution of \eqref{Ns_1}–\eqref{Ns_3}, we prescribe the source terms accordingly, yielding
\[
\begin{cases}
p_{1,i,j,n}^{\text{exact}}&= \sin(t_n)\left(\sin(\pi y_j) - \frac{2}{\pi}\right), \\
u_{1,i+\frac{1}{2},j,n}^{\text{exact}}& = \sin(t_n)\sin^2(\pi x_{i+\frac{1}{2}})\sin(2\pi y_j), \\
v_{1,i,j+\frac{1}{2},n}^{\text{exact}} &= -\sin(t_n)\sin(2\pi x_{i})\sin^2(\pi y_{j+\frac{1}{2}}).
\end{cases}
\]
Homogeneous Neumann boundary conditions are imposed on the pressure field, while Dirichlet boundary conditions are applied to the velocity components.

To minimize spatial discretization errors, the test is conducted in the 2D domain $\Omega=[0,1]^2$ using a fixed fine grid of $N_1=M_1=4096$, with time interval $J=[0,1]$ and Reynolds number $Re=10$. With $t_n=n\Delta t$, the user-defined parameters in Table~\ref{tab:FASsymbols} are set to $tol=10^{-9}$, $k_{\text{Max}}=20$, $s=2$, and $meshLevel=11$. The asymptotic-in-time convergence results for the memory-efficient first- and second-order schemes are reported in Tables~\ref{first-convergence-test} and \ref{second-convergence-test}, demonstrating that the expected first- and second-order accuracies are indeed achieved.

\begin{table}[H]
\tabcolsep=0.25cm
\centering%把表居中
\begin{tabular}{ccccccc}%四个c代表该表一共四列，内容全部居中
\\
\toprule%第一道横线
$\Delta t$&Error $(u)$&Order&Error $(v)$&Order&Error $(p)$&Order\\
\midrule%第二道横线 
$\frac{1}{10}$&$3.82E-03$&-&$4.53E-03$&-&$3.18E-02$&-\\
$\frac{1}{20}$&$1.57E-03$&$1.28$&$1.68E-03$&$1.43$&$1.52E-02$&$1.06$\\
$\frac{1}{40}$&$7.53E-04$&$1.06$&$7.70E-04$&$1.13$&$7.49E-03$&$1.02$\\
$\frac{1}{80}$&$3.73E-04$&$1.01$&$3.76E-04$&$1.03$&$3.73E-03$&$1.00$\\
$\frac{1}{160}$&$1.86E-04$&$1.00$&$1.87E-04$&$1.00$&$1.87E-03$&$1.00$\\
\bottomrule%第三道横线
\end{tabular}
\caption{Asymptotic-in-time convergence tests of the memory-efficient first-order scheme on a 2D $4096^2$ grid with time steps $\Delta t$ ranging from $1/10$ to $1/160$. The error is measured as $ \| *_{1,i,j,n}^{\mathrm{exact}} - *_{1,i,j}^{n} \|_{L^2(\Omega_1)}$ at $t_n=1$ for all variables. The results demonstrate first-order convergence rate, consistent with Li~\cite{li2022}. See \S~\ref{sec4.4.3} for details. }
\label{first-convergence-test}
\end{table}

\begin{table}[H]

\tabcolsep=0.25cm
\centering%把表居中
\begin{tabular}{ccccccc}%四个c代表该表一共四列，内容全部居中
\\
\toprule%第一道横线
$\Delta t$&Error $(u)$&Order&Error $(v)$&Order&Error $(p)$&Order\\
\midrule%第二道横线 
$\frac{1}{10}$&$1.24E-03$&-&$1.59E-03$&-&$2.06E-02$&-\\
$\frac{1}{20}$&$2.90E-03$&$2.09$&$3.70E-04$&$2.11$&$7.43E-02$&$1.47$\\
$\frac{1}{40}$&$7.07E-05$&$2.04$&$8.9E-05$&$2.05$&$2.92E-03$&$1.35$\\
$\frac{1}{80}$&$1.75E-05$&$2.02$&$2.20E-05$&$2.02$&$1.25E-03$&$1.22$\\
$\frac{1}{160}$&$4.35E-06$&$2.01$&$5.48E-06$&$2.01$&$5.73E-03$&$1.13$\\
\bottomrule%第三道横线
\end{tabular}
\caption{Asymptotic-in-time convergence tests of the memory-efficient second-order scheme on a 2D $4096^2$ grid with time steps $\Delta t$ ranging from $1/10$ to $1/160$. The error is measured as $ \| *_{1,i,j,n}^{\mathrm{exact}} - *_{1,i,j}^{n} \|_{L^2(\Omega_1)}$ at $t_n=1$ for all variables. The results exhibit second-order convergence for the velocity components $(u,v)$ and first-order convergence for the pressure $p$, consistent with Li~\cite{li2022}. See \S~\ref{sec4.4.3} for details. }
\label{second-convergence-test}
\end{table}
\subsubsection{Benchmark validations}\label{sec4.4.4}
The classical 2D and 3D lid-driven cavity flow problems, first introduced by Ghia et al.~\cite{ghia1982}, have become standard benchmarks for validating Navier–Stokes solvers at various Reynolds numbers (see, e.g., \cite{ku1987,shi2020}). In the 2D case, the computational domain is $\Omega = [0,1]^2$, where the top lid at $y=1$ moves uniformly in the $x$-direction while the remaining three walls are stationary and impermeable. In the 3D case, the domain is $\Omega = [0,1]^3$, with the top lid at $z=1$ translating uniformly in the $x$-direction, and the other five walls fixed with no-slip, impermeable boundary conditions. The initial velocity is set to zero, corresponding to a quiescent fluid.

Simulations are performed using Algorithm~\ref{two-grid algorithm} with the user-defined parameters in Table~\ref{tab:FASsymbols}: $tol = 10^{-10}$, $k_{\text{Max}} = 20$, and $s = 2$. For 2D, computations are carried out on a $4096^2$ grid with $meshLevel = 11$, and for 3D on a $512^3$ grid with $meshLevel = 8$. In both cases, the time step is set to $\Delta t = 0.001$.

\begin{figure}[H] 
\centering
\includegraphics[width=1\linewidth]{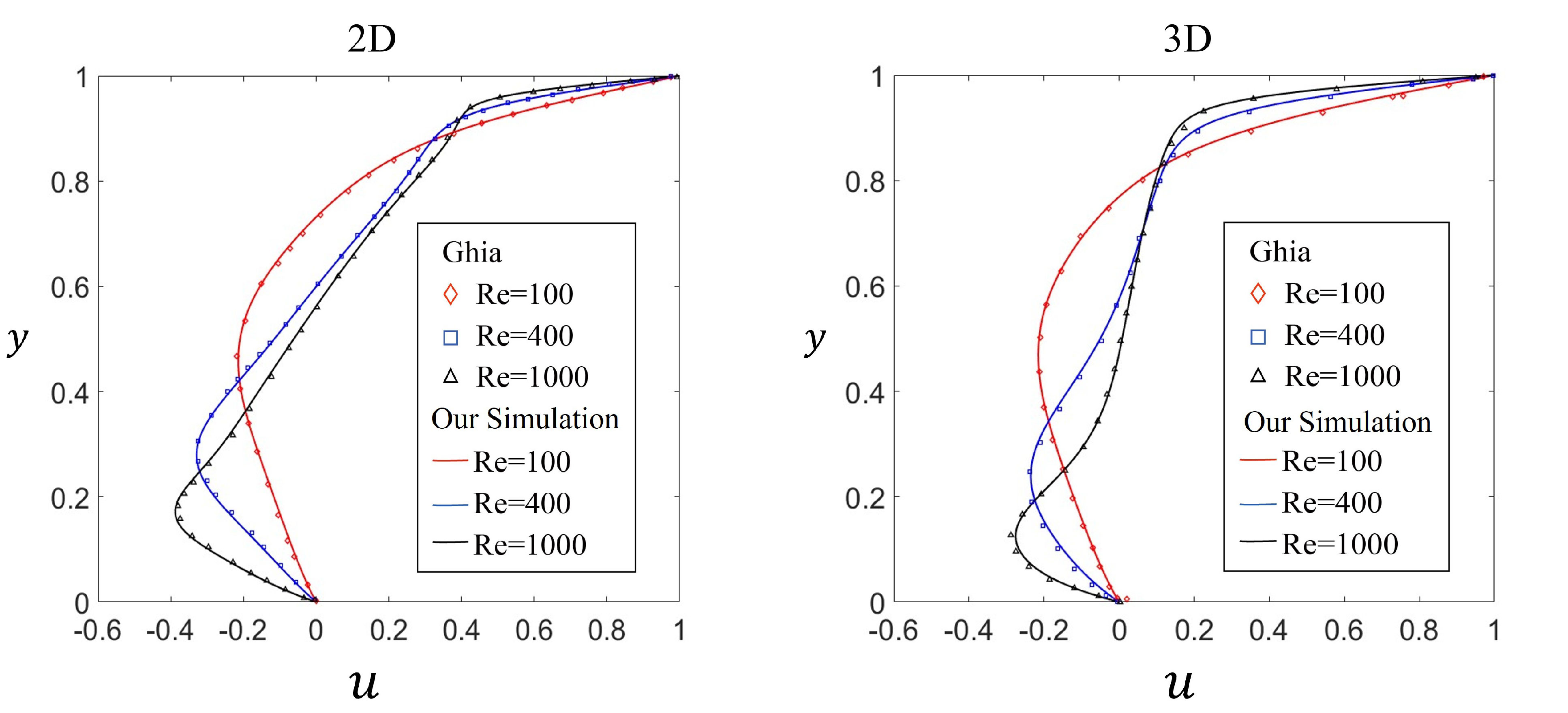}
\caption{Steady-state solutions of the memory-efficient implementations (lines) compared with the benchmark data of Ghia et al.~\cite{ghia1982} (symbols) for the 2D ($4096^2$, left) and 3D ($512^3$, right) lid-driven cavity flowsat Re = 100, 400, and 1000. Since the steady-state results of the first- and second-order schemes are nearly identical, only the second-order ones are shown. See \S~\ref{sec4.4.4} for details.}
\label{GHIA-test}
\end{figure}

Figure~\ref{GHIA-test} compares our numerical results with the benchmark data of Ghia et al.~\cite{ghia1982}. The benchmark data are shown as symbols, and our simulations are plotted as lines: red for $Re=100$, blue for $Re=400$, and black for $Re=1000$. Steady-state solutions are extracted after sufficiently long integration times. Since the steady-state profiles obtained from the first- and second-order schemes are nearly identical, only the second-order memory-efficient results are shown. As evident in the figure, the numerical solutions exhibit excellent agreement with the benchmark data, confirming the accuracy and reliability of the proposed implementation.

\begin{figure}[H] 
\centering
\includegraphics[width=1\linewidth]{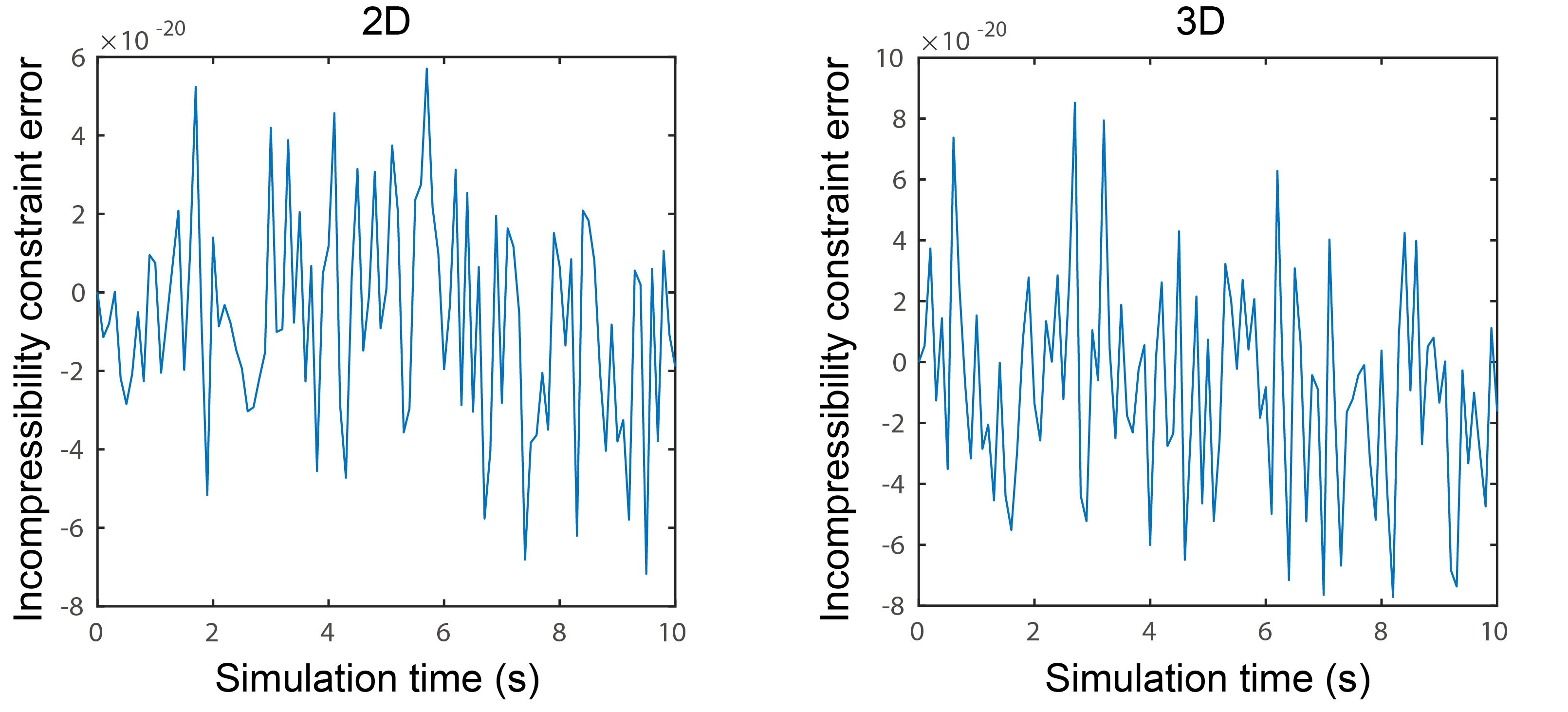}
\caption{Evolution of the incompressibility constraint error over simulation time. Left: 2D cavity flow; Right: 3D cavity flow. The errors remain below $10^{-19}$ at all time steps, indicating that the incompressibility condition is preserved to machine precision.}
\label{div-u}
\end{figure}

To verify that our solver faithfully preserves the incompressibility constraint error, we also measure the discrete divergence of the velocity field in the benchmark cavity flow problem on the finest grid ($level=L$) with spacing $h_{L}$:

\begin{align}
\text{2D:} \quad \int_{\Omega} \boldsymbol{\nabla} \cdot \mathbf{u}^{n+1}~d{\boldsymbol{x}} &\approx h_{L}^2\sum_{i,j} (\boldsymbol{\nabla} \cdot \mathbf{u}^{n+1})_{L,i,j},\quad \text{for}~\forall n\geq 0,\label{div2D}\\
\text{3D:} \quad \int_{\Omega} \boldsymbol{\nabla} \cdot \mathbf{u}^{n+1}~d{\boldsymbol{x}} &\approx h_{L}^3\sum_{i,j,l} (\boldsymbol{\nabla} \cdot \mathbf{u}^{n+1})_{L,i,j,l} ,\quad \text{for}~\forall n\geq 0.\label{div3D}
\end{align}
Figure~\ref{div-u} shows the evolution of the discrete divergence of the velocity field over the simulation time, computed according to Equations~\eqref{div2D}-\eqref{div3D}. As can be seen, in both the 2D and 3D cases, the values remain below $10^{-19}$ at all time steps, well within the double-precision round-off errors. This confirms that the incompressibility condition is preserved to machine precision, demonstrating the accuracy and robustness of our solver.

\subsection{ Air-in-water two-bubble coalescence problem}\label{sec4.5}
 To demonstrate the practical utility of our memory-efficient implementation, we reproduce the experimental observations of two air bubbles coalescing in water from \cite{brereton1991}. The simulations are based on the quasi-incompressible Navier–Stokes Cahn–Hilliard (q-NSCH) model \cite{guo2022, denefle2014} for two-phase flows with variable density and viscosity:
\begin{align}
\rho \frac{\partial \mathbf{u}}{\partial t} + \rho \mathbf{u} \cdot \boldsymbol{\nabla} \mathbf{u} + \boldsymbol{\nabla} p &= \boldsymbol{\nabla} \cdot \left( \nu(c) \boldsymbol{\nabla} \mathbf{u} \right) + \frac{1}{3} \boldsymbol{\nabla} \left( \nu(c) \boldsymbol{\nabla} \cdot \mathbf{u} \right) - \tilde{\sigma} c \boldsymbol{\nabla} \mu - \rho g \mathbf{z},  \\
\boldsymbol{\nabla} \cdot \mathbf{u} &= \alpha \boldsymbol{\nabla} \cdot \left( \tilde{\sigma} M(c) \boldsymbol{\nabla} \mu \right) + \alpha^2 \boldsymbol{\nabla} \cdot \left( \tilde{\sigma} M(c) \boldsymbol{\nabla} p \right),  \\
\frac{\partial c}{\partial t} + \boldsymbol{\nabla} \cdot (\mathbf{u} c) &= \boldsymbol{\nabla} \cdot \left( \tilde{\sigma} M(c) \boldsymbol{\nabla} \mu \right) + \alpha \boldsymbol{\nabla} \cdot \left( \tilde{\sigma} M(c) \boldsymbol{\nabla} p \right),  \\
\mu &= f(c) - \epsilon^2\Delta c. 
\end{align}
Here, $\mathbf{u}$ denotes the velocity, $p$ the pressure, and $c \in [0,1]$ the phase variable, with $c = 1$ and $c = 0$ corresponding to air and water, respectively. $\mu$ denotes the chemical potential, $\epsilon$ the diffuse interface thickness, $g$ the gravitational acceleration, and $\mathbf{z}$ the unit vector in the vertical direction. The effective surface tension $\tilde{\sigma}$ is related to the physical surface tension $\sigma$ by $\tilde{\sigma} = 6\sqrt{2}\, \sigma$. The density and viscosity vary with $c$ as
\begin{align}
\rho(c) &= \rho_1 c + \rho_2 (1 - c), \\
\nu(c) &= \nu_1 c + \nu_2 (1 - c),
\end{align}
where $\rho_i$ and $\nu_i$ ($i = 1,2$) are the constant density and viscosity of fluid $i$. 
The relative density ratio is defined as $\alpha = (\rho_2 - \rho_1)/\rho_2$. The nonlinear term $F(c) = c^2(1 - c)^2/4$ is the standard double-well potential, with $f(c) = F'(c)$. 
The degenerate mobility is given by $M(c) = c^2(1 - c)^2 + \epsilon$.

In the simulation, two initially spherical bubbles of radius $R=0.005$ are placed in an oblique alignment within the domain $\Omega = [0,0.04] \times [0,0.04] \times [0,0.16]$, discretized by a uniform $256 \times 256 \times 1024$ grid and solved on a single RTX~4090 (24~GB). No-slip boundary conditions are imposed for $\mathbf{u}$ on all domain boundaries, while $c$ and $\mu$ satisfy no-flux conditions.

The time step is $\Delta t = 2.5 \times 10^{-5}$. The numerical parameters are set as $ \varepsilon = 0.0003, \; g = 9.8, \; \rho_1 = 1, \; \rho_2 = 1000, \; \nu_1 = 2.15 \times 10^{-4}, \; \nu_2 = 2.15 \times 10^{-2}, \; \sigma = 0.0052 $. 

\begin{figure}[H] 
    \centering
    \includegraphics[width=1\linewidth]{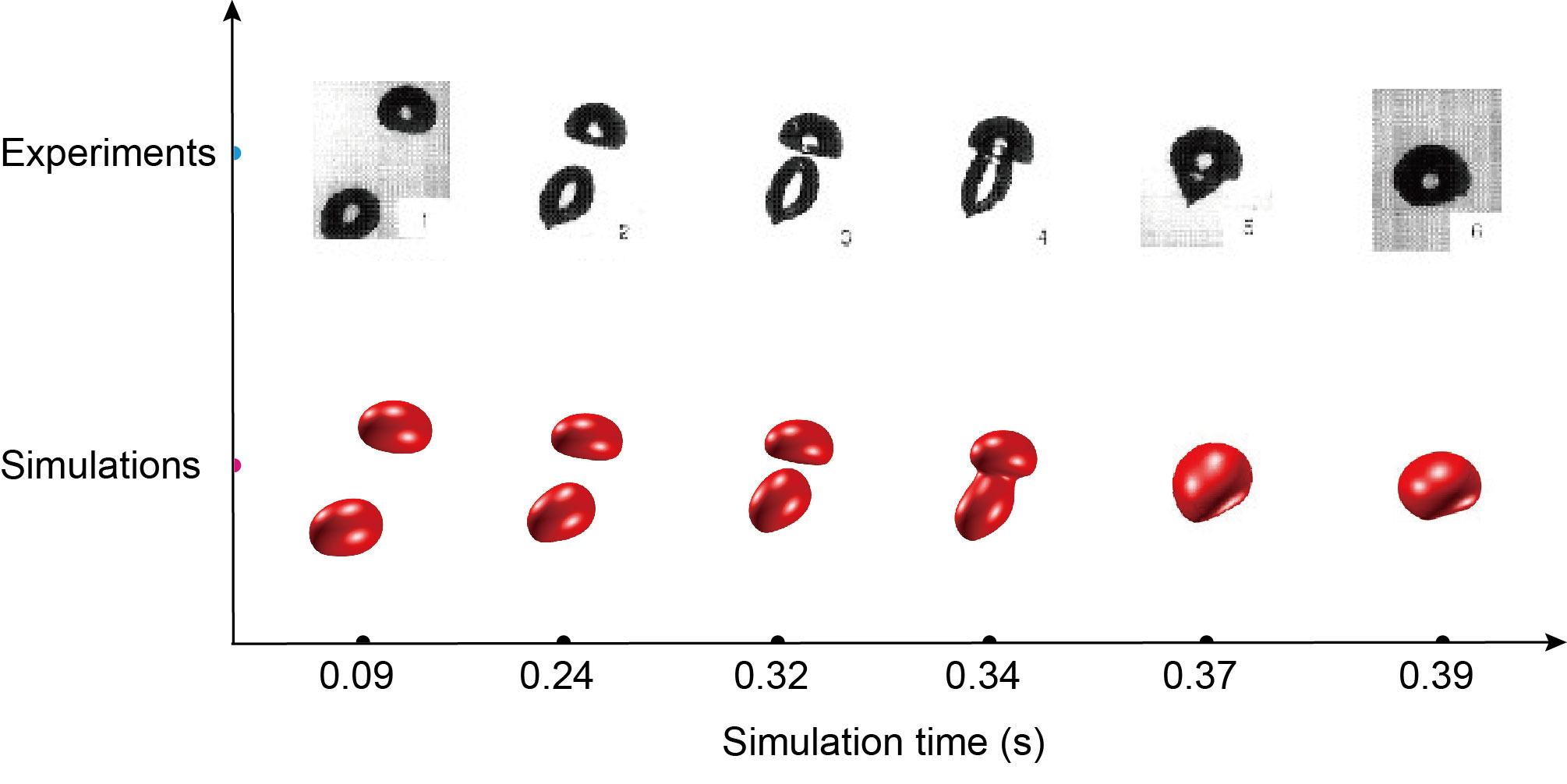}
    \caption{Comparison between experimental observations~\cite{brereton1991} and numerical predictions (iso-surface $c=0.5$) for the oblique coalescence of two bubbles, computed on a $256 \times 256 \times 1024$ grid. See \S~\ref{sec4.5} for details. }
    \label{bubble3D}
\end{figure}
The numerical results obtained with the memory-efficient implementations (Figure~\ref{bubble3D}) show excellent agreement with the experimental data of Brereton~\cite{brereton1991}. In particular, the lower bubble deforms and accelerates in the $z$-direction as it enters the wake of the upper bubble, a phenomenon attributable to the smaller drag force acting on it prior to coalescence. These findings not only reproduce the experimental dynamics but also validate the physical fidelity of the q-NSCH model and the effectiveness of the proposed GPU-based solver for large-scale two-phase flow simulations.

\begin{figure}[H] 
    \centering
    \includegraphics[scale=0.55]{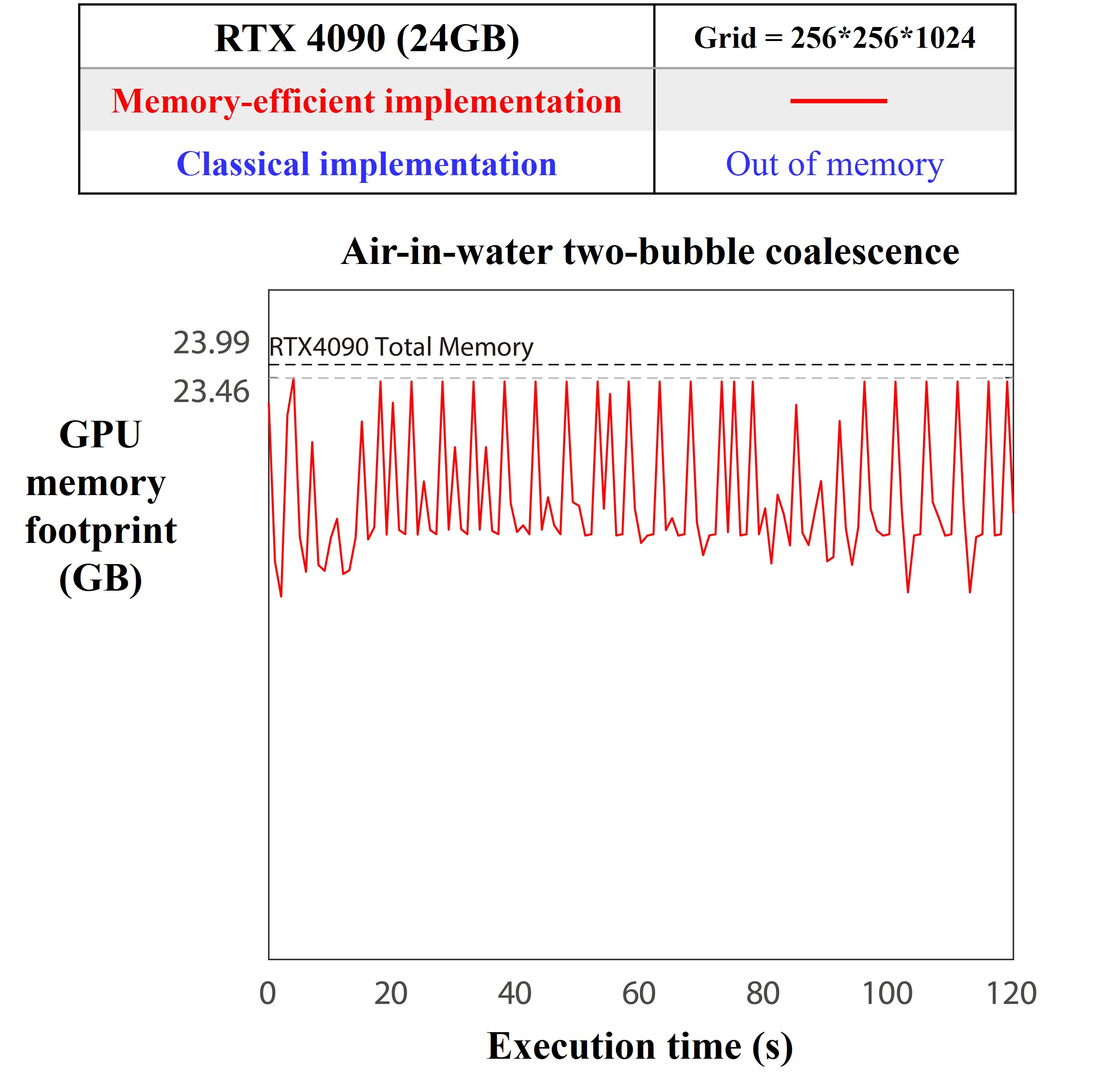}
    \caption{GPU memory footprint versus execution time on RTX~4090 (24~GB) for Air-in-water two-bubble coalescence problem. The memory-efficient implementations of second-order scheme successfully perform the $256\times256\times1024$ case within 24~GB memory limit, while the classical counterpart exceeds available memory capacity. See \S~\ref{sec4.4} for details.}
    \label{bubble-memory}
\end{figure}

Finally, Figure~\ref{bubble-memory} shows the memory footprint of our second-order memory-efficient implementation for the q-NSCH equations on a single RTX 4090 (24 GB) using a $256\times256\times1024$ grid. The results demonstrate that the memory-efficient implementation can handle this large-scale simulation within the GPU memory limit, while the classical implementation exceeds available memory.

\section{Conclusion}
\label{sec5}
In this work, we have developed a matrix-free Full Approximation Storage (FAS) multigrid solver based on staggered finite differences and implemented in MATLAB on GPUs. The matrix-free formulation together with intermediate-variable reuse substantially reduces the memory footprint, enabling large-scale 2D and 3D simulations. To further optimize GPU performance, we propose an X-shaped Multi-Color Gauss–Seidel (X-MCGS) smoother, in which the grid is partitioned into four disjoint sets in 2D (or eight in 3D) following an X-shape ordering. This decomposition eliminates explicit loops and conditional branching (e.g., \texttt{if}-statements), yielding more regular memory access, avoiding thread divergence, and enabling fully vectorized updates. Although each iteration requires more sub-steps than RBGS, the improved memory regularity and parallel efficiency offset this overhead, making X-MCGS particularly effective for large-scale GPU-accelerated multigrid solvers. Moreover, we show that the X-shape scheme requires only about half as many iterations as the U-shape ordering and one-third as many as the Z-shape ordering to achieve the same convergence. Additional GPU acceleration is applied to restriction and prolongation operations, further enhancing overall efficiency.

Algebraic and asymptotic convergence tests confirm the solver’s robustness and accuracy. Benchmark results demonstrate substantial acceleration: on an $8192^2$ grid in 2D, the RTX~4090, RTX~4060, and RTX~3070 achieve speedups of about 61×, 14×, and 13×, respectively, over a single CPU core, while in 3D on a $512^3$ grid the RTX~4090 reaches 47×; for a smaller $256^3$ grid, the RTX~4060 and RTX~3070 deliver 17× and 12×, respectively. These results highlight both the efficiency of the solver and its consistent scalability across GPU platforms.Applications to grain growth and multicomponent vesicle phase separation further validate the solver’s capability for complex, large-scale systems, with results in good agreement with theoretical scaling laws and experimental observations.

To further alleviate GPU memory limitations and enhance performance, we develop memory-efficient implementations of both first- and second-order projection schemes for the incompressible Navier–Stokes equations. By employing a dynamic reuse strategy, the number of GPU-resident variables is reduced from 12 (first-order) and 15 (second-order) to only eight, lowering the memory footprint and improving execution efficiency by 20–30\%. This optimization enables large-scale simulations up to $512^3$ on a single NVIDIA RTX~4090 GPU (24~GB), where classical implementations would exceed device memory. The successful reproduction of air–water two-bubble coalescence experiments on a $256 \times 256 \times 1024$ grid highlights both the physical fidelity of the models and the practical effectiveness of the solver.

Overall, the proposed framework provides a scalable and memory-efficient tool for large-scale multiphysics simulations on modern GPUs, laying a solid foundation for future extensions to more complex flow problems, including large-scale fluid dynamics and multiphase flow equations, all solvable efficiently on a single-GPU platform. Moreover, the framework also paves the way for extensions to multi-GPU cooperative computing, enabling even larger-scale and higher-dimensional simulations with enhanced adaptability and computational potential.

\section*{Acknowledgements}
J.L., S.T. and Z.G. acknowledge support by the National Natural Science Foundation of China (Grant No. 12371387), and S.M. acknowledges partial support by the National Science Foundation of the USA through grants NSF-DMS 2309547.

\appendix
\section{Red-Black Gauss-Seidel smoothing operator in 2D}
\label{app1}
\begin{figure}[H] 
    \centering
    \includegraphics[scale=0.55]{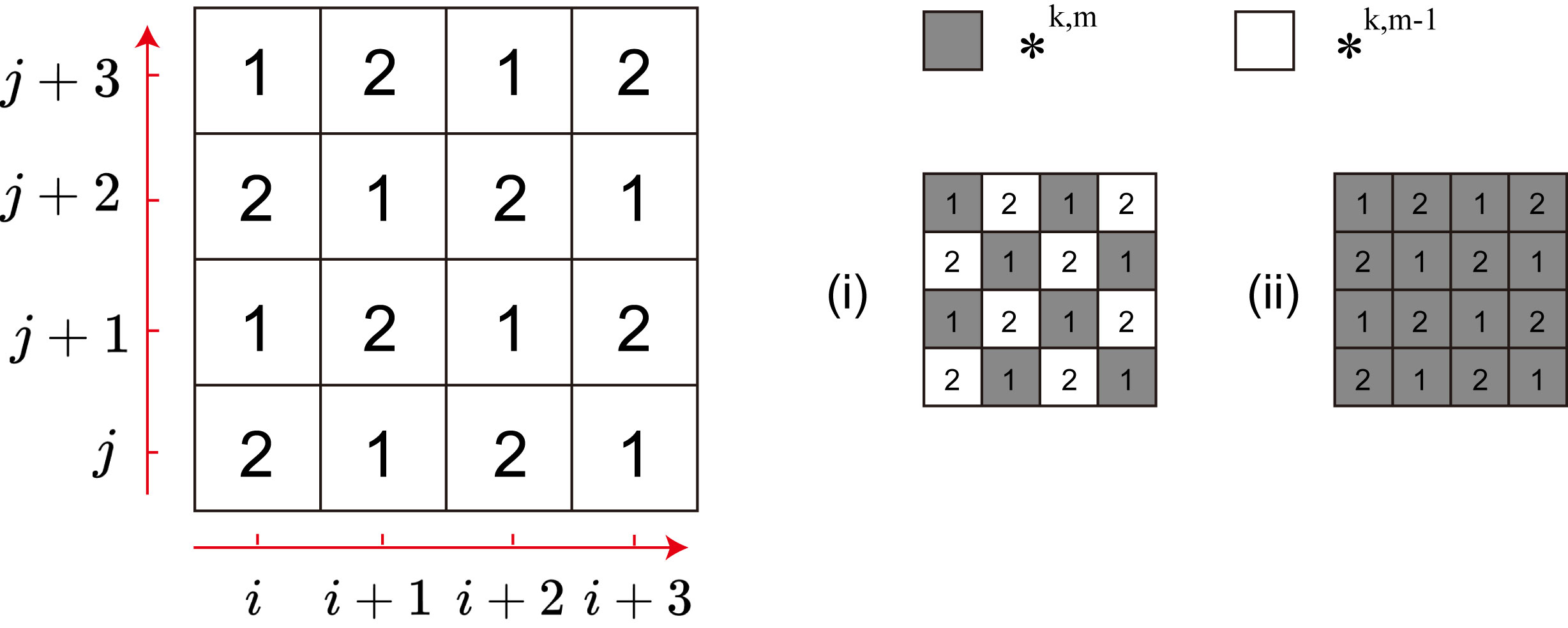}
    \caption{Illustration of the RBGS partitioning. \textbf{Left:} The computational grid is decomposed into two disjoint sets (1--2) based on the parity of the row and column indices. 
    \textbf{Right:} Sequential update order: (i) set~1, (ii) set~2. White cells correspond to the values $*^{k,m-1}$ from the preceding smoothing step, whereas gray cells denote the new values $*^{k,m}$ updated from them during the current $m$-th smoothing of the $k$-th V-cycle. See \S\ref{app1} for details.}
    \label{RBGS}
\end{figure}

We briefly review the Red–Black Gauss–Seidel (RBGS) method \cite{wesseling1992,briggs2000}. Although RBGS exhibits satisfactory convergence, its sequential nature limits parallelization in MATLAB. Grid points are grouped into black and red sets according to the parity of $i+j$, where $i$ and $j$ are the indices of $p_{1,i,j}$ (see Figure~\ref{RBGS}). At each grid level, all black points (set 1) are updated simultaneously using the latest red values, followed by an update of the red points (set 2).

All updates on the fine grid are performed in a matrix-free manner based on the discrete Equation~\eqref{PDE_discret}. The corresponding update rules for each color set are specified below. Specifically, in the $k$-th multigrid V-cycle, once $p_{1,i,j}^{k,m-1}$ has been computed after the $(m-1)$-th smoothing step, the new value $p_{1,i,j}^{k,m}$ is obtained by:

\begin{enumerate}[label=(\roman*)]
    \item Parallel update of set 1:
    \begin{equation}\label{BWGS-B}
    \begin{minipage}{0.15\linewidth}
    \includegraphics[width=\linewidth]{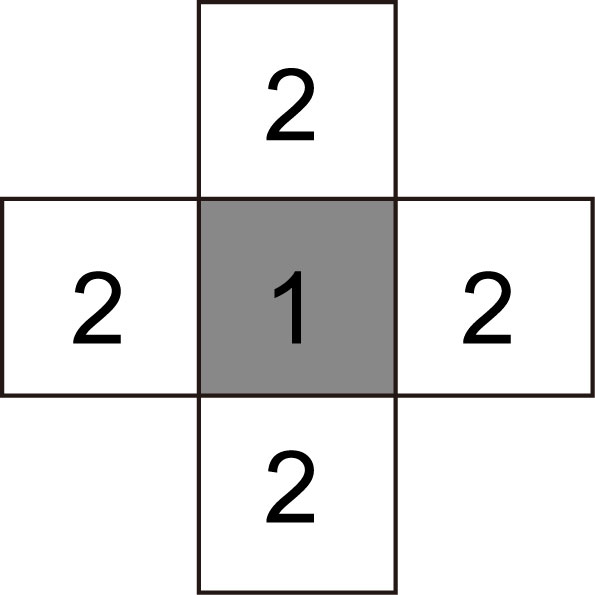} 
    \end{minipage}%
    \begin{minipage}{0.75\linewidth}
    \vspace{-1.4em}
    \[
    \left\{
    \begin{aligned}
    p_{1,i,j}^{k,m} &= \frac{p_{1,i+1,j}^{k,m-1} + p_{1,i-1,j}^{k,m-1} + p_{1,i,j+1}^{k,m-1} + p_{1,i,j-1}^{k,m-1} + h_1^2 f_{1,i,j}}{4 + h_1^2}, \\
    i+j &= \text{odd}.
    \end{aligned}
    \right.
    \]
    \end{minipage}
    \end{equation}

    \item Parallel update of set 2:
    \begin{equation}\label{BWGS-W}
    \begin{minipage}{0.15\linewidth}
    \includegraphics[width=\linewidth]{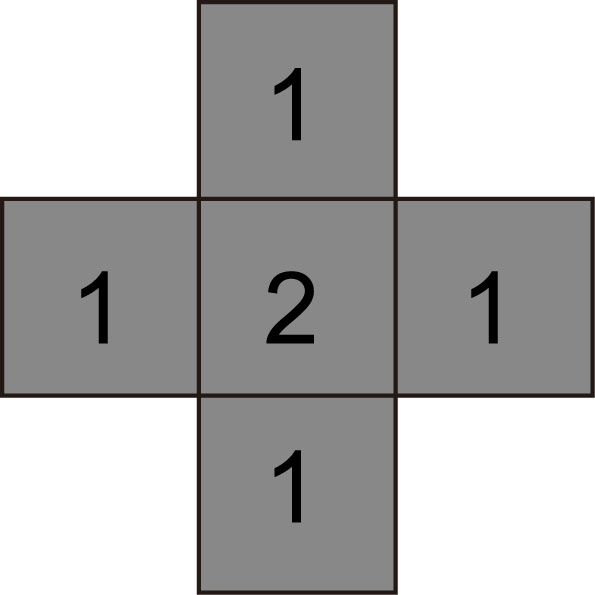} 
    \end{minipage}%
    \begin{minipage}{0.75\linewidth}
    \vspace{-1.4em}
    \[
    \left\{
    \begin{aligned}
    p_{1,i,j}^{k,m} &= \frac{p_{1,i+1,j}^{k,m} + p_{1,i-1,j}^{k,m} + p_{1,i,j+1}^{k,m} + p_{1,i,j-1}^{k,m} + h_1^2 f_{1,i,j}}{4 + h_1^2}, \\
    i+j &= \text{even}.
    \end{aligned}
    \right.
    \]
    \end{minipage}
    \end{equation}
\end{enumerate}

Updates for the two sets (i)–(ii) are executed sequentially, A MATLAB implementation of the RBGS method is provided in Listing~\ref{RBGS smoothing}:
\begin{lstlisting}[caption={MATLAB implementation of the smoothing operator $\mathrm{S}$ using the Red-Black Gauss–Seidel scheme in 2D}, label={RBGS smoothing}]
function p = Smoothing(f, p, M_1, N_1)
% Smoothing for Red-Black Gauss-Seidel 

for color = 1:2   % 1 = black, 2 = Red
    for i = 1:M_1
        for j = 1:N_1
            if mod(i+j,2) == color
                % update single points (i,j)
                p(i,j) = GaussSeidel(f, p, i, j);     
            end
        end
    end
end

end
\end{lstlisting}

\section{MATLAB implementation of the smoothing operator using the X-shape multi-color Gauss–Seidel scheme in 3D}
\begin{lstlisting}[caption={MATLAB implementation of the smoothing operator $\mathrm{S}$ using the X-shape Multi-Color Gauss–Seidel scheme in 3D}]
function p = Smoothing(f, p, M_1, N_1, L_1)
% Smoothing for X-shape Multi-Color Gauss-Seidel in 3D, M_1, N_1, L_1 represent the numbers of grid intervals in the horizontal, vertical, and depth directions on the fine grid. 

% Parallel update of set 1
vector_i=1:2:M_1-1; vector_j=2:2:N_1; vector_l=2:2:L_1;
  p(vector_i,vector_j,vector_l) 
= GaussSeidel(f, p, vector_i, vector_j, vector_l); 

% Parallel update of set 2
vector_i=2:2:M_1; vector_j=1:2:N_1-1; vector_l=2:2:L_1;
  p(vector_i,vector_j,vector_l) 
= GaussSeidel(f, p, vector_i, vector_j, vector_l); 

% Parallel update of set 3
vector_i=2:2:M_1-1; vector_j=2:2:N_1; vector_l=1:2:L_1-1;
  p(vector_i,vector_j,vector_l) 
= GaussSeidel(f, p, vector_i, vector_j, vector_l); 

% Parallel update of set 4
vector_i=1:2:M_1-1; vector_j=1:2:N_1-1; vector_l=1:2:L_1-1;
  p(vector_i,vector_j,vector_l) 
= GaussSeidel(f, p, vector_i, vector_j, vector_l); 

% Parallel update of set 5
vector_i=1:2:M_1-1; vector_j=2:2:N_1; vector_l=1:2:L_1-1;
  p(vector_i,vector_j,vector_l) 
= GaussSeidel(f, p, vector_i, vector_j, vector_l);  

% Parallel update of set 6
vector_i=2:2:M_1; vector_j=1:2:N_1-1; vector_l=1:2:L_1-1;
  p(vector_i,vector_j,vector_l) 
= GaussSeidel(f, p, vector_i, vector_j, vector_l); 

% Parallel update of set 7
vector_i=2:2:M_1-1; vector_j=2:2:N_1; vector_l=2:2:L_1;
  p(vector_i,vector_j,vector_l) 
= GaussSeidel(f, p, vector_i, vector_j, vector_l);  

% Parallel update of set 8
vector_i=1:2:M_1-1; vector_j=1:2:N_1-1; vector_l=2:2:L_1;
  p(vector_i,vector_j,vector_l) 
= GaussSeidel(f, p, vector_i, vector_j, vector_l); 
end
\end{lstlisting}
\label{app2}
\section{Restriction $\mathrm{R}$ and Prolongation $\mathrm{P}$ operators for cell-centered grids in 2D}\label{app3}
We first introduce the restriction $\mathrm{R}$ and prolongation $\mathrm{P}$ operators in Algorithm~\ref{two-grid algorithm} for cell-centered grids in 2D, which serve as the foundation for multigrid transfer between coarse and fine grids. These operators are applied only to interior points and are defined under the assumption of a uniform grid structure, where each fine grid patch is square and centered within its corresponding coarse grid patch. For simplicity, we adopt the standard four-point averaging $\mathrm{R}$ operator and the corresponding $\mathrm{P}$ operator:
\begin{align}
\begin{aligned}
&p_{0,i,j} =\mathrm{R}(p_{1,i,j})
=\frac{p_{1,2i-1,2j-1}+p_{1,2i-1,2j}+p_{1,2i,2j-1}+p_{1,2i,2j}}{4},\\
&p_{0,i,j}\in \mathcal{C}_{\Omega_0},\,\, p_{1,i,j}\in \mathcal{C}_{\Omega_1},
\end{aligned}
\end{align}
and
\begin{align} 
p_{1,i,j}= \mathrm{P}(p_{0,i,j})\quad\Longrightarrow \quad
\left.\begin{aligned}p_{1,2i-1,2j-1}=&p_{0,i,j},\\
p_{1,2i-1,2j}=&p_{0,i,j},\\
p_{1,2i,2j-1}=&p_{0,i,j},\\
p_{1,2i,2j}=&p_{0,i,j},
\end{aligned}\right\}\quad p_{0,i,j}\in \mathcal{C}_{\Omega_0},\,\, p_{1,i,j}\in \mathcal{C}_{\Omega_1}.
\end{align}
These standard $\mathrm{R}$ and $\mathrm{P}$ operators in the 2D are sufficient to achieve second-order accuracy. Their 3D counterparts are defined analogously, and our framework can also be naturally extended to higher-order interpolation operators. See \cite{dipietro2023,hemker1990} for details. Finally, we emphasize that restriction and prolongation are inherently matrix operations, making them particularly well-suited for GPU-based parallelization, thereby improving computational efficiency and reducing execution time. 

\section{Restriction and Prolongation operator for edge-centered grids in 2D}\label{app4}
For edge-centered variables, the definition of $\mathrm{R}$ and $\mathrm{P}$ operators in Algorithm~\ref{two-grid algorithm} is inherently more involved due to the different variable locations. As an example, we now present the standard $\mathrm{R}$ operator for east–west (EW) edge-centered variables in 2D:
\begin{align}
\begin{aligned}
u_{0,i+\frac{1}{2},j}=&\mathrm{R}(u_{1,i+\frac{1}{2},j})\\
=&\frac{u_{1,2i-\frac{1}{2},2j-1}+2u_{1,2i-\frac{1}{2},2j}+u_{1,2i-\frac{1}{2},2j+1}}{8}\\
+&\frac{u_{1,2i+\frac{1}{2},2j-1}+2u_{1,2i+\frac{1}{2},2j}+u_{1,2i+\frac{1}{2},2j+1}}{8},\\
 u_{0,i+\frac{1}{2},j}&\in \mathcal{E}^{EW}_{\Omega_0},\,\,u_{1,i+\frac{1}{2},j}\in \mathcal{E}^{EW}_{\Omega_1}.
\end{aligned}
\end{align}
The corresponding $\mathrm{P}$ operator for EW edge-centered variables in 2D is defined by injection:
\begin{align} 
&u_{1,i+\frac{1}{2},j}=\mathrm{P}(u_{0,i+\frac{1}{2},j})\nonumber\\
\Longrightarrow
&\left.\begin{aligned}
u_{1,2i+\frac{1}{2},2j-1}=&\frac{3u_{0,i+\frac{1}{2},j}+u_{0,i+\frac{1}{2},j-1}}{4},\\
u_{1,2i+\frac{1}{2},2j}=&\frac{3u_{0,i+\frac{1}{2},j}+u_{0,i+\frac{1}{2},j+1}}{4},\\
u_{1,2i+\frac{3}{2},2j-1}=&\frac{u_{1,2i+\frac{1}{2},2j-1}+u_{1,2i+\frac{5}{2},2j-1}}{2},\\
u_{1,2i+\frac{3}{2},2j}=&\frac{u_{1,2i+\frac{1}{2},2j}+u_{1,2i+\frac{5}{2},2j}}{2},
\end{aligned}\right\}u_{0,i+\frac{1}{2},j}\in \bar{\mathcal{E}}^{EW}_{\Omega_0},\,\,u_{1,i+\frac{1}{2},j}\in \bar{\mathcal{E}}^{EW}_{\Omega_1}.
\end{align}

These standard $\mathrm{R}$ and $\mathrm{P}$ operators in the 2D are sufficient to achieve second-order accuracy. Their 3D counterparts are defined analogously, and our framework can also be naturally extended to higher-order interpolation operators. See \cite{dipietro2023,hemker1990} for details. Finally, we emphasize that restriction and prolongation are inherently matrix operations, making them particularly well-suited for GPU-based parallelization, thereby improving computational efficiency and reducing execution time.


%%       \cite{<label>} ==> [1]

%%
%Example citation, See \cite{lamport94}.

 \bibliographystyle{elsarticle-num} 
 \bibliography{cas-refs}

\begin{thebibliography}{10}
\expandafter\ifx\csname url\endcsname\relax
  \def\url#1{\texttt{#1}}\fi
\expandafter\ifx\csname urlprefix\endcsname\relax\def\urlprefix{URL }\fi
\expandafter\ifx\csname href\endcsname\relax
  \def\href#1#2{#2} \def\path#1{#1}\fi

\bibitem{bakhvalov1966}
N.~S. Bakhvalov, On the convergence of a relaxation method with natural constraints on the elliptic operator, USSR Comput. Math. Math. Phys. 6~(5) (1966) 101--135.
\newblock \href {https://doi.org/10.1016/0041-5553(66)90118-2} {\path{doi:10.1016/0041-5553(66)90118-2}}.

\bibitem{fedorenko1962}
R.~P. Fedorenko, A relaxation method for solving elliptic difference equations, USSR Comput. Math. Math. Phys. 1~(4) (1962) 1092--1096.
\newblock \href {https://doi.org/10.1016/0041-5553(62)90031-9} {\path{doi:10.1016/0041-5553(62)90031-9}}.

\bibitem{fedorenko1964}
R.~P. Fedorenko, The speed of convergence of one iterative process, USSR Comput. Math. Math. Phys. 4~(3) (1964) 559--559.
\newblock \href {https://doi.org/10.1016/0041-5553(64)90253-8} {\path{doi:10.1016/0041-5553(64)90253-8}}.

\bibitem{brandt1977}
A.~Brandt, Multi-level adaptive solutions to boundary-value problems, Math. Comput. 31~(138) (1977) 333--390.
\newblock \href {https://doi.org/10.2307/2006422} {\path{doi:10.2307/2006422}}.

\bibitem{hackbusch1977}
U.~Trottenberg, C.~W. Oosterlee, A.~Schuller, Multigrid methods, Academic press, 1977.

\bibitem{hackbusch1984}
W.~Hackbusch, Multi-grid methods and applications, Springer, 1984.
\newblock \href {https://doi.org/https://doi.org/10.1007/978-3-662-02427-0} {\path{doi:https://doi.org/10.1007/978-3-662-02427-0}}.

\bibitem{fulton1986}
S.~R. Fulton, P.~E. Ciesielski, W.~H. Schubert, Multigrid methods for elliptic problems: A review, Mon. Weather Rev. 114~(5) (1986) 943--959.
\newblock \href {https://doi.org/https://doi.org/10.1175/1520-0493(1986)114<0943:MMFEPA>2.0.CO;2} {\path{doi:https://doi.org/10.1175/1520-0493(1986)114<0943:MMFEPA>2.0.CO;2}}.

\bibitem{trottenberg2001}
U.~Trottenberg, C.~W. Oosterlee, A.~Schüller, Multigrid, Elsevier, San Diego, CA, USA, 2000.

\bibitem{Bolz2003}
J.~Bolz, I.~Farmer, E.~Grinspun, P.~Schr{\"o}der, Sparse matrix solvers on the gpu: conjugate gradients and multigrid, ACM T. Graphics 22~(3) (2003) 917--924.
\newblock \href {https://doi.org/https://doi.org/10.1145/882262.882364} {\path{doi:https://doi.org/10.1145/882262.882364}}.

\bibitem{ljungkvist2017}
K.~Ljungkvist, M.~Kronbichler, Multigrid for matrix-free finite element computations on graphics processors, Tech. rep., Uppsala University, Department of Information Technology, Division of Scientific Computing, preprint or technical report (2017).
\newblock \href {https://doi.org/https://dl.acm.org/doi/10.5555/3108096.3108097} {\path{doi:https://dl.acm.org/doi/10.5555/3108096.3108097}}.

\bibitem{feng2014numerical}
C.~Feng, S.~Shu, J.~Xu, C.-S. Zhang, Numerical study of geometric multigrid methods on cpu-gpu heterogeneous computers, Adv. Appl. Math. and Mech. 6~(1) (2014) 1--23.
\newblock \href {https://doi.org/https://doi.org/10.1017/S2070073300002411} {\path{doi:https://doi.org/10.1017/S2070073300002411}}.

\bibitem{antepara2024}
O.~Antepara, S.~Williams, H.~Johansen, M.~Hall, High-performance, scalable geometric multigrid via fine-grain data blocking for gpus, in: Proceedings of the SC24 Workshops, Atlanta, GA, USA, 2024, pp. 1177--1191.
\newblock \href {https://doi.org/10.1109/SCW63240.2024.00159} {\path{doi:10.1109/SCW63240.2024.00159}}.

\bibitem{cui2025}
C.~Cui, G.~Kanschat, Multigrid methods for the stokes problem on gpu systems, Comput. Fluids 299 (2025) 106703.
\newblock \href {https://doi.org/https://doi.org/10.1016/j.compfluid.2025.106703} {\path{doi:https://doi.org/10.1016/j.compfluid.2025.106703}}.

\bibitem{chen2024}
A.~Chen, B.~A. Erickson, J.~E. Kozdon, J.~Choi, Matrix-free sbp-sat finite difference methods and the multigrid preconditioner on gpus, in: Proceedings of the 38th ACM International Conference on Supercomputing, 2024, pp. 400--412.
\newblock \href {https://doi.org/https://doi.org/10.1145/3650200.3656614} {\path{doi:https://doi.org/10.1145/3650200.3656614}}.

\bibitem{shi2020}
X.~Shi, T.~Agrawal, C.-A. Lin, F.-N. Hwang, T.-H. Chiu, A parallel nonlinear multigrid solver for unsteady incompressible flow simulation on multi-gpu cluster, J. Comput. Phys. 414 (2020) 109447.
\newblock \href {https://doi.org/https://doi.org/10.1016/j.jcp.2020.109447} {\path{doi:https://doi.org/10.1016/j.jcp.2020.109447}}.

\bibitem{BraedstrupEgholm2014}
C.~F. Br{\ae}dstrup, A.~Damsgaard, D.~L. Egholm, Ice-sheet modelling accelerated by graphics cards, Comput. Geosci. 72 (2014) 210--220.
\newblock \href {https://doi.org/https://doi.org/10.1016/j.cageo.2014.07.019} {\path{doi:https://doi.org/10.1016/j.cageo.2014.07.019}}.

\bibitem{Gorobets2024}
A.~V. Gorobets, S.~Soukov, A.~Magomedov, Heterogeneous parallel implementation of a multigrid method with full approximation in the noisette code, Math. Models Comput. Simul. 16~(4) (2024) 609--619.
\newblock \href {https://doi.org/10.1134/S2070048224700261} {\path{doi:10.1134/S2070048224700261}}.

\bibitem{Ansari2025}
M.~Q. Ansari, M.~Q. Ansari, Accelerating matrix multiplication: A performance comparison between multi-core cpu and gpu, arXiv preprint arXiv:2507.19723 (2025).
\newblock \href {https://doi.org/https://doi.org/10.48550/arXiv.2507.19723} {\path{doi:https://doi.org/10.48550/arXiv.2507.19723}}.

\bibitem{Recasens2025}
P.~G. Recasens, F.~Agullo, Y.~Zhu, C.~Wang, E.~K. Lee, O.~Tardieu, J.~Torres, J.~L. Berral, Mind the memory gap: Unveiling gpu bottlenecks in large-batch llm inference, arXiv preprint arXiv:2503.08311 (2025).
\newblock \href {https://doi.org/https://doi.org/10.48550/arXiv.2503.08311} {\path{doi:https://doi.org/10.48550/arXiv.2503.08311}}.

\bibitem{griebel2010}
M.~Griebel, P.~Zaspel, A multi-gpu accelerated solver for the three-dimensional two-phase incompressible navier-stokes equations, Comput. Sci. Res. Dev. 25~(1) (2010) 65--73.
\newblock \href {https://doi.org/https://doi.org/10.1007/s00450-010-0111-7} {\path{doi:https://doi.org/10.1007/s00450-010-0111-7}}.

\bibitem{Kashi2025}
A.~Kashi, S.~Nadarajah, On the effectiveness of fine-grain parallel linear iterations for computational aerodynamics on structured grids for graphics processing units, Comput. Fluids 299 (2025) 106714.
\newblock \href {https://doi.org/10.1016/j.compfluid.2025.106714} {\path{doi:10.1016/j.compfluid.2025.106714}}.

\bibitem{saad2003}
Y.~Saad, Iterative methods for sparse linear systems, SIAM, 2003.
\newblock \href {https://doi.org/https://doi.org/10.1137/1.9780898718003} {\path{doi:https://doi.org/10.1137/1.9780898718003}}.

\bibitem{mazumder2016}
S.~Mazumder, Numerical Methods for Partial Differential Equations: Finite Difference and Finite Volume Methods, Academic Press, Cambridge, MA, 2016.

\bibitem{wesseling1992}
P.~Wesseling, An Introduction to Multigrid Methods, John Wiley, Chichester, UK, 1992.

\bibitem{briggs2000}
W.~L. Briggs, V.~E. Henson, S.~F. McCormick, A multigrid tutorial, SIAM, 2000.
\newblock \href {https://doi.org/https://doi.org/10.1137/1.9780898719505} {\path{doi:https://doi.org/10.1137/1.9780898719505}}.

\bibitem{ZHAO2022}
L.~Zhao, C.~Feng, C.-S. Zhang, S.~Shu, Parallel multi-stage preconditioners with adaptive setup for the black oil model, Comput. Geosci. 168 (2022) 105230.
\newblock \href {https://doi.org/https://doi.org/10.1016/j.cageo.2022.105230} {\path{doi:https://doi.org/10.1016/j.cageo.2022.105230}}.

\bibitem{li2020}
J.~Li, J.~Liu, G.~D. Egbert, R.~Liu, R.~Guo, K.~Pan, An efficient preconditioner for 3-d finite difference modeling of the electromagnetic diffusion process in the frequency domain, IEEE Trans. Geosci. and Remote Sens. 58~(1) (2019) 500--509.
\newblock \href {https://doi.org/10.1109/TGRS.2019.2937742} {\path{doi:10.1109/TGRS.2019.2937742}}.

\bibitem{fan1997}
D.~Fan, L.-Q. Chen, Computer simulation of grain growth using a continuum field model, Acta Mater. 45~(2) (1997) 611--622.
\newblock \href {https://doi.org/10.1016/S1359-6454(96)00200-5} {\path{doi:10.1016/S1359-6454(96)00200-5}}.

\bibitem{wang2022lipid}
Y.~Wang, Y.~Palzhanov, A.~Quaini, M.~Olshanskii, S.~Majd, Lipid domain coarsening and fluidity in multicomponent lipid vesicles: A continuum based model and its experimental validation, Biochim. Biophys. Acta Biomembr. 1864~(7) (2022) 183898.
\newblock \href {https://doi.org/10.1016/j.bbamem.2022.183898} {\path{doi:10.1016/j.bbamem.2022.183898}}.

\bibitem{ratz_pdes_2006}
A.~Rätz, A.~Voigt, {PDE}'s on surfaces --- a diffuse interface approach, Commun. Math. Sci. 4~(3) (2006) 575--590.
\newblock \href {https://doi.org/10.4310/CMS.2006.v4.n3.a5} {\path{doi:10.4310/CMS.2006.v4.n3.a5}}.

\bibitem{teigen2009diffuse}
K.~E. Teigen, X.~Li, J.~Lowengrub, F.~Wang, A.~Voigt, A diffuse-interface approach for modeling transport, diffusion and adsorption/desorption of material quantities on a deformable interface, Commun. Math. Sci. 4~(7) (2009) 1009.
\newblock \href {https://doi.org/10.4310/cms.2009.v7.n4.a10} {\path{doi:10.4310/cms.2009.v7.n4.a10}}.

\bibitem{li2009solving}
X.~Li, J.~Lowengrub, A.~R{\"a}tz, A.~Voigt, Solving pdes in complex geometries: a diffuse domain approach, Commun. Math. Sci. 7~(1) (2009) 81.
\newblock \href {https://doi.org/10.4310/cms.2009.v7.n1.a4} {\path{doi:10.4310/cms.2009.v7.n1.a4}}.

\bibitem{guo2022}
Z.~Guo, Q.~Cheng, P.~Lin, C.~Liu, J.~Lowengrub, Second order approximation for a quasi-incompressible navier-stokes cahn-hilliard system of two-phase flows with variable density, J. Comput. Phys. 448 (2022) 110727.
\newblock \href {https://doi.org/https://doi.org/10.1016/j.jcp.2021.110727} {\path{doi:https://doi.org/10.1016/j.jcp.2021.110727}}.

\bibitem{denefle2014}
R.~Den{\`e}fle, S.~Mimouni, J.~Caltagirone, S.~Vincent, Multifield hybrid method applied to bubble rising and coalescence, Int. J. Comput. Methods Exp. Meas. 2~(1) (2014) 46--57.
\newblock \href {https://doi.org/10.2495/CMEM-V2-N1-46-57} {\path{doi:10.2495/CMEM-V2-N1-46-57}}.

\bibitem{brereton1991}
G.~Brereton, D.~Korotney, Coaxial and oblique coalescence of two rising bubbles, in: Dynamics of Bubbles and Vortices Near a Free Surface, Vol. 119, ASME, 1991, pp. 50--73.

\bibitem{FENG2018}
W.~Feng, Z.~Guo, J.~S. Lowengrub, S.~M. Wise, A mass-conservative adaptive fas multigrid solver for cell-centered finite difference methods on block-structured, locally-cartesian grids, Journal of Computational Physics 352 (2018) 463--497.
\newblock \href {https://doi.org/https://doi.org/10.1016/j.jcp.2017.09.065} {\path{doi:https://doi.org/10.1016/j.jcp.2017.09.065}}.

\bibitem{alvarez2022}
G.~B. Alvarez, D.~C. Lobao, W.~A. de~Menezes, The m-order jacobi, gauss–seidel and symmetric gauss–seidel methods, Pesqui. Ensino Ci{\^e}nc. Exatas Nat. 6~(1) (2022) 10.
\newblock \href {https://doi.org/https://doi.org/10.29215/pecen.v6i0.1773} {\path{doi:https://doi.org/10.29215/pecen.v6i0.1773}}.

\bibitem{bhatti2023}
N.~Bhatti, Niketa, Comparative study of symmetric gauss--seidel methods and preconditioned symmetric gauss--seidel methods for linear system, Int. J. Sci. Res. Arch. 08~(01) (2023) 940--947.
\newblock \href {https://doi.org/10.30574/ijsra.2023.8.1.0155} {\path{doi:10.30574/ijsra.2023.8.1.0155}}.

\bibitem{li2022}
X.~Li, J.~Shen, Z.~Liu, New sav-pressure correction methods for the navier-stokes equations: stability and error analysis, Math. Comput. 91~(333) (2022) 141--167.
\newblock \href {https://doi.org/https://doi.org/10.1090/mcom/3651} {\path{doi:https://doi.org/10.1090/mcom/3651}}.

\bibitem{liu1994}
X.-D. Liu, S.~Osher, T.~Chan, Weighted essentially non-oscillatory schemes, J. Comput. Phys. 115~(1) (1994) 200--212.
\newblock \href {https://doi.org/https://doi.org/10.1006/jcph.1994.1187} {\path{doi:https://doi.org/10.1006/jcph.1994.1187}}.

\bibitem{jiang1996}
G.-S. Jiang, C.-W. Shu, Efficient implementation of weighted eno schemes, J. Comput. Phys. 126~(1) (1996) 202--228.
\newblock \href {https://doi.org/https://doi.org/10.1006/jcph.1996.0130} {\path{doi:https://doi.org/10.1006/jcph.1996.0130}}.

\bibitem{shu1999}
C.-W. Shu, High order eno and weno schemes for computational fluid dynamics, in: C.~Fetecau, C.~Groth, R.~Mittal (Eds.), High-order Methods for Computational Physics, Springer, Berlin, 1999, pp. 439--582.
\newblock \href {https://doi.org/https://doi.org/10.1007/978-3-662-03882-6\_5} {\path{doi:https://doi.org/10.1007/978-3-662-03882-6\_5}}.

\bibitem{ghia1982}
U.~Ghia, K.~N. Ghia, C.~Shin, High-re solutions for incompressible flow using the navier-stokes equations and a multigrid method, J. Comput. Phys. 48~(3) (1982) 387--411.
\newblock \href {https://doi.org/10.1016/0021-9991(82)90058-4} {\path{doi:10.1016/0021-9991(82)90058-4}}.

\bibitem{ku1987}
H.~C. Ku, R.~S. Hirsh, T.~D. Taylor, A pseudospectral method for solution of the three-dimensional incompressible navier-stokes equations, J. Comput. Phys. 70~(2) (1987) 439--462.
\newblock \href {https://doi.org/10.1016/0021-9991(87)90190-2} {\path{doi:10.1016/0021-9991(87)90190-2}}.

\bibitem{dipietro2023}
D.~A. Di~Pietro, P.~Matalon, P.~Mycek, U.~R{\"u}de, High-order multigrid strategies for hybrid high-order discretizations of elliptic equations, Numer. Linear Algebra Appl. 30~(1) (2023) e2456.
\newblock \href {https://doi.org/https://doi.org/10.1002/nla.2456} {\path{doi:https://doi.org/10.1002/nla.2456}}.

\bibitem{hemker1990}
P.~Hemker, On the order of prolongations and restrictions in multigrid procedures, J. Comput. Appl. Math. 32~(3) (1990) 423--429.
\newblock \href {https://doi.org/https://doi.org/10.1016/0377-0427(90)90047-4} {\path{doi:https://doi.org/10.1016/0377-0427(90)90047-4}}.

\end{thebibliography}



%%%%%%%%%%%%%%%%%%%%%%%
% \cite{}

% \bibliographystyle{unsrt}
% \bibliography{cas-refs}
%%%%%%%%%%%%%%%%%%%%%%%%%5

%% For numbered reference style

% \bibitem{wang2022lipid}
% Y. Wang, Y. Palzhanov, A. Quaini, M. Olshanskii, S. Majd,
% \textit{Lipid domain coarsening and fluidity in multicomponent lipid vesicles: A continuum based model and its experimental validation},
% Biochim. Biophys. Acta Biomembr., \textbf{1864} (7), 2022, 183898.

\end{document}